\documentclass[11pt,twoside]{article}
\usepackage{amsmath,amssymb}
\usepackage{graphics}
\usepackage[usenames]{color}
\baselineskip 3.1pt \lineskip 3.1pt \numberwithin{equation}{section}
\def\AAu{{{\cal A}_u}}
\def\BBu{{{\cal B}_{uv}}}

\def\tr{{\rm tr}}
\def\QJ{{quasi-Jacobi}}
\def\QJ{{Jacobi }}

\def\ign{_{\rm igno}}
\def\de{e}
\def\Dhat{\ddot}
\def\Plar{\raisebox{0pt}{\mbox{${\sc\prod\limits_{\leftharpoondown}{\!\ssc\,}}$}}}
\def\Prar{\raisebox{0pt}{\mbox{${\sc\prod\limits_{\rightharpoondown}{\!\ssc\,}}$}}}
\def\pri#1{{\textbf{\textit{#1}}}}

\def\UPo#1{_{\la#1\ra}}
\def\UP#1{_{\la-#1\ra}}
\def\Comp#1{_{\la#1\ra}}
\def\aaa{\sigma}
\def\gcd{{\rm gcd}}
\def\Coeff{{\rm C_{oeff}}}
\def\comp#1{{}_{[#1]}{}}

\def\ptl{\partial}
\def\Supp{{\rm Supp}}
\def\res{{\rm Res}}
\def\deg{{\rm deg}}
\def\max{{\rm max}}

\def\ad{{\ssc\,}{\rm ad}}
\def\DD{\hbox{$I\hskip -4pt D$}}
\def \N{\hbox{$I\hskip -3pt N$}}
\def \Z{\hbox{$Z\hskip -5.2pt Z$}}

\def \Q{\hbox{$Q\hskip -5pt \vrule height 6pt depth 0pt\hskip 6pt$}}
\def \R{\mathbb{R}}
\def \C{\hbox{$C\hskip -5pt \vrule height 6pt depth 0pt \hskip 6pt$}}

\def\qed{\ \ \ifhmode\unskip\nobreak\fi\ifmmode\ifinner
         \else\hskip5pt\fi\fi
 \hbox{\hskip5pt\vrule width4pt height6pt depth1.5pt\hskip 1 pt}}
\def\a{\alpha}
\def\b{\beta}
\def\d{\delta}

\def\g{\gamma}

\def\l{\lambda}

\def\si{\sigma}
\def\sc{\scriptstyle}
\def\ssc{\scriptscriptstyle}
\def\dis{\displaystyle}
\def\cl{\centerline}
\def\DD{{\cal D}}

\def\nl{\newline}
\def\ol{\overline}

\def\Res{{\rm Res}}

\def\Lra{\Longleftrightarrow}
\def\bs{\backslash}

\def\rb{\raisebox}
\def\vs{\vspace*}
\def\VS#1{}
\def\PAR{}
\def\NOindent{}
\def\ra{\rangle}
\def\la{\langle}
\def\PP{{\cal P}}
\def\ni{\noindent}

\def\barx{x}
\def\bary{y}

\def\WW{{\cal W}}
\def\AA{{\cal A}}
\def\BB{{\cal B}}
\def\CC{{\cal C}}

\def\C{\mathbb{C}}
\def\F{\mathbb{C}}
\def\Z{\mathbb{Z}}
\def\N{\mathbb{N}}
\def\Q{\mathbb{Q}}

\newtheorem{theo}{Theorem}[section]
\newtheorem{lemm}[theo]{Lemma}

\newtheorem{rema}[theo]{Remark}
\newtheorem{exam}[theo]{Example}

\newtheorem{coro}[theo]{Corollary}

\newtheorem{defi}[theo]{Definition}
\def\equa#1#2{\begin{equation}\label{#1}\mbox{$#2$}\end{equation}}
\def\equan#1#2{\begin{equation*}\mbox{$#2$}\end{equation*}}

\begin{document}
\def\OUTLINE{
\setcounter{page}{1}\pagenumbering{roman}
\noindent\underline{\bf Outline of Approach:}\vskip4pt
\par
\def\NP{{\rm NP}}
The approach is {\color{red}purely algebraic}, which  is mainly divided into 2 steps:\vskip4pt
\par
\noindent{\bf Step 1} (Lemma \ref{Nerw}):
Let $(F,G)$ be a {\it Jacobi pair}, i.e., $F,G\in\C[x,y]$ such that \equa{Jacobi-con}{J(F,G):=(\ptl_xF)(\ptl_yG)-(\ptl_yF)(\ptl_xG)\in\C\bs\{0\}.}
It is a well-known fact that the 2-dimensional Jacobian conjecture is equivalent to
the following statement:
\par
``For every Jacobi pair $(F,G)$ with $\deg_yF\ge2$, the Newton polygon $\NP(F)$ of $F$ is a triangle.''
Here, we define
\begin{eqnarray*}&&
\NP(F)=\mbox{ the minimal polygon which surrounds }(\Supp\,F)\cup\{0\},\\
&&
\Supp\,F=\{(i,j)\,|\,f_{ij}\ne0\}\mbox{ if $F=\sum\limits_{i,j}$}f_{ij}x^iy^j,\ f_{ij}\in\C.
\end{eqnarray*}
Thus we can start from a Jacobi pair $(F,G)\in\C[x,y]\times\C[x,y]$ (later on, after we apply some variable changes, we shall regard $F,G$ as in $\C[x^{\pm\frac1N},y]$ for a fixed sufficient large $N$ such that all elements under consideration are in $\C[x^{\pm\frac1N},y]$), such that $\Supp\,F$ (resp., $\Supp\,G$) has a vertex
$(m_0,m)\in\Q_+\times\Z_+$ (resp., $(n_0,n)$) and (cf.~\eqref{Pick0.1})
$$
0<m_0<m,\ \ \frac{m_0}{m}=\frac{n_0}{n},\ \ m,n\ge2,\ \ m\not|n,\ \ n\not|m.
$$
Let $L$ (resp., $L'$) be the side of $\Supp\,F$ (resp., $\Supp\,G$) with top vertex
$(m_0,m)$ (resp., $(n_0,n)$). We can suppose $L$ and $L'$ have the same slope $\frac{q}{p}>\frac{m}{m_0}>1$ (cf.~\S\ref{S3.4}), where $p$ can be zero. Let $F_L$ (resp., $G_L$) be the part of $F$ (resp., $G$) whose support is $L$ (resp., $L'$).
\equa{Pick0.1}{\ \hspace*{-150pt}
\put(0,50){$\line(-1,-4){20}
\put(17,0){$\sc\circ\ (m_0,m):\, 0<m_0<m,\,m_0\in\Q_+,\,m\in\Z_+$}
\put(20,-35){$\sc L:\mbox{ slope }\frac{q}{p},\ \ \  F_L:\ \Supp\,F_L=L$}
\put(-6,-85){$\sc\circ  $}
$}}
We always have  (cf.~\eqref{p0-q0})
\begin{eqnarray}\label{1-p0-q0}&&G_L=(F_L)^{\frac{n}{m}} \mbox{ and } \frac{q}{p}>\frac{m+n-1}{m_0+n_0-1}, \mbox{ or else,}\\
\label{2-p0-q0}&&(F_L,G_L)\mbox{ is a \QJ  pair and } \frac{q}{p}=\frac{m+n-1}{m_0+n_0-1}.
\end{eqnarray}
Here we say a pair $(F,G)$ is a {\it \QJ  pair} if $J(F,G)\in\C\bs\{0\}$ but $F,G$ are not necessarily in $\C[x,y]$.
\par
{\color{blue}Note that in case $F,G\in\C[x,y]$, we always have \eqref{1-p0-q0}.
If we only consider a side of $\NP(F)$ such that \eqref{1-p0-q0} holds, we find out that it is very difficult to effectively make use of the condition of the nonzero constant Jacobian determinant \eqref{Jacobi-con}. This suggests us to apply some variable changes (rather than degree reductions) so that we can find a {\color{red}very special} side of $\NP(F)$ such that \eqref{2-p0-q0} holds.}
\par
The {\color{red}most important feature} of Lemma \ref{Nerw} shows that by a variable change of the form
(cf.~\eqref{VAR}) \equan{MAMSM}{y\mapsto
y+\l_2x^{-\frac{p_1}{q_1}}+\cdots+\l_kx^{-\frac{p_k}{q_k}},}
for some $\l_i\in\C$, $p_i,q_i\in\N$ with $\frac{p_i}{q_i}<\frac{p_1}{q_1}<1$ for $i\ge1$, $F,G$ become some elements in $\C[x^{\pm\frac1N},y]$ such that $\NP(F)$ has the {\color{red}special side  $L$} (which is in fact the lowest side)
{\color{red}satisfying \eqref{2-p0-q0}}:
\equa{1Pick0.1}{\color{blue}\ \hspace*{-150pt}
\put(1,5){$\line(1,5){10}$}\put(11,55){$\sc\circ\ (m^0_0,m^0)$}
\put(12,25){$\sc L^0:\mbox{ slope }\frac{q_0}{p_0},\ \ G_{L^0}=(F_{L^0})^{\frac{n}{m}}$}
{\color{red}\line(-1,-3){20}}
\put(17,0){\color{red}$\sc\circ\ (m_0,m):\, 0<m_0<m,\,m_0\in\Q_+,\,m\in\Z_+$}
\put(20,-35){\color{red}$\sc L:\mbox{ slope }\frac{q}{p},\ \ \  F_L:\ \Supp\,F_L=L$}
\put(20,-50){\color{red}$\sc (F_L,G_L)\mbox{ is a \QJ  pair}$}
\put(-6,-65){\color{red}$\sc\circ $}
}
{\color{blue}By examining many examples of $F,G\in\C[x^{\pm\frac1N},y]$, we find out that in order to rule out some possibilities, it is not enough to only consider one side of $\NP(F)$. Thus we shall not only need to consider the {\color{red}lowest side $L$}
but also the {\color{blue}side $L^0$} which has the common vertex $(m_0,m)$ with $L$.}
\par
By changing variables
$$
(x,y)\mapsto(x^{\frac{q}{q-p}},x^{-\frac{p}{q-p}}y),
$$
the $\NP(F)$ becomes the one which looks as follows
(cf.~\eqref{new-supp-FCOM-region}, {\color{blue}\eqref{f0===}}):
\equa{1-new-supp-FCOM-region}{%
\mbox{$\!\!\!\!\!\!\!\!\!\!\!\!\!\!\!\!\!\!\!\!\!\!\!\!\!\!\!\!$} \put(-5,3){$\vector(1,0){130}$}\put(72,-5){$\sc0$}\put(125,-5){$\sc
 x$}
 \put(76,-3){$\vector(0,1){110}$}\put(70,108){$\sc y$}
\put(5,70){\put(0,-67){\line(-1,0){150}}
\put(-22,23){$\sc\circ$}
\put(-48,0){\line(-1,-1){48}}\put(-48,0){\line(1,1){25}}
\put(-100,-49){$\sc\circ$}\put(-96,-47){\line(3,-2){30}}\put(-67,-69){$\sc\circ$}
$\put(-20,25){$\line(3,-1){79}$}$}\put(65,70){\color{blue}$\line(1,-1){35}$}
\put(64,68){$\sc\circ$}\put(58,75){\color{blue}$\sc\,\tau^0$}
\put(99,32){$\sc\circ$}\put(104,30){\color{red}
$\sc\,\tau =(\frac{m }{m +n },m )$}
\put(102,16){\color{red}\
$\sc\,L:\ \mbox{\color{red}write }\color{red}\ F_L=x^{\frac{m}{m+n}}f\mbox{ for some }f\in\C[y]$}
\put(100,35){\color{red}
$\line(0,-1){32}$}
%
%
%
%
\put(15,40){$\sc\Supp\,F$} \put(82,56){\color{blue}$\sc L^{0}:\
\mbox{slope $-\frac1{\a},$ write }
F_{L^0}=x^{\frac{m}{m+n}}y^m\prod\limits_{i=1}^d(a_ix^{-\a}y+1)^{m_i}$}
 \ \hspace*{80pt}\
}\VS{-19pt}\PAR\NOindent
and $\NP(G)$ is similar to $\NP(F)$ with the data $(m,n)$ replaced by $(n,m)$.
\par
\vskip4pt
\par
 \noindent{\bf Step 2} (Section 4): The space $\PP:=\C[y]((x^{-\frac1N}))$ has a natural structure of Poisson algebra.
 Regard $F,G$ as in $\PP$.
\par
{\color{blue}The main strategy is to give $\PP$ a different filtrations and consider automorphisms of the Poisson algebra $\PP$ which respect the filtrations and apply these automorphisms to $F,G$ so that we can observe
different properties of $F,G$ ({\color{red}this amounts to view things from different angles}). In this way, we can find some properties of $F,G$ which are inconsistent.}
\par
Lemma \ref{MMM-lemm} says that by writing $F,G$ (as in \eqref{1-new-supp-FCOM-region}) as
combinations of rational powers of $x$ with coefficients
in $\C[y]$ ({\color{red}this amounts to give a filtration
 of $\PP$ by $x$-degrees of its elements}):
\equa{F-G=i}{F=x^{\frac{m}{m+n}}(f+\sum\limits_{i=1}^{M_1}x^{-\frac{i}{N}}f_i),\ G=x^{\frac{n}{m+n}}(g+\sum\limits_{i=1}^{M_2}x^{-\frac{i}{N}}g_i),\ f_i,g_i\in\C[y], M_1,M_2\in\N,}
and using nonzero constant Jacobian determinant condition
 \eqref{Jacobi-con}, we see that $F,G$ can be obtained from $F_L,G_L$ by an automorphism.
More precisely, let $i_0>0$ be smallest such that $(f_{i_0},g_{i_0})\ne(0,0)$, then there
exists $Q_{i_0}\in x^{1-\frac{i_0}{N}}\C[y]$ such that ({\color{blue}with an exception, see \eqref{tau-c}}) \equa{q-f-g--}{
x^{\frac{m}{m+n}-\frac{i_0}{N}}f_{i_0}=[Q_{i_0},x^{\frac{m}{m+n}}f],\ \ \ \ \ \ \
x^{\frac{n}{m+n}-\frac{i_0}{N}}g_{i_0}=[Q_{i_0},x^{\frac{n}{m+n}}g].}
Thus if we apply the automorphism $e^{\ad_{-Q_{i_0}}}$ to $F,G$, then $f_{i_0},\,g_{i_0}$
become zero. This means that
there exists an automorphism $\si$ of the Poisson algebra $\PP$, which is a (possibly infinite but well-defined) product of the automorphisms of forms  $e^{\ad_{Q_i}}$ for $Q_i\in x^{1-\frac{i}N}\C[y]$ for $i\ge1$ and $\tau_c$ defined in \eqref{tau-c},  such that \equa{FAGA}{F=\si(x^{\frac{m}{m+n}}f)=H^{\frac{m}{m+n}}f(K),\
G=\si(x^{\frac{n}{m+n}}g)=H^{\frac{n}{m+n}}g(K),\mbox{ where }H=\si(x),K=\si(y).}
Note that although $F,G$ are in $\C[x^{\pm\frac1N},y]$, {\color{blue}$H,K$ are} in general not in $\C[x^{\pm\frac1N},y]$ but {\color{blue}in $\PP$.}
Furthermore, $\Supp\,H,\Supp\,K$ are below some line with negative slope. More precisely,
\par
\equa{?new-supp-FCOM-region}{%
\put(185,100){\color{blue}$\a'\le\a$}
\put(40,-30){\color{blue}$J(H,K)=J(H\UP0,K\UP0)=1,\ \ \#\Supp\,H\UP0\ge2,\ \#\Supp\,K\UP0\ge2$}
\put(80,40){\color{red}$\sc L_h: \mbox{ slope }-\frac1{\a'},\ H\UP0:=H_{L_h}$}
\put(95,0){\color{red}$\line(-1,2){55}$}\put(95,0){$\sc\circ$}\put(92,-10){$\sc(1,0)$}
\put(40,15){$\sc\Supp\,H$}
 \put(-5,3){$\vector(1,0){130}$}\put(72,-5){$\sc0$}\put(125,-5){$\sc
 x$}
 \put(76,-3){$\vector(0,1){110}$}\put(70,108){$\sc y$}
 \ \hspace*{220pt}\
 \put(-5,3){$\vector(1,0){130}$}\put(72,-5){$\sc0$}\put(125,-5){$\sc
 x$}
 \put(76,-3){$\vector(0,1){110}$}\put(70,108){$\sc y$}
\put(78,20){\color{red}$\line(-1,2){45}$}\put(78,20){$\sc\circ\ (0,1)$}
\put(25,25){$\sc\Supp\,K$}
\put(63,60){\color{red}$\sc L_k: \mbox{ slope }-\frac1{\a'},\ K\UP0:=K_{L_k}$}
 \ \hspace*{150pt}\
}
For $P=F,G,H$ or $K$, we can always write $P$ as $P=\sum_{i=0}^\infty P\UP{i}$, where $P\UP{i}$ (called the {\it $i$-th component} of $P$) is the part of $P$ whose support $\Supp\,P\UP{i}$ is on a line with slope $-\frac1{\a'}$
({\color{red}this amounts to give a filtration of $\PP$ by lines with slope $-\frac1{\a'}$}). Then
\equa{F---0000}{F\UP0=F_{L^0}\mbox{ if $\a'=\a,$ \ \ or $F\UP0=x^{\frac{m}{m+n}}y^m$ if $\a'<\a$ (cf.~\eqref{1-new-supp-FCOM-region} and \eqref{?new-supp-FCOM-region})}.}
Let $\nu >0$ be smallest such that $(H\UP{\nu },K\UP{\nu })\ne(0,0)$, then as in \eqref{q-f-g--}, there exists $Q_{\nu }\in\PP$ such that
\equa{H-nu=0-k}{H\UP{\nu }=[Q_{\nu },H\UP0],\ \ K\UP{\nu }=[Q_{\nu },K\UP0].}
We observe an {\color{red}important fact} that although $H,K$ are in general {\color{blue}not rational functions} on $y$ (cf.~\eqref{H-K-UP0=}), the following elements
\equa{111-ptl-x-K}{\dis H\ptl_xK=\frac{1}{(m+n)J}(mF\ptl_xG-nG\ptl_xF),\ \ \ H\ptl_yK=\frac{1}{(m+n)J}(mF\ptl_yG-nG\ptl_yF),}
are {\color{blue}polynomials} on $y$. Using this and the fact that $H^{\frac{m}{m+n}}f(K)=F$ is a polynomial on $y$, we obtain in Lemma \ref{Rational} that \equa{Q-nuuuuu}{\mbox{\color{red}$Q_\nu$ is a rational function on $y$}} ({\color{blue}except the special case $\nu=\nu_0$, which needs to be treated separately, where $\nu_0=1+\frac1{\a'}$}).
\par
In order to be able to use our machinery, we allow $F,G$ to be  rational functions on $y$. Then we can apply the automorphism $e^{\ad_{-Q_\nu}}$ to $F,G$ such that $Q_\nu$ becomes zero.
{\color{blue}First assume that the special case $\nu=\nu_0$ never occurs.} Then the above shows that we can suppose $H=H\UP0,$ $K=K\UP0$. In particular,
$$F=H\UP0^{\frac{m}{m+n}}f(K\UP0),\ \ \ G=H\UP0^{\frac{n}{m+n}}g(K\UP0).$$
{}From this and the fact that $f,g$ are coprime, by computing either the $m$-th component of $F$ or $n$-th component of $G$, we obtain \equan{AAMMM}{\mbox{\color{red}either $H\UP0^{\frac{m}{m+n}}$ or  $K\UP0^{\frac{n}{m+n}}$ is rational.}} However by \eqref{H-K-UP0=}, we see this is impossible. Thus we eventually obtain a contradiction
({\color{blue}the case $\nu=\nu_0$ is treated in details in the paper}).
\par
{\color{red}Remark: We can construct lots of examples of the \QJ  pairs $(F,G)\in\C[x^{\pm\frac1N},y]\times \C[x^{\pm\frac1N},y]$ with $\NP(F)$ as in \eqref{1-new-supp-FCOM-region} but the slope of $L^0$ is non-negative. Therefore the condition that the slope of $L^0$ is negative is necessary.}\vskip7pt
\par
\par
\ \par\ \par
Remark: Since Section 4 contains some computations, the following 3 pages are 3 {\it Mathematica} programs to verify
Equation \eqref{AMAMAMA}, and $r_1\ne0$ in \eqref{r-1=====0}, and $\tilde r_i,\,\tilde r\ne0$ in \eqref{t-R-7=====}.
\par
{\parskip0pt
\newpage
{\noindent\sl Pa20110808.nb}\hfill 1\ \ {}\\[-7pt] $\line(1,0){460}$\\
\par
{\footnotesize\tt\noindent\
\"{}The\ following\ is\ to\ prove\ (4.36),\ where\ H0,K0,Qnu,ap are 
H\_{\rm\{}$\la0\ra${\rm\}},K\_{\rm\{}$\la0\ra${\rm\}},Q\_nu,alpha'\"{}\par\noindent  \
Clear[H0]\par\noindent  \
Clear[K0]\par\noindent  \
Clear[Qnu]\par\noindent  \
Clear[ap]\par\noindent  \
\"{}The following 2 definitions are from (4.35)\"{}\par\noindent  \   
Hnu[x\_] = 1/(ap y)(H0[x] Qnu'[x] - (1 + ap(1 - nu))Qnu[x] H0'[x])\par\noindent  \
Knu[x\_] = 1/(ap y)(ap K0[x] Qnu'[x] - (1 + ap(1 - nu))Qnu[x] K0'[x])\par\noindent  \
\"{}The following is (4.33)\"{}\par\noindent  \   
R1[x\_] = (1 - nu)H0[x] Knu[x] + Hnu[x]K0[x]\par\noindent  \
\"{}Use (4.30) to verify whether the following is zero\"{}\par\noindent  \
f1 := ap R1'[x] - (1 + ap(1 - nu))(H0[x] Knu'[x] + Hnu[x] K0'[x])\par\noindent  \
\ \hspace*{10pt} \ \ \ /. K0'[x] {\tt->} (ap K0[x]H0'[x]+ap y)/H0[x]\par\noindent  \
\ \hspace*{10pt} \ \ \ /. K0''[x] {\tt->} ap(K0'[x]H0'[x]+K0[x]H0''[x])/H0[x]-
(ap K0[x]H0'[x]+ap y)H0'[x]/H0[x]\^{}2\par\noindent  \
\ \hspace*{10pt} \ \ \ /. K0'[x] {\tt->} (ap K0[x]H0'[x]+ap y)/H0[x]\par\noindent  \
Factor[f1]\par\noindent  \
Clear[H0]\par\noindent  \
Clear[K0]\par\noindent  \
Clear[Qnu]\par\noindent  \
Clear[ap]\par\noindent\
}
\par{\small
\noindent\ The following is to prove (4.36), where H0,K0,Qnu,ap are H\_\{$\la0\ra$\},K\_\{$\la0\ra$\},Q\_nu,alpha'\par\noindent  \
The following 2 definitions are from (4.35)\par\noindent  \
$\dis\rm\frac{-(1 + ap\,(1 - nu))\,Qnu[x]\ H0'[x]+H0[x]\ Qnu'[x]}{ap\ y}$\par\noindent  \
$\dis\rm\frac{-(1 + ap\,(1 - nu))\,Qnu[x]\ K0'[x]+ap\ K0[x]\ Qnu'[x]}{ap\ y}$\par\noindent  \
The following is (4.33)\par\noindent  \
$\dis\rm
\frac{K0[x]\, (-(1 + ap\, (1 - nu))\,Qnu[x]\ H0'[x]+H0[x]\ Qnu'[x])}{ap\ y}+$\par\noindent  \
\ \hspace*{6pt}$\dis\rm\ \ \frac{(1-nu)\,H0[x]\,(-(1+ap\,(1-nu))\,Qnu[x]\ K0'[x]+ap\ K0[x]\ Qnu'[x])}{ap\ y}$\par\noindent  \
Use (4.30) to verify whether the following is zero\par\noindent  \
0
}
\par
\newpage
\par
{\noindent\sl Pa20110808.nb}\hfill 2\ \ {}\\[-7pt] $\line(1,0){460}$\par\noindent
\par
{\footnotesize\tt\noindent\
\"{}The following program is to prove (4.44), i.e., r\_1 not= 0\"{}\par\noindent  \
Clear[H0]\par\noindent  \
Clear[K0]\par\noindent  \
Clear[qu]\par\noindent  \
Clear[ap]\par\noindent  \
Clear[bet]\par\noindent  \
\"{}The following is (4.40)\"{}\par\noindent  \
Qnu[x\_] = qn[x] + bet H0[x] K0[x]{\rm\^{}}(1 - nu)\par\noindent  \
Remark := \"{}The following 2 definitions are from (4.35)\"{}\par\noindent  \
Hnu[x\_] := 1/(ap y)(H0[x] Qnu'[x] - (1 + ap(1 - nu))Qnu[x] H0'[x])\par\noindent  \
Knu[x\_] := 1/(ap y)(ap K0[x] Qnu'[x] - (1 + ap(1 - nu))Qnu[x] K0'[x])\par\noindent  \
Remark := \"{}The following is (4.43)\"{}\par\noindent  \
R7[x\_] := (1\,-\,nu)\,ap\,(Knu[x]\,Hnu'[x]\,+\,Hnu[x]\,Knu'[x])\,-\,(1\,+\,ap(1\,-\,2nu))Hnu[x]\,Knu'[x]\par\noindent  \
r0 = Factor[Coefficient[R7[x], bet{\rm\^{}}2]bet{\rm\^{}}2/K0[x]{\rm\^{}}(-2\,nu)];\par\noindent  \
r1 = Factor[Coefficient[R7[x], bet]bet/K0[x]{\rm\^{}}(-nu)];\par\noindent  \
r2 = Factor[Coefficient[R7[x] bet, bet]];\par\noindent  \
qn[x\_] := -bet x y{\rm\^{}}(1 - nu)\par\noindent  \
H0[x\_] := x\par\noindent  \
K0[x\_] := y\par\noindent  \
\"{}Note: The lowest terms of r\_0, r\_1, r\_2 are resp.~y, y{\rm\^{}}(1-nu), y{\rm\^{}}(1-2\,nu)\"{}\par\noindent  \
\"{}\phantom{Note: }The following shows that r\_0, r\_1, r\_2 are nonzero\"{}\par\noindent  \
Factor[r0]\par\noindent  \
Factor[r1]\par\noindent  \
Factor[r2]\par\noindent  \
Clear[qn]\par\noindent  \
Clear[H0]\par\noindent  \
Clear[K0]\par\noindent  \
Clear[ap]\par\noindent  \
Clear[bet]\par\noindent \
}
\par{\small
\noindent\  The following program is to prove (4.44), i.e., r\_1 not= 0\par\noindent  \
The following is (4.40)\par\noindent  \
$\rm bet\ H0[x]\ K0[x]^{1 - nu} + qn[x]$\par\noindent  \
Note: The lowest terms of r\_0, r\_1, r\_2 are resp.~y, y{\rm\^{}}(1-nu), y{\rm\^{}}(1-2\,nu)\par\noindent  \
\phantom{Note: }The following shows that r\_0, r\_1, r\_2 are nonzero\par\noindent  \
$\rm -ap\ bet^2\,(-1 + nu)^2\,y$\par\noindent  \
$\rm 2\,ap\ bet^2\,(-1 + nu)^2\,y^{1-nu}$\par\noindent  \
$\rm -ap\ bet^2\,(-1 + nu)^2\,y^{1-2\,nu}$\par\noindent
}
\par
\newpage
\par
{\noindent\sl Pa20110808.nb}\hfill 3\ \ {}\\[-7pt] $\line(1,0){460}$\par\noindent
\par
{\footnotesize\tt\noindent\
\"{}The following is to prove all r\_i not= 0 in (4.59), and tilde r not= 0\"{}\par\noindent  \
\"{}Note:~xQn0[x]\,=\,Qnu0'[x],\,tqn0[x]\,=\,tilde~qnu0[x],\,tc0\,=\,tilde~cnu0,\,tQk[x]\,=\,tilde~Qknu0\"{}\par\noindent  \
xQn0[x\_] := tqn0[x] + tc0 y K0[x]{\rm\^{}}(-nu0)\par\noindent  \
Hn0[x\_] := 1/(ap y)H0[x] xQn0[x]\par\noindent  \
Kn0[x\_] := 1/ y K0[x] xQn0[x]\par\noindent  \
tQk[x\_] := tq[x] + gam H0[x] K0[x]{\rm\^{}}(1-k nu0)\par\noindent  \
tHk[x\_] := 1/(ap y)(H0[x] tQk'[x] - (1 + ap(1 - k nu0)) tQk[x] H0'[x])\par\noindent  \
tKk[x\_] := 1/(ap y)(ap K0[x] tQk'[x] - (1 + ap(1 - k nu0)) tQk[x] K0'[x])\par\noindent  \
tR7[x\_] := (1-k\,nu0)\,ap(tKk[x]\,Hn0'[x]+Hn0[x]\,tKk'[x])+(1-nu0)ap(Kn0[x]\,tHk'[x]+\par\noindent  \
\ \hspace*{15pt}\          tHk[x]Kn0'[x])-(1+ap(1-(k nu0+nu0)))(Hn0[x] tKk'[x]+tHk[x]Kn0'[x])\par\noindent  \
s0 = Factor[Coefficient[tR7[x],tc0]tc0];\par\noindent  \
s1 = Factor[Coefficient[tR7[x] tc0,tc0]];\par\noindent  \
tr0 = Factor[Coefficient[s0,gam]gam/K0[x]{\rm\^{}}(-nu0-k nu0)];\par\noindent  \
tr1 = Factor[Coefficient[s0 gam,gam]/K0[x]{\rm\^{}}(-nu0)];\par\noindent  \
tr2 = Factor[Coefficient[s1,gam]gam/K0[x]{\rm\^{}}(-k nu0)];\par\noindent  \
tr3 = Factor[Coefficient[s1 gam,gam]];\par\noindent  \
\"{}Lowest terms of tilde r\_i${\ssc\!}$,i=0${\ssc\!}$,1${\ssc\!}$,2${\ssc\!}$,3~are~y${\ssc\!}$,y${\ssc\!}${\rm\^{}}(1-k\,nu0)${\ssc\!}$,y${\ssc\!}${\rm\^{}}(1-nu0)${\ssc\!}$,y${\ssc\!}${\rm\^{}}(1-k\,nu0+nu0)${\ssc\!}$\"{}\par\noindent  \
tqn0[x\_] := -tc0 y{\rm\^{}}(1-nu0)\par\noindent  \
tq[x\_] := -gam x y{\rm\^{}}(1-k\,nu0);\par\noindent  \
H0[x\_] := x;\par\noindent  \
K0[x\_] := y;\par\noindent  \
ap := 1/(nu0-1)\par\noindent  \
\"{}The following shows that tilde r\_i, i=0,1,2,3 are nonzero\"{}\par\noindent  \
Factor[tr0]\par\noindent  \
Factor[tr1]\par\noindent  \
Factor[tr2]\par\noindent  \
Factor[tr3]\par\noindent  \
\"{}Below shows tilde r not= 0, where tilde r is in paragraph after (4.59)\"{}\par\noindent  \
Factor[tr1 + tr2 K0[x]{\rm\^{}}(nu0-k\,nu0)]\par\noindent  \
Clear[tqn0]\par\noindent  \
Clear[tq]\par\noindent  \
Clear[ap]\par\noindent  \
Clear[H0]\par\noindent  \
Clear[K0]\par\noindent \
}
\par{\small
\noindent\  The following is to prove all r\_i not= 0 in (4.59), and tilde r not= 0\par\noindent  \
Note: xQn0[x] = Qnu0'[x], tqn0[x] = tilde qnu0[x], tc0 = tilde cnu0, tQk[x] = tilde Qknu0\par\noindent  \
Lowest terms of tilde r\_i,i=0,1,2,3 are  y,y{\rm\^{}}(1-k\,nu0),y{\rm\^{}}(1-nu0),y{\rm\^{}}(1-k\,nu0+nu0)\par\noindent  \
The following shows that tilde r\_i, i=0,1,2,3 are nonzero\par\noindent  \
$\rm-2\, gam\,(-1+k\, nu0)\, tc0\ y$\par\noindent  \
$\rm2\,gam\,(-1+ k\,nu0)\, tc0\ y^{1 - k\,nu0}$\par\noindent  \
$\rm2\,gam\,(-1+ k\,nu0)\, tc0\ y^{1 - nu0}$\par\noindent  \
$\rm-2\,gam\,(-1+ k\,nu0)\, tc0\ y^{1 - nu0- k\,nu0}$\par\noindent  \
Below shows tilde r not= 0, where tilde r is in paragraph after (4.59)\par\noindent  \
$\rm4\,gam\,(-1+ k\,nu0)\, tc0\ y^{1 - k\, nu0}$
}
}
\par
\newpage\setcounter{page}{1}\pagenumbering{arabic}
\par}
\cl
 {{
 \bf  Poisson algebras, Weyl algebras and Jacobi \vspace*{-5pt}pairs 
 }\,\footnote{Supported by NSF grant 10825101 of
China, the Fundamental Research Funds for the Central Universities\vs{4pt}\nl\hspace*{4.5ex}{\it Mathematics Subject
Classification (2000):} 14R15, 14E20, 13B10, 13B25, 17B63}}
 \vs{5pt}
 \cl{{Yucai Su}} \cl{\small\it Department of Mathematics, Tongji University, Shanghai 200092,
China}
 \cl{\small\it Email: ycsu@tongji.edu.cn} \vskip3pt
%
\par

\def\P{{\cal W}}\def\NP{{\rm NP}}\noindent{\small\footnotesize
{\bf Abstract.~}~We study Jacobi pairs in details and obtained some properties.
We also study the natural Poisson algebra structure $(\PP,[\cdot,\cdot],\cdot)$ on
the space $\PP:=\C[y]((x^{-\frac1N}))$ for some sufficient large $N$, and
introduce
some automorphisms of $(\PP,[\cdot,\cdot],\cdot)$ which are (possibly infinite but well-defined) products of the automorphisms of forms  $e^{\ad_H}$ for $H\in x^{1-\frac1N}\C[y][[x^{-\frac1N}]]$ and $\tau_c:(x,y)\mapsto(x,y-cx^{-1})$ for some $c\in\C$. These automorphisms are used as tools to study Jacobi pairs in $\PP$. In particular, starting from
 a Jacobi pair $(F,G)$ in $\C[x,y]$  which violates the two-dimensional Jacobian conjecture, by
applying some variable change $(x,y)\mapsto\big(x^{b},x^{1-b}(y+a_1 x^{-b_1}+\cdots+a_kx^{-b_k})\big)$ for some $b,b_i\in\Q_+,a_i\in\C$ with $b_i<1<b$, we obtain a \QJ  pair still denoted by $(F,G)$ in $\C[x^{\pm\frac1N},y]$ with the form $F=x^{\frac{m}{m+n}}(f+F_0)$, $G=x^{\frac{n}{m+n}}(g+G_0)$ for some positive integers $m,n$,
and $f,g\in\C[y]$, $F_0,G_0\in x^{-\frac1N}\C[x^{-\frac1N},y]$, such that $F,G$ satisfy some additional conditions.
Then we generalize the results to the Weyl algebra $\WW=\C[v]((u^{-\frac1N}))$ with relation $[u,v]=1$, and obtain some properties of  pairs $(F,G)$ satisfying $[F,G]=1$,
referred to as Dixmier pairs. 
\vskip3pt

\vskip3pt

}

{\small\parskip1pt
\tableofcontents
\vskip4pt
\noindent{\bf Acknowledgements\hfill\pageref{ACKN}}
\vskip4pt
\noindent{\bf References\hfill\pageref{referencessss}}}\newpage

\section{Introduction}
\par
\def\C{{\mathbb{C}}}

It is a well known fact that if $n$ polynomials $f_1,...,f_n$ are
generators of the polynomial ring $\C[x_1,...,x_n]$, then the {\it
Jacobian determinant} $ 
J({\ssc\!}f_1,...,f_n{\ssc\!})$ $= {\rm det\,}A \in \C\bs\{0\} $ 
is a nonzero constant, where
$A=(\ptl_{x_j}f_i)_{i,j=1}^n$ is the $n\times n$ {\it
Jacobian matrix} of $f_1,...,f_n$. One of the major unsolved
problems of mathematics \cite{S} (see also \cite{B, CM, V2}), viz.~the {\it
Jacobian conjecture}, states that the reverse of the above statement
also holds, namely, if  $J(f_1,...,f_n)\in\C\bs\{0\}$, then
$f_1,...,f_n$ are
generators of $\C[x_1,...,x_n]$. 

This conjecture relates to many aspects of mathematics
[\ref{A}, \ref{DV}--\ref{H}, \ref{R}--\ref{SY}] 
and has attracted great attention in mathematical and physical
literatures during the past 60 years and there have been a various
ways of approaches toward the proof or disproof of this conjecture
(here we simply give a short random list of references
[\ref{Abh}, \ref{BCW}, \ref{CCS}, \ref{D}, \ref{J}--\ref{Ki}, \ref{M1}, \ref{V1}--\ref{W}]). 
Hundreds 
of papers have appeared
in connection with this conjecture,
even for the simplest case $n=2$ \cite{AO, N, No}.
However this conjecture remains
unsolved even for the case $n=2$.

Let $W_n$ 
be the {\it rank $n$ Weyl
algebra}, which is the associative unital algebra generated by $2n$
generators $u_1,...,u_n,v_1,...,v_n$ satisfying the relations
$[u_i,u_j]=[v_i,v_j]=0,\,[v_j,u_i]=\d_{ij}$, where the commutator $[\cdot,\cdot]$ is defined by
\equa{Brak}{[a,b]=ab-ba\mbox{ \ for }a,b\in W_n.} Under the commutator, $W_n$ becomes a Lie algebra,
denoted by $W^L_n$, called the {\it rank $n$ Weyl Lie algebra}. With
a history of 40 years, the {\it Dixmier conjecture} \cite{Di} states that
every nonzero endomorphism $W_n$ is an automorphism. This conjecture remains
open for $n\ge1$. It is well known \cite{AV, BK} that the rank $n$
Dixmier conjecture implies the $n$-dimensional Jacobi conjecture and
the $2n$-dimensional Jacobi conjecture implies the rank $n$ Dixmier
conjecture.

In this paper, we study Jacobi pairs in details and obtained some properties.
We also study the natural Poisson algebra structure $(\PP,[\cdot,\cdot],\cdot)$ on
the space $\PP:=\C[y]((x^{-\frac1N}))$ for some sufficient large $N$, and
introduce
some automorphisms of $(\PP,[\cdot,\cdot],\cdot)$ which are (possibly infinite but well-defined) products of the automorphisms of forms  $e^{\ad_H}$ for $H\in x^{1-\frac1N}\C[y][[x^{-\frac1N}]]$ and $\tau_c:(x,y)\mapsto(x,y-cx^{-1})$ for some $c\in\C$. These automorphisms are used as tools to study Jacobi pairs in $\PP$. In particular, starting from
 a Jacobi pair $(F,G)$ in $\C[x,y]$  which violates the two-dimensional Jacobian conjecture, by
applying some variable change $(x,y)\mapsto\big(x^{b},x^{1-b}(y+a_1 x^{-b_1}+\cdots+a_kx^{-b_k})\big)$ for some $b,b_i\in\Q_+,a_i\in\C$ with $b_i<1<b$, we obtain a \QJ  pair still denoted by $(F,G)$ in $\C[x^{\pm\frac1N},y]$ with the form $F=x^{\frac{m}{m+n}}(f+F_0)$, $G=x^{\frac{n}{m+n}}(g+G_0)$ for some positive integers $m,n$,
and $f,g\in\C[y]$, $F_0,G_0\in x^{-\frac1N}\C[x^{-\frac1N},y]$, such that $F,G$ satisfy some additional conditions (see Theorem \ref{Nerw}).

Then we generalize the results to the Weyl algebra $\WW=\C[v]((u^{-\frac1N}))$ with relation $[u,v]=1$, and obtain some properties of  pairs $(F,G)$ satisfying $[F,G]=1$,
referred to as {\it Dixmier pairs}. In particular,
one can define the Newton polygon $\NP(F)$ of $F$ as for the case of Jacobi pairs (cf.~Subsection \ref{S-NewtonP} and arguments after \eqref{p==Dix}).
 We can  suppose $\NP(F)$ has a vertex $(m_0,m)$ with $0<m_0<m$.
 First (as in the case of  Jacobi pairs), from the pair $(F,G)$, by applying some automorphism, we  obtain a
  Dixmier pair still denoted by $(F,G)$ in $\C[u^{\pm\frac1N},v]$
 with the form $F=u^{\frac{m}{m+n}}(f+F_0)$, $G=u^{\frac{n}{m+n}}(g+G_0)$ for some positive integers $m,n$,
and $f,g\in\C[v]$, $F_0,G_0\in u^{-\frac1N}\C[u^{-\frac1N},u]$, such that $F,G$ satisfy some additional conditions (cf.~Theorem \ref{Dixmier-lemm}).

 By applying the automorphism $(u,v)\mapsto(v,-u)$, we can assume $F$ has a vertex
 $(m_0, m)$ with $ m_0> m>0$, and we can further assume that the slope of the edge located at the right bottom side of $\Supp\, F$ and with the top vertex being $( m_0, m)$ is positive (as for the case of Jacobi pairs), which turns out to be $1$ (cf.~\eqref{p==Dix} and Theorem \ref{a-b====}). Furthermore, the coefficient of the term $u^{ m_0-1}v^{ m-1}$ is always $\frac{ m_0 m}{2}$ (if we assume $\Coeff(F,u^{m_0}v^{ m})=1$).
We remark that this is the place where the great difference between Newton polygons of Jacobi pairs and Dixmier pairs occurs; for the Jacobi pairs, an edge of the Newton polygon can never have slope $1$ (cf.~Theorem \ref{C[x,y]'}).

The main results in the present paper are summarized in  Theorems \ref{Jacobian-el}, \ref{C[x,y]'}, \ref{Nerw}, \ref{MMM-lemm}, Corollary \ref{Rema4.1} and Theorems \ref{Dixmier-lemm}, \ref{a-b====}.

\section{Definition of the prime degree $p$, notations and preliminaries}
\vs{0pt}\par
In this section, we first give some notations and definitions, then we present some preliminary results.
\subsection{Notations and definitions}
Denote by $\Z,\Z_+,\Z_-,\N,\Q,\Q_+,\C$ the sets of integers,
non-negative integers, negative integers, positive integers, rational numbers, non-negative rational numbers, complex numbers
respectively. Let $\C(x,y)=\{\frac{P}{Q}\,|\,P,Q\in\C[x,y]\}$ be the
field of rational functions in two variables. We use $\AA,\,\BB,\CC$
to denote the following rings (they are in fact fields and $\AA,\CC$
are algebraically closed fields):
\begin{eqnarray}\nonumber
\!\!\!\!&\!\!\!\!&\mbox{$\AA=\{f=\sum\limits_{i\in\Q}f_ix^i\,|\,f_i\in\C,\,\Supp_xf\subset
\a-\dis\frac{1}{\b}\Z_+\mbox{ \ for some }\a,\b\in\Z,\,\b>0\}, $}
\\\label{ring-ABC1--BB}
\!\!\!\!&\!\!\!\!&\mbox{$
\BB=\AA((y^{-1}))=
\{F=\sum\limits_{j\in\Z}F_jy^j\,|\,F_j\in\AA,\,\Supp_yF\subset\a-\Z_+\mbox{
\ for some }\a\in\Z\}, $}
 \\[-3pt]
\label{CC-TAU}
 \!\!\!\!&\!\!\!\!&\mbox{$
\CC=\{F\!=\!\sum\limits_{j\in\Q}F_jy^j\,|\,F_j\in\AA,\,\Supp_yF\subset\a\!-\!\dis\frac{1}{\b}\Z_+\mbox{
 for some }\a,\b\!\in\!\Z,\,\b\!>\!0\}, $}
\end{eqnarray}
where \equan{smaller-supp}{\begin{array}{ll}
\Supp_xf=\{i\in\Q\,|\,f_i\ne0\}&\mbox{(the {\it support of $f$}),}\\[5pt]
\Supp_yF=\{j\in\Z\,|\,F_j\ne0\}&\mbox{(the {\it support of $F$ with
respect to $y$}).}\end{array}}
 Denote $\ptl_x=\frac{\ptl}{\ptl x}$ or
$\frac{d}{dx},\,\ptl_y=\frac{\ptl}{\ptl y}$. For
$F=\sum_{i\in\Q,j\in\Z}f_{ij}x^iy^j\in\BB$, we define
\begin{eqnarray}\label{supp-F}&& \Supp\,F=\{(i,j)\,|\,f_{ij}\ne0\}\mbox{ \
(called the {\it support
of } $F$)},
\end{eqnarray}and define \begin{eqnarray}
&&\deg_xF=\max\{i\in\Q\,|\,f_{ij}\ne0\mbox{ \ for some }j\}
\mbox{ \ (called the {\it $x$-degree} of $F$)},\nonumber\\
&&\deg_yF=\max\{j\in\Z\,|\,f_{ij}\ne0\mbox{ \ for some }i\}
\mbox{ \ (called the {\it $y$-degree} of $F$)},\nonumber\\
&&\deg\,F=\max\{i+j\in\Q\,|\,f_{ij}\ne0\}
\mbox{ \ (called the {\it total degree} of $F$)}.\nonumber\end{eqnarray}
 Note
that a degree can be $-\infty$ (for instance, $F=0$), or $+\infty$
(for instance, $\deg_xF=\deg\,F=+\infty$ if $F=\sum_{i=0}^\infty
x^{2i}y^{-i}$).
An element $F=\sum_{i=0}^\infty f_iy^{m-i}$ is {\it
monic} if $f_0=1$.

For any $h=x^{\a}+\sum_{i=1}^\infty
h_ix^{\a-\frac{i}{\b}}\in\AA$ with $\a\in\Q$, $\b\in\N$, and for any $a\in\Q$, we define
$h^a$ to be the unique element in $\AA$:
\begin{equation}\label{h-a-power}
\mbox{$h^a=x^{a\a}\big(1+\sum\limits_{i=1}^\infty
h_ix^{-\frac{i}{\b}}\big)^a=x^{a\a} \sum\limits_{j=0}^\infty\big(\
^{^{\dis a}}_{_{\dis j}}\,\big)\big(\sum\limits_{i=1}^\infty
h_ix^{-\frac{i}{\b}}\big)^j=x^{a\a}+\sum\limits_{i=1}^\infty
h_{a,i}x^{a\a-\frac i\b},$}
\end{equation}
where the coefficient of $x^{a\a-\frac{i}{\b}}$, denoted by
$\Coeff(h^a,x^{a\a-\frac{i}{\b}})$, is
\begin{equation*}
\Coeff(h^a,x^{a\a-\frac{i}{\b}})=h_{a,i}=\mbox{$\sum\limits_{^{\sc
r_1+2r_2+\cdots+ir_i=i}_{\sc\ \ \ \,r_1,...,r_i\ge0}}$}
\Big(\begin{array}{c}a\\
r_1,r_2,...,r_i\end{array}\Big) h_1^{r_1}h_2^{r_2}\cdots h_i^{r_i},
\end{equation*}
and\begin{equation*}
\Big(\begin{array}{c}a\\
r_1,r_2,...,r_i\end{array}\Big)=
\frac{a(a-1)\cdots(a-(r_1+\cdots+r_i)+1)}{r_1!\cdots r_i!},
\end{equation*}
is a {\it multi-nomial coefficient}. Note that if
$\Supp_xh\subset\a-\frac{1}{\b}\Z_+$, then $\Supp_xh^a\subset
a\a-\frac1{\b}\Z_+$. Similarly, for any $F=\sum_{i=0}^\infty
f_iy^{m-i}\in\BB$ with $f_0\ne0$, and for any $a,b\in\Z,\,b\ne0$ with
$b|am$, we can define $F^{\frac{a}{b}}$ to be the unique element in
$\BB$ (note that if $b\mbox{$\not|$}\,am$, we can still define $F^{\frac{a}{b}}$, but in this case it is in $\CC$ instead of $\BB$):
\begin{eqnarray}
\label{equa2.2}
\!\!\!\!\!\!\!\!\!\!&\!\!\!\!\!\!\!\!\!\!F^{\frac{a}{b}}
\!\!\!&=f_0^{\frac{a}{b}}y^{\frac{am}{b}}\big(1+\mbox{$\sum\limits_{i=1}^\infty$}
f_0^{-1}f_{i}y^{-i}\big)^{\frac{a}{b}}
=
f_0^{\frac{a}{b}}y^{\frac{am}{b}}+\mbox{$\sum\limits_{j=1}^\infty$}
f_{a,b,j}y^{\frac{am}{b}-j},
\end{eqnarray}
\VS{-7pt}where
\begin{equation}
\label{equa2.3}
f_{a,b,j}=\Coeff(F^{\frac{a}{b}},y^{\frac{am}{b}-j})=\mbox{$\sum\limits_{^{\sc
r_1+2r_2+\cdots+jr_j=j}_{\sc\ \ \ \ r_1,...,r_j\ge0}}$}
\Big(\begin{array}{c}\frac{a}{b}\\
r_1,r_2,...,r_j\end{array}\Big) f_1^{r_1}\cdots
f_j^{r_j}f_0^{\frac{a}{b}-(r_1+\cdots+r_j)}.
\end{equation}

\begin{defi}
\label{defi2.1}\rm (cf.~Remark \ref{rema-prime}) Let $p\in\Q$ and $F=\mbox{$\sum_{i=0}^\infty$}
f_iy^{m-i}\in\BB$ with $f_0\ne0$. If
\begin{equation}\label{p-degree}\mbox{ $\deg_xf_i\le\deg_xf_0+pi$ for all
$i$ with equality holds for at least one $i\ge1$,}\end{equation}
then $p$, denoted by $p(F)$, is called the {\it prime degree of
$F$}. We set $p(F)=-\infty$ if $F=f_0y^m$, or set $p(F)=+\infty$ if
it does not exist (clearly, $p(F)<+\infty$ if $F$ is a polynomial).
\end{defi}

Note that the definition of $p:=p(F)$ shows that the support
$\Supp\,F$  of $F$, regarded as a subset of the   plane $\R^2$, is
located at the left side of  the {\it prime line}
$L_F:=\{(m_0,m)+z(p,-1)\,|\,z\in\R\}$ (where $m_0=\deg_xf_0$)
passing the point $(m_0,m)$ and at least another point $(i,j)$ of $\Supp\,F$
(thus $-p^{-1}$ is in fact the slope of the prime line $L_F$)\VS{-5pt}:
\begin{equation}\label{p-line}\ \put(-50,30){$\hspace*{-125pt}
(m_0,m)\put(90,-30){$L_F$}
\put(33,-15){$\circ$}\put(14,-22){$(i,j)$}\!\!\circ\!\!\line(3,-1){80}\put(-50,-50){$\Supp\
F$}\put(-120,-20){$\circ\ \ \circ\ \ \circ$}\put(-120,-35){$\circ\ \
\circ\ \ \circ\ \ \circ$}\put(-120,-50){$\circ\ \ \circ\ \ \circ\ \
\circ$}\ \hspace*{80pt}\
\put(0,0){\put(90,-30){$F\comp{0}$}$\line(3,-1){80}$}
\put(0,-15){\put(90,-30){$F\comp{r}$}$\line(3,-1){80}$}
\put(0,-30){$\line(3,-1){80}\put(10,-35){Components}$}$}
\end{equation}
\begin{rema}\label{rema-prime}\rm
It may be more proper to define the prime degree $p$ to be \begin{equation}\label{p=tobe}
\raisebox{-8pt}{$\stackrel{\dis{\rm sup}}{\sc i\ge1}$}\,\frac{\deg_xf_i-\deg_xf_0}{i}.\end{equation}
Then in case $F\!\in\!\C[x^{\pm1},y^{\pm1}]$, both definitions coincide. However, when $F\!\notin\!\C[x^{\pm1},y^{\pm1}]$, it is possible that the prime line $L_F$ defined as above only passes through one point of $\Supp\,F$.
Since we do not like such a case to happen when we consider \QJ  pairs in later sections, we use  (\ref{p-degree}) to define $p$ instead of
(\ref{p=tobe}). (For example, for $F=y+\sum_{i=0}^\infty x^iy^{-i}$, if we use (\ref{p=tobe}) to define $p$,
it would be $1$; but if we use (\ref{p-degree}) to define $p$, it is $+\infty$.)
\end{rema}

Let $p\ne\pm\infty$ be a fixed rational number. We always assume {\bf all} elements under consideration below have prime degrees $\le p$ (and in the next section, we always
take $p=p(F)$).

\begin{defi}\label{comp}\rm\begin{itemize}\parskip-5pt\item[(1)] Let \equan{F=========}
{F=\sum\limits_{i=0}^\infty f_iy^{m-i}\in\BB \mbox{ with
$f_0\ne0$}.}
In what follows, we always use $m$ to denote $m=\deg_yF$ and use $m_0$ to denote
$\deg_xf_0=m_0$  {until \eqref{x-yyyyy}}. We call $x^{m_0}y^m$ the {\it first term of $F$}.
Suppose $p(F)\le p$. For $r\in\Q$, we define the {\it $p$-type
$r$-th component} (or simply the $r$-th component) of $F$ to be
\begin{equation}\label{eq-comp}\mbox{$F\comp{r}=\sum\limits_{i=0}^\infty\Coeff(f_i,x^{m_0+r+ip})x^{m_0+r+ip}y^{m-i}$},\end{equation} which simply collects those terms $f_{ij}x^iy^j$ of $F$ with
$(i,j)$ located in a line parallel to the prime line $($cf.~(\ref{p-line}).
One
immediately sees that $F\comp{r}=0$ if $r>0$ and $F\comp{0}\ne0$. We
remark that if $p(F)>p$ then it is possible that $F\comp{r}\ne0$ for
$r>0.$
\item[(2)] Suppose $F\!=\!\sum_{i=0}^m f_iy^{m-i}\!\in\! \F((x^{-1}))[y]$ with $f_0$ being a monic Laurent polynomial of $x$ and $p(F)=p$.
 \begin{description}\parskip-3pt\item[{\rm(i)}]
We  call $F\comp{0}$  the {\it leading polynomial} of $F$, and
$F\comp{<0}:=\sum_{r<0}F\comp{r}$ the {\it ignored polynomial} of
$F$.
\item[{\rm(ii)}]
We always use the bold symbol $\pri{F}\in\C[x^{\pm1}][y]$ to denote
the unique monic polynomial of $y$ $($with coefficients being
Laurent polynomials of $x)$ such that
$F\comp{0}=x^{m_0}\pri{F}^{m'}$ with $m'\in\N$ maximal $($we always
use $m'$ to denote this integer$)$. Then (a polynomial satisfying
(\ref{pri-pro}) is usually called a {\it power free polynomial}),
\begin{equation}\label{pri-pro}\mbox{$\pri{F}\ne H^k$ \ \ for any \ $H\in\F[x^{\pm1}][y]$ and
$k>1$.}\end{equation} We call $\pri{F}$ the {\it primary polynomial}
of $F$. We always use $d$ to denote
\begin{equation}\label{denote-d}d=\deg_y\pri{F},\mbox{\ \ thus \ }
m'=\frac{m}{d}.
\end{equation}\end{description}\item[(3)]
 An element of the form
$$H=\mbox{$\sum\limits_{i=0}^{\infty}$}
c_ix^{r_0+ip}y^{\a_0-i}\in\BB\mbox{ \  with \ } c_i\in\F,\,\
r_0\in\Q,\,\ \a_0\in\Z_+,$$ $($i.e., its support is located in a
line$)$ is called a {\it $p$-type quasi-homogenous element}
(q.h.e.), and it is called a {\it $p$-type quasi-homogenous
polynomial} (q.h.p.) if it is a polynomial.
\end{itemize}
\end{defi}

\begin{lemm}
\label{lemm2.02}
\begin{itemize}\parskip-2pt\item[$(1)$] $p(F)=p(F^{\frac{a}{b}})$ if $a,b\ne0$.
\item[$(2)$] $p(FG)\le\max\{p(F),p(G)\}$ with equality holds if
one of the following conditions holds:
\begin{itemize}\parskip-2pt\item[{\rm(i)}] $p(F)\ne p(G)$ or \item[{\rm(ii)}] $F,G$ both are polynomials, or
\item[{\rm(iii)}] $F,G$ are $p$-type q.h.e.~such that $FG$ is not a
monomial.
\end{itemize}
\item[{\rm(3)}] Suppose $p(F)\le p$. Then $p(F)=p$ if and only if
$F\comp{0}$ is not a monomial.
\end{itemize}
\end{lemm}
\noindent{\it Proof.~}~(1) We immediately see from (\ref{equa2.3})
that $\deg_x f_{a,b,i}\le \deg_x f_0^{\frac{a}{b}}+i{\ssc\,} p(F)$, i.e.,
$p(F^{\frac{a}{b}})\le p(F)$. Thus also, $p(F)=p\big((F^{\frac{a}{b}})^{\frac{b}{a}}\big)\le p(F^{\frac{a}{b}})$.
Hence $p(F^{\frac{a}{b}})=p(F).$

(2) and (3) are straightforward to verify.\hfill$\Box$\vskip7pt

Note that the equality in Lemma \ref{lemm2.02}(2) does not
necessarily hold in general; for instance,
$$F=(y+x^2)(y+x^3)=y^2+(x^2+x^3)y+x^5,\ \ \ G=(y+x^3)^{-1},$$ with
$p(F)=p(G)=3$, but $FG=y+x^2$ with $p(FG)=2$.

\subsection{Some preliminary results}
We remark that the requirement that any element under consideration
has prime degree $\le p$ is necessary, otherwise it is possible that
in (\ref{product-HK}), there exist infinite many $r>0$ with
$H\comp{r}\ne0$ and the right-hand side becomes an infinite sum.
\def\tildem{\tilde m}\begin{lemm}\label{lemm-product-KH}
\begin{itemize}\parskip-4pt
\item[{\rm(1)}]
Suppose $F=\sum_{i=0}^\infty F_i$, where $F_i=\sum_{j=0}^\infty
\tilde f_{ij}y^{\tildem_i-j}\in\BB$ $($and all have prime degree
$\le p)$ with $\tilde f_{i0}\ne0$, $\deg_x\tilde
f_{i0}\!=\!\tildem_{i0}$ such that $\tildem_0\!>\!\tildem_1\!>\!...$
$($in this case $\sum_{i=0}^\infty F_i$ is called
{\textbf{summable}}$)$. Then
\begin{equation}\label{eq-lemm-1}\mbox{$F\comp{r}=\sum\limits_{i=0}^\infty
(F_i)\comp{r+\tildem_{00}-\tildem_{i0}+p(\tildem_0-\tildem_i)}.$}\end{equation}
\item[{\rm(2)}] An element is a $p$-type q.h.e. $\Lra$ it is a component of itself.
\item[{\rm(3)}] If $F$ is a $p$-type q.h.e.~with $\deg_yF=m$, then
$F^{\frac{a}{b}}$ is a $p$-type q.h.e.~if $b|am$.\item[{\rm(4)}]
 Let
$H,K\in\BB$ with prime degrees $\le p$, and  $r\in\Q$. Then $($where ``\,${}\comp{r}$\,'' is defined as in \eqref{eq-comp}$)$
\begin{equation}\label{product-HK}(HK)\comp{r}= \mbox{$\sum\limits_{r_1+r_2=r}$}
H\comp{r_1}K\comp{r_2}.\end{equation}
\item[{\rm(5)}] Let $F$ be as in Definition $\ref{comp}(2)$. Let $\ell \in\Z$ with $d|\ell$. Then
for all $r\in\Q$,
$x^{-\frac{m_0\ell}{m}}(F^{\frac{\ell}{m}})\comp{r}$ $\in\C(x,y)$ is
a rational function of the form $\pri{F}^aP$ for some $a\in\Z$ and
$P\in\C[x^{\pm1}][y]$.
\item[{\rm(6)}] Let $P,Q\in\BB$ with prime degrees $\le p$ and $Q\ne0$.
Then each $p$-type component of the rational function $R=\frac{P}{Q}$ is a rational function. Furthermore, there exists some $\a\in\N$ such that the $p$-type $r$-th component $R\comp{r}=0$ if
$r\notin\frac1{\a}\Z$.
\end{itemize}
\end{lemm}
\noindent{\it Proof.~}~Using (\ref{eq-comp}), (\ref{equa2.2}) and
(\ref{equa2.3}), it is straightforward to verify (1)--(3).

 (4)
Suppose $H=\sum_{i=0}^\infty h_iy^{\a-i},\,K=\sum_{i=0}^\infty
k_iy^{\b-i}$ with $\deg_xh_0=\a_0,\,\deg_xk_0=\b_0$. Then $HK=
\mbox{$\sum_{i=0}^\infty$} \chi_i y^{\a+\b-i}$ with $\chi_i=
\sum_{s_1+s_2=i}h_{s_1}k_{s_2}$. Thus
$$\mbox{$\Coeff(\chi_i,x^{\a_0+\b_0+r+ip})=\sum\limits_{r_1+r_2=r}\
\sum\limits_{s_1+s_2=i}\Coeff(h_{s_1},x^{\a_0+r_1+s_1p})\Coeff(k_{s_2},x^{\b_0+r_2+s_2p})$.}$$
Hence we have (4).

 (5) We have $F=\sum_{0\ge s\in\Q}F\comp{s}$ and
$F\comp{0}=x^{m_0}\pri{F}^{m'}.$ Thus
$$
\begin{array}{ll}
x^{-\frac{m_0\ell}{m}}F^{\frac{\ell }{m}}\!\!\! &
= \pri{F}^{\frac{\ell }{d}} \bigl(1+\sum\limits_{s<0}
x^{-m_0}F\comp{s}\pri{F}^{-m'}\bigr)^{\frac{\ell }{m}}
\\[12pt]
&= \sum\limits_{0\ge r\in\Q}\ \sum\limits_{{}^{^{\sc
s_1t_1+\cdots+s_kt_k=r}_{\sc\ \ \,0>s_1>\cdots>s_k}}_{\sc \ \ \
t_i\in\Z_+,\,k\ge0}}
\Big(\begin{array}{c}\frac{\ell }{m}\\
t_1,...,t_k\end{array}\Big)
x^{-(t_1+\cdots+t_k)m_0}\prod\limits_{i=1}^k\bigl(F\comp{s_i}\bigr)^{t_i}\pri{F}^{\frac{\ell}{d}
-(t_1+\cdots+t_k)m'}.
\end{array}
$$
By (4), if the $j$-th component of $(F\comp{s_i})^{t_i}$ is nonzero,
then $j=s_it_i$. By (2) and (3), if the $j$-th component of
$\pri{F}^{\frac{\ell}{d} -(t_1+\cdots+t_k)m'}$ is nonzero, then
$j=0.$ Thus by (1) and (4), the $r$-th component of $x^{-\frac{m_0\ell}{m}}F^{\frac{\ell
}{m}}$ for $r\le0$ is
\begin{equation}\label{equa-f-xy-r}
x^{-\frac{m_0\ell}{m}}(F^{\frac{\ell }{m}})\comp{r} =
\mbox{$\sum\limits_{{}^{^{\sc s_1t_1+\cdots+s_kt_k=r}_{\sc\ \
\,0>s_1>\cdots>s_k}}_{\sc\ \ \ t_i\in\Z_+,\,k\ge0}}$}
\Big(\begin{array}{c}\frac{\ell }{m}\\
t_1,...,t_k\end{array}\Big)x^{-(t_1+\cdots+t_k)m_0}
\mbox{$\prod\limits_{i=1}^k$}\bigl(F\comp{s_i}\bigr)^{t_i}\pri{F}^{\frac{\ell}{d}
-(t_1+\cdots+t_k)m'},
\end{equation}
which is a finite sum of rational functions of $y$ with coefficients
in $\C[x^{\pm1}]$ by noting that every component $F\comp{s_i}$ is a
polynomial of $y$ and that the powers of $\pri{F}$ in
(\ref{equa-f-xy-r}) are integers and $\{r\,|\,F\comp{r}\ne0\}$ is a
finite set.

(6) By 
(4) and (5),
\begin{equation}\label{FFFF}R\comp{r}=\mbox{$\sum\limits_{r_1+r_2=r}$}P\comp{r_1}(Q^{-1})\comp{r_2}.\end{equation}
Since every component of the polynomial $P$ is a polynomial, and $P$
has only finite nonzero components, thus the sum in (\ref{FFFF}) is
finite. By (5), $(Q^{-1})\comp{r_2}$ is a rational function. Thus we
have the first statement of (6). The second statement follows from
(\ref{equa-f-xy-r}) and (\ref{FFFF}). \hfill$\Box$\vskip7pt

Equation (\ref{equa-f-xy-r}) in particular gives
\begin{equation}\label{0-th-com}
(F^{\frac{\ell}{m}})\comp{0}=x^{\frac{m_0\ell}{m}}\pri{F}^{\frac{\ell}{d}}.
\end{equation}
\begin{rema}\label{rational=F}\rm Suppose $F=\frac{P}{Q}$ is a rational function  such that
$\deg_yF\ne0$ and $P,Q$ have the same prime degree $p$ and the same
primary polynomial $\pri{F}$, then
$P\comp{0}=x^{a}\pri{F}^{b},\,Q\comp{0}=x^{c}\pri{F}^d$ for some
$a,b,c,d\in\Z$, so $F\comp{0}=x^{a-c}\pri{F}^{b-d}$. In this case we
also call $\pri{F}$ the primary polynomial of $F$ if $b\ne d$.\end{rema}

 The
result in Lemma \ref{lemm-product-KH}(5) can be extended to rational
functions as follows.
\begin{lemm}\label{ration-comp}
Let $F,G\in\C((x^{-1}))[y]$ with prime degree $p$ and primary polynomial
$\pri{F}$. Let $x^{m_0}y^m,$ $x^{n_0}y^n$ be the first terms of $F$
and $G$. Let $a,b\in\Z,\,\check F=F^aG^b\in\C((x^{-1}))(y)$ with $\check
m:=\deg_y\check F=am+bn\ne0$. Set $\check m_0=am_0+bn_0$. Let
$\ell\in\Z$ with $d|\ell$. Then for all $r\in\Q$, $x^{-\frac{\check
m_0\ell}{\check m}}(\check F^{\frac{\ell}{\check m}})\comp{r}
\in\C(x,y)$ is a rational function of the form $\pri{F}^{{\ssc\,}\a} P$ for some
$\a\in\Z$ and $P\in\C[x^{\pm1}][y]$.
\end{lemm}\ni{\it Proof.~}~By Lemma \ref{lemm-product-KH}(5) (by taking $\ell=am$ or $bm$), each
component of $F^a,\,G^b$ is a rational function of the form
$\pri{F}^{{\ssc\,}\a} P$. Thus the ``\,$\check{\,}$\,'' version of
(\ref{equa-f-xy-r}) (which is still a finite sum by Lemma
\ref{lemm-product-KH}(6)) shows that we have the result.
\hfill$\Box$\vskip7pt

The following result generalized from linear algebra will be used in
the next section.
\begin{lemm}
\label{ration-function} Suppose $H\in\C((x^{-1}))[y]$ such that
$\deg_yH>0$. Let $m\in\N$. Suppose there exists a finite nonzero
combination
\begin{equation}\label{TATA}\mbox{$P:=\sum\limits_{i\in\Z}$}p_iH^{\frac{i}{m}}\in\C((x^{-1}))(y)
\mbox{ \ for some \ }p_i\in\C((x^{-1})).
\end{equation}
Then $H=H_1^{\frac{m}{d}}$ is the $\frac{m}{d}$-th power of some
polynomial $H_1\in\C((x^{-1}))[y]$, where $d=\gcd(m,i\,|\,p_i$ $\ne0)$
is the greatest common divisor of the integer set
$\{m,i\,|\,p_i\ne0\}$.
\end{lemm}
\noindent{\it Proof.~}~We thank Dr.~Victor Zurkowski who
suggested the following simple proof.
Applying operator $\frac{H}{\ptl_yH}\ptl_y$ to (\ref{TATA})
iteratively, we obtain a system of equations (regarding
$p_iH^{\frac{i}{m}}$ \VS{-3pt}as unknown variables).~Solving the system
shows $p_iH^{\frac{i}{m}}\!\in\!\C((x^{-1}))(y)$ for all $i.$ Then induction on
$\#A$ (where $A\!=\!\{i\,|\,p_i\!\ne\!0\}$) gives the result (when $\#A\!=\!1$,
it is a well-known result of linear algebra, one can also simply
prove as follows: Suppose $H^{\frac{n}{m}}\!=\!\frac{R}{Q}$ for some
coprime polynomials $R,Q$, then $Q^mH^n\!=\!R^m$. Since
$\C((x^{-1}))[y]$ is a uniquely factorial domain, by decomposing
each polynomial into the product of its irreducible polynomials,
 we see that $H$ is the $\frac{m}{d}\,$-th power of some
polynomial).
 \hfill$\Box$\vskip7pt

We shall also need the following lemma.\def\mmm{\infty}
\begin{lemm}\label{+new-lamm}
Let $\b\in\C\bs\{0\}$. Let $A_\mmm=\C[z_1,z_2,...]$ be the polynomial
ring with $\mmm$ variables, and denote $z_0=1$. \VS{-5pt}Set
\begin{equation}\label{h-m-i-0} \mbox{$
h_{\b}(x)=\Big(1+\sum\limits_{i=1}^\mmm z_i
x^i\Big)^\b=\sum\limits_{j=0}^\infty h_{\b,j}x^j\in A_\mmm[[x]]$},
\end{equation} \VS{-5pt}where 
\begin{equation*}
\mbox{$h_{\b,i} ={\Coeff}(h_{\b}(x),x^i)=\sum\limits_{\stackrel{\sc
i_1+2i_2+\cdots+mi_m=i}{^{\,}_{\sc m\in\N,\,i_1,...,i_m\in\Z_+}}}
\Big(\begin{array}{cc}\b\\
i_1,...,i_m\end{array}\Big)z_1^{i_1}\cdots z_m^{i_m}\mbox{ \ for \
}i\in\Z_+. $} \VS{-6pt}\end{equation*} Then for all $r\in\Z_+$, we
\VS{-5pt}have
\begin{eqnarray}\label{h-m-i-1}&&
\mbox{$\sum\limits_{s=0}^\mmm(s(\b+1)-r)z_sh_{\b,r-s}=0 \mbox{ \ \ and
\ }\sum\limits_{s=0}^\mmm z_sh_{\b,r-s}=h_{\b+1,r}
 $}
.
\end{eqnarray}
\end{lemm} \ni{\it Proof.~}~Taking derivative with respect to $x$ in (\ref{h-m-i-0}), \VS{-5pt}we
have
$$\mbox{$
\sum\limits_{j=0}^\infty
jh_{\b,j}x^{j-1}=\b\Big(\sum\limits_{i=0}^\mmm
z_ix^i\Big)^{\b-1}\Big(\sum\limits_{j=0}^\mmm jz_jx^{j-1}\Big)$}.
$$
Multiplying it by $\sum_{i=0}^\mmm z_ix^i$ and comparing the
coefficients of $x^{r-1}$ in both sides, we obtain
$$\mbox{$\sum\limits_{s=0}^\mmm(r-s)z_sh_{\b,r-s}=\b\sum\limits_{s=0}^\mmm s{\ssc\,}z_sh_{\b,r-s}$},
$$
i.e., the first equation of (\ref{h-m-i-1}) holds. To prove the
second, simply write $h_{\b+1}(x)$ as $h_{\b+1}(x)=h_{\b}(x)h_{1}(x)$
and compare the coefficients of $x^r$. \hfill$\Box$\vskip7pt

Let $F,G\in\BB$ such that
\begin{equation}\label{monic-F-G}
\mbox{$F=\sum\limits_{i=0}^\infty f_iy^{m-i},\ \ \
G=\sum\limits_{j=0}^\infty g_jy^{n-j}\mbox{ \ \ for some \
$f_i,\,g_j\in\AA$, \ $f_0g_0\ne0$}\VS{-7pt},$}\end{equation}
with $m=\deg_yF>0$ and $n=\deg_yG$,
we can express $G$ \VS{-7pt}as
\begin{equation}
\label{equa2.5} G=\mbox{$\sum\limits_{i=0}^\infty$}
b_iF^{\frac{n-i}{m}} \mbox{ \ \ for some \ }b_i\in\AA,
\end{equation}
where by comparing the coefficients of $y^{n-i}$, $b_i$ can be
inductively determined by the following (cf.~(\ref{equa2.3}))\VS{-5pt}:
\begin{eqnarray}
\label{equa2.6-}\!\!\!\!\!\!&\!\!\!\!\!\!\!\!\!\!\!\!\!\!\!\!&
b_i=h^{-\frac{m'_0(n-i)}{m}}\big(g_i-\mbox{$\sum\limits_{j=0}^{i-1}$}b_j f_{n-j,m,i-j}\big),\mbox{ \ \ or more precisely,}\\[-3pt]
\label{equa2.6}\!\!\!\!\!&\!\!\!\!\!\!\!\!\!\!\!\!\!\!\!\!&
b_i\!=\!h^{-\frac{m'_0(n-i)}{m}}\!\Big(\!g_i\!-\!\mbox{$\sum\limits_{j=0}^{i-1}$}b_j\!
\mbox{$\sum\limits_{^{ r_1+2r_2+\cdots+mr_m=i-j}_{\ \ \ \
r_1,...,r_m\ge0}}$}\! \big(\!\begin{array}{c}\frac{n-j}{m}\\
r_1,r_2,...,r_m\end{array}\!\!\big) f_1^{r_1}\cdots
f_m^{r_m}h^{\frac{m'_0(n-j)}{m}-(r_1+\cdots+r_m)m'_0}\!\Big).\ \ \
\end{eqnarray}
Similarly, we can express the polynomial $y$ \VS{-5pt}as\def\varphi{{\bar b}}
\begin{eqnarray}
\label{equa2.7} y= \mbox{$\sum\limits_{i=0}^\infty$} \varphi_i
F^{\frac{1-i}{m}},
\end{eqnarray}
where $\bar b_i\in\AA$ is determined \VS{-5pt}by
\begin{equation}
\label{equa2.7+}
\varphi_i\!=\!h^{-\frac{m'_0(1-i)}{m}}\!\Big(\!\d_{i,0}\!-\!\mbox{$\sum\limits_{j=0}^{i-1}$}\,\varphi_j
\mbox{$\sum\limits_{^{r_1+2r_2+\cdots+mr_{m}=i-j}_{\ \ \ \ \ \
r_1,...,r_{m}\ge0}}$}
\!\!\big(\!\!\begin{array}{c}\frac{1-j}{m}\\
r_1,r_2,...,r_{m}\end{array}\!\!\big)\! f_1^{r_1}\cdots
f_m^{r_m}h^{\frac{m'_0(1-j)}{m}-(r_1+\cdots+r_m)m'_0}\!\Big).
\end{equation}
 We
 observe a simple fact that  if $F,G\in\AA[y]$ then $F^{-1}$ does not appear in the
expression of $G$ in (\ref{equa2.5}) (we would like to thank
Professor Leonid Makar-Limanov, who told us the following more
general fact and the suggestion for the simple proof).

 \setcounter{clai}{0}
\begin{lemm}\label{clai-b+m+n}
Let $F,G\in\AA[y]$ be any polynomials
with $\deg_yF\!=m$, $\deg_yG=n$. Express $G$ as
in $(\ref{equa2.5})$. We always have $b_{i m+n}=0$ for $i\ge1$,
namely, all negative integral power of $F$ cannot appear in the
expression.
\end{lemm}
\ni{\it Proof.~}~The proof can be obtained by regarding $F$ and $G$
as polynomials in $y$, from an observation that $\int G{\ssc\,} dF$ is a
polynomial while the term with $F^{-1}$ would require the
logarithmic term in integration. Iteration of this observation of
course shows that the coefficients with all negative integral powers
of $F$ are zeros.\hfill$\Box$ \vskip7pt
More generally, we have
\begin{lemm}\label{f-f-1=0}
Let $F\in\BB$ 
with $\deg_yF=m>0$ and   $k,\ell\in\Z$ with $m\ne-k$. Then in the
expression  \equa{AQAQ}{\mbox{$y^\ell
F^{\frac{k}{m}}=\sum\limits_{i=0}^\infty
c_{\ell,i}F^{\frac{k+\ell-i}{m}}$ with $c_{\ell,i}\in\AA$, }} the
coefficient of $F^{-1}$ \VS{-5pt}is
\begin{equation}\label{coeff-F-1}
c_{\ell,m+k+\ell}=-\frac{\ell}{m+k}\Coeff(F^{\frac{m+k}{m}},y^{-\ell}).\end{equation}
In particular, the coefficient $c_{\ell,i}$ in
$y^\ell=\sum_{i=0}^\infty c_{\ell,i}F^{\frac{\ell-i}{m}}$ is
\begin{equation}\label{COEFF-c}
c_{\ell,i}=-\frac{\ell}{i-\ell}\Coeff(F^{\frac{i-\ell}{m}},y^{-\ell})\mbox{
\ \ for all \ } \ell\ne i.
\end{equation}
We also have $c_{0,0}=1$ and
\begin{equation}\label{COEFF-c-0}
c_{\ell,\ell}\!=\! -\mbox{${\dis\frac{1}{m}}\sum\limits_{s=0}^\infty$}
s{\sc\,}\Coeff(F,y^{m-s})\Coeff(F^{-1},y^{-m-\ell+s})\!=\!\frac{1}{m}\Coeff(F^{-1}\ptl_yF,y^{-\ell})
\ \mbox{if  }\ell\!>\!0.
\end{equation}
\end{lemm}\ni{\it Proof.~}~For any $A\in\BB$ with
$\deg_yA=a$, we always  use (until the end of this proof) $A_r$ to
denote the homogenous part of $A$ of $y$-degree $a-r$, i.e.,
\begin{equation}\label{denote-A-r}A_r=\Coeff(A,y^{a-r})y^{a-r},\ \ \ r\in\Z.\end{equation}
Comparing the homogenous parts of $y$-degree $k+\ell-(r-s)$ in (\ref{AQAQ}),
 multiplying the result by $(s(\frac{k+\ell}{m}+1)-r)F_s$ and taking sum over all $s$, we obtain (note that both sides are homogenous
polynomials of $y$ of degree $m+k+\ell-r$)
\begin{eqnarray}\label{aaa-aaa}
\mbox{$\sum\limits_{s\in\Z_+}\!\!\Big(s\big({\dis\frac{k\!+\!\ell}{m}}\!+\!1\big)\!-\!r\Big)F_s(y^\ell
F^{\frac{k}{m}})_{r-s}$}
\!=\! \sum\limits_{i=1}^\infty
c_{\ell,i}\sum\limits_{s\in\Z_+}\Big(s\big(\frac{k\!+\!\ell}{m}\!+\!1\big)\!-\!r\Big)F_s(F^{\frac{k+\ell-i}{m}})_{r-i-s},
\end{eqnarray}
where on the right-hand side, the terms with $i=0$ vanish  by (\ref{h-m-i-1}) (
we take
$z_r=\frac{F_r}{F_0},\,\b=\frac{k+\ell}{m}$ in (\ref{h-m-i-1})), and $(y^\ell
F^{\frac{k}{m}})_{r-s}$, $(F^{\frac{k+\ell-i}{m}})_{r-i-s}$ are the
notations as in (\ref{denote-A-r}). Note from definition
(\ref{denote-A-r}) that $(y^\ell
F^{\frac{k}{m}})_{r-s}=y^\ell(F^{\frac{k}{m}})_{r-s}$. For $i\ne
m+k+\ell$, we write $s(\frac{k+\ell}{m}+1\big)-r$ in  both sides of
(\ref{aaa-aaa}) respectively \VS{-5pt}as
\begin{eqnarray}\label{WAWASSS}&\!\!\!\!\!\!\!\!\!\!\!\!\!\!\!\!\!\!\!\!\!\!\!\!\!\!\!\!\!\!\!\!&
s(\frac{k+\ell}{m}+1\big)-r=\frac{m+k+\ell}{m+k}\Big(s\big(\frac{k}{m}+1\big)-r\Big)+\frac{r\ell}{m+k},
\\\label{WAWASSS1}&\!\!\!\!\!\!\!\!\!\!\!\!\!\!\!\!\!\!\!\!\!\!\!\!\!\!\!\!\!\!\!\!&
s(\frac{k+\ell}{m}+1\big)-r=\frac{m+k+\ell}{m+k+\ell-i}\Big(s\big(\frac{k+\ell-i}{m}+1\big)-(r-i)\Big)+\frac{(r-m-k-\ell)i}{m+k+\ell-i},
\end{eqnarray}
and substitute them into (\ref{aaa-aaa}) (using (\ref{h-m-i-1})
again, the first terms in the right-hand sides of (\ref{WAWASSS})
and (\ref{WAWASSS1}) then become vanishing), we obtain (where the first
equality follows from the second equation of (\ref{h-m-i-1})\VS{-5pt})
\begin{eqnarray}\label{aaa-aaa+1}
&\!\!\!\!\!\!\!\!\!\!\!\!
\mbox{${\dis\frac{r\ell}{m\!+\!k}}y^\ell(F^{\frac{m+k}{m}})_r$}\!\!&=
\mbox{${\dis\frac{r\ell}{m\!+\!k}}y^\ell\sum\limits_{s\in\Z_+}
F_s(F^{\frac{k}{m}})_{r-\ell-s}$}\nonumber\\
&\!\!\!\!\!\!\!\!\!\!\!\!&
\mbox{$= \sum\limits_{1\le i\ne m+k+\ell}
c_{\ell,i}\sum\limits_{s\in\Z_+}{\dis\frac{(r\!-\!m\!-\!k\!-\!\ell)i}{m\!+\!k\!+\!\ell\!-\!i}}F_s(F^{\frac{k+\ell-i}{m}})_{r-i-s}$}
\nonumber\\
\!\!\!\!\!\!\!\!\!\!\!\!&\!\!\!\!\!\!\!\!\!\!\!\!&
\mbox{$\phantom{=}\,+\,c_{\ell,m+k+\ell}\sum\limits_{s\in\Z_+}
\Big(s\big({\dis\frac{k+\ell}{m}}+1\big)-r\Big)F_s(F^{-1})_{r-m-k-\ell-s}$}.
\end{eqnarray}
Take $r=m+k+\ell$ in  (\ref{aaa-aaa+1}), then the first summand in
the right-hand side is $0$. As for the last summand, if $s\ne0$,
then $r-m-k-\ell-s<0$ and so $(F^{-1})_{r-m-k-\ell-s}=0$ (the
definition of (\ref{denote-A-r}) shows $A_r=0$ if $r<0$). Thus
(\ref{aaa-aaa+1}) \VS{-5pt}gives
\begin{equation}\label{LAST-ee}
\frac{(m+k+\ell)\ell}{m+k}y^\ell(F^{\frac{m+k}{m}})_{m+k+\ell}=-c_{\ell,m+k+\ell}(m+k+\ell)F_0(F^{-1})_0=-(m+k+\ell)
c_{\ell,m+k+\ell}.
\end{equation} Since $m+k+\ell>0$ by our assumption (cf.~\eqref{aaa-aaa+1}), (\ref{LAST-ee})
 together with
definition (\ref{denote-A-r}) shows we have (\ref{coeff-F-1}). To
prove (\ref{COEFF-c-0}), assume $\ell>0$. We take $k=-m,\,r=\ell$ in
(\ref{aaa-aaa}) (but we do not substitute (\ref{WAWASSS}) into the
left-hand side of (\ref{aaa-aaa})), then the left-hand side of
(\ref{aaa-aaa}) becomes (using the fact $\sum_{s\in\Z_+}
F_s(F^{-1})_{\ell-s}=(F^{1-1})_\ell=0$ by the second equation of
(\ref{h-m-i-1})\VS{-5pt})\begin{equation}\label{M13}\mbox{$y^\ell
\sum\limits_{s\in\Z_+}
s{\dis\frac{\ell}{m}}F_s(F^{-1})_{\ell-s}$}.\end{equation} Thus the
left-hand side of (\ref{LAST-ee}) should be (\ref{M13}). Hence we
have
 (\ref{COEFF-c-0}).
 \hfill$\Box$

\section{Jacobi pairs and Jacobian elements}
In this section, we discuss properties of Jacobi pairs and Jacobian elements in details. The main results of this section are Theorems \ref{Jacobian-el}, \ref{C[x,y]'} and \ref{Nerw}.
\subsection{General discussions on \QJ  pairs in $\BB$}
\begin{defi}\label{gen-JP}\rm \begin{itemize}\parskip-3pt
\item[(1)]
Define the Lie bracket $[\cdot,\cdot]$ on $\BB$ by the Jacobian determinant
\equa{Lie-b}{[F,G]=J(F,G)=(\ptl_{x} F)(\ptl_{y}G)-(\ptl_{y}F)( \ptl_{x}G)
\mbox{ \ for }F,G\in\C[x,y].}
It is well-known (e.g., \cite{MTU,SX}) that the triple $(\BB,[\cdot,\cdot],\cdot)$ (where $\cdot$ is the usual product in $\BB$) is a {\it Poisson  algebra}, namely $(\BB,[\cdot,\cdot])$ is a Lie algebra, $(\BB,\cdot)$ is a commutative associative algebra, and the following compatible Leibniz rule holds:
\equa{Leib}{[F,GH]=[F,G]H+G[F,H]\mbox{ \ for }F,G,H\in\BB.}
\item[(2)] A pair $(F,G)$  is called a {\it quasi-Jacobi pair} (or simply a {\it \QJ  pair}), if the following conditions
are satisfied\VS{-5pt}:\begin{itemize}\parskip-3pt\item[(i)]
$F,G\in\BB$ 
are of the form \eqref{monic-F-G}
such that $[F,G]\in
\C\bs\{0\}$, \item[(ii)] $F$ has prime degree
$p\ne\pm\infty$, 
\item[(iii)]
$p(G)\le p$.
\end{itemize} \end{itemize} \end{defi}
\begin{rema}\label{MAM3.1}\rm\begin{itemize}\parskip-3pt\item[(1)]
We always denote (until \eqref{x-yyyyy})\equa{MAMSMSMSMS}{\mbox{$m:=\deg_yF,\ \ \ \,m_0:=\deg_xf_0$,\ \ \
$n:=\deg_yG,\ \ \ \,n_0:=\deg_xg_0$,}} and use notations $b_i,\bar b_i$ in \eqref{equa2.5}--\eqref{equa2.7+}.
We always assume $m\ge 1$ and $m_0\ge0$ (but not
necessarily $m_0\in\Z$). Note that we always have $m_0\ne m$ (by Lemma \ref{Jacobian-el} since $(-\frac{m_0}{m},-1)\in\Supp\,F^{-\frac{1}{m}}$).
Also note that $n$ can be negative, but the
nonzero Jacobian determinant requires that $m+n\ge1$ (cf.~Lemma
\ref{lemm2.3}). 
\item[(2)]
If $F,G\in\C[x,y]$ with $[F,G]\in\C\bs\{0\}$, then $(F,G)$ is
called a (usual) {\it Jacobi pair}. In this case if necessary by
exchanging $F$ and $G$, we can always suppose $p(G)\le p$ (in fact
$p(G)=p$ if $m+n\ge2$ by Lemma \ref{C[x,y]}). Thus a usual Jacobi pair
is necessarily a quasi-\QJ  pair.
\item[(3)]We always suppose $f_0,g_0$ are monic (i.e., the coefficients of the highest powers of $x$ in $f_0,g_0$ are $1$).
If  $f_0\in\C[x]$ and  $f_0\ne1$, let $x=a$ be a root of $f_0$. By
applying the automorphism $(x,y)\mapsto(x\!+\!a,y)$, we can suppose
$f_0$ does not contain the constant term. In this case, we also let
$h\!\in\!\C[x]$ be the unique monic  polynomial such that $f_0\!=\!h^{m'_0}$
with $m'_0\!\in\!\N$ maximal, and set $d_0\!=\!\frac{m_0}{m'_0}$. \VS{-5pt}Then
\begin{equation}\label{h---}h\ne h_1^k\mbox{ \ for any }h_1\in\C[x],\,k>1.\end{equation}
If $f_0\!=\!1$, we set $h\!=\!x,\,d_0\!=\!1,\,m'_0\!=\!m_0\!=\!0$. If $f_0\!\notin\!\C[x]$,
we simply set $h\!=\!f_0,\,d_0\!=\!m_0,\,m'_0\!=\!1$.
\item[(4)]
We always fix notations $h,\,m'_0,\,d_0$.\end{itemize}\end{rema}
\begin{exam}\label{example-1}\rm For any $m_0,m\in\N,\,a,b\in\C$ with
$m_0>
m$ and $a,b\ne0$, the pairs
$$\begin{array}{ll}
(F,G)=\big(x^{m_0}y^m+a
x^{m_0-1}y^{m-1},(x^{m_0}y^{m}+ax^{m_0-1}y^{m-1})^2+b
x^{1-m_0}y^{1-m}\big),\\[6pt]
(F_1,G_1)=\big(x^4(y^2+x^{-1})^3,\,x^{-2}(y^2+x^{-1})^{-2}y\big),
\end{array}$$ are \QJ  pairs with
$p(F)=-1$ (but if $m_0<m$, then $p(G)=-\frac{3m_0-1}{3m-1}>-1$), and
$p(F_1)=p(G_1)=-\frac{1}{2}.$ Note that $p(G)$ can be smaller than
$p(F)$; for instance, if we replace $G$ by $x^{1-m_0}y^{1-m}$ then
$p(G)=-\infty<p(F)$ (but in this case $\deg_yF+\deg_yG=1$).
\end{exam}

\begin{defi}\rm\label{trace-jaco}\begin{itemize}\parskip-3pt\item[(1)]
We introduce a notion, referred to as the {\it trace}, of an element
$F\in\BB$ to be
(cf.~\eqref{tr-F-def} and \eqref{A-tr} for the reason why we call it ``trace'')
\equa{Res-x-y}{\tr\,F=\Res_x\Res_y\,F=\Coeff(F,(xy)^{-1})\mbox{ for }F\in\BB.}
If $\tr \,F=0$, we say $F$ has the {\it vanishing trace property}.\item[(2)]An element $F\in\BB$ is called a {\it Jacobian element} if there is
$G\in\BB$ such that $[F,G]\in\C\bs\{0\}$.
\end{itemize}\end{defi}

\begin{rema}\label{residue=0}\rm
Note that it is important to consider the residue  of an element in a specific space, otherwise the  residue may be different; for instance, if we regard the element $H=\frac1{xy+(xy)^2}$ as in $\C[x^{\pm1}]((y))$ then $H=-\sum_{i=-1}^\infty(-xy)^{i}$ and its residue (with respect to $x$) is $y^{-1}$, but if
we consider it as in the space $\C[x^{\pm1}]((y^{-1}))$ then $H=\sum_{i=2}^\infty(-xy)^{-i}$ and its residue is zero.
\end{rema}

The following Theorem \ref{Jacobian-el}(2)
characterizes Jacobian elements.
\begin{theo}\label{Jacobian-el}\begin{itemize}\parskip-3pt\item[\rm(1)]
Assume $(F,G)$ is a \QJ  pair in $\AA[y]$. Suppose
\equa{Res-x-y=0} {\res_x(F\ptl_xG):=\Coeff(F\ptl_xG,x^{-1})=0.}
Let $x=u(t),y=v(t)\in\C((t^{-1}))$, then
$\res_t(F\frac{dG}{dt})=J{\ssc\,}\res_t(u\frac{d v}{dt}).$\item[\rm(2)]
 An element $F\in\BB$ is a
Jacobian element if and only if $F\notin\C$ and
$(F-c)^a$ has the vanishing trace property for all $a\in\Q$ and $c\in\C$
$($note that $(F-c)^a$ might not be in $\BB$ but in $\CC{\sc\,})$.\item[\rm(3)]Let $F\in\BB$ be a Jacobian element, and $H\in\BB$. There exists some $K\in\BB$ such that $H=[K,F]$ if and only if $HF^a$ has the vanishing trace property for all $a\in\Q$.\item[\rm(4)]Let $(F,G)$ be a \QJ pair in $\BB$ with $[F,G]=1$.  One has
\equa{Traaa}{\tr([H,F]G)=\tr\,H\mbox{ for any $H\in\BB$}.}
\end{itemize}\end{theo}
\ni{\it Proof.~}~(1) Assume
 $F=\sum_{i\in\Q,j\in\Z_+}f_{ij}x^iy^j$ and
 $G=\sum_{k\in\Q,\ell\in\Z_+}g_{k\ell}x^ky^\ell$. Then
$[F,G]=J$ and \eqref{Res-x-y=0} imply
\equa{f-ij-g-kl=}{\sum\limits_{i+k=a,\,j+\ell=b}(i\ell-jk)f_{ij}g_{k\ell}=\d_{a,1}\d_{b,1}J,\ \ \ \ \ \
\sum\limits_{i+k=0,\,j+\ell=b}kf_{ij}g_{k\ell}=0,\ \ \forall\,a,b\in\C.}
Denoting $\ptl_t=\frac{d}{dt}$, and considering 4 cases (noting that $j=\ell=0$ if $j+\ell=0$): (i) $i+k=0=j+\ell$; (ii) $i+k=0\ne j+\ell$; (iii) $i+k\ne0=j+\ell$; and (iv) $i+k\ne0\ne j+\ell$, we obtain from \eqref{f-ij-g-kl=},
\begin{eqnarray*}
\!\!\!\!\!\!\!&\!\!\!\!\!\!\!\!\!\!\!\!\!\!\res_t(F\ptl_tG)
\!\!\!&=\res_t\Big(\mbox{$\sum\limits_{i,j,k,\ell}$}
f_{ij}g_{k\ell}(ku^{i+k-1}v^{j+\ell}\ptl_tu+\ell u^{i+k}v^{j+\ell-1}\ptl_tv)\Big)\nonumber\\
\!\!\!\!\!\!\!&\!\!\!\!\!\!\!&
=\res_t\Big(\mbox{$\sum\limits_{i+k\ne0,\,j+\ell>0}$}f_{ij}g_{k\ell}\Big(
\frac{\ell}{j+\ell}\ptl_t(u^{i+k}v^{j+\ell})+
\frac{jk-i\ell}{j+\ell}u^{i+k-1}v^{j+\ell}\ptl_tu\Big)\Big)\nonumber\\
\!\!\!\!\!\!\!&\!\!\!\!\!\!\!&=J{\ssc\,}\res_x(v\ptl_tu).
\end{eqnarray*}
One can also give the following simple proof as suggested by Dr.~Victor Zurkowski: \eqref{Res-x-y=0} means
 $F\ptl_xG =\ptl_x U$ for some $U\in\AA[y]$.
Then
$\ptl_y\ptl_xU = - J + \ptl_x(F\ptl_yG)$, so
$\ptl_yU + Jx - F\ptl_yG$ is a polynomial in $y$, denoted by $g$, and we can write
$F\ptl_yG -\ptl_y U = Jx + \ptl_yg$. Thus
$F\ptl_t G = F\ptl_xG\ptl_t u + F\ptl_yG\ptl_t v =\ptl_t ( U + g)  + J u \ptl_tv$.
Taking the residue with respect to $t$, the term from the exact derivative drops and one gets the result.

(2) ``\,$\Longrightarrow$\,'': First note that
 by computing the coefficient of $x^{-1}y^{-1}$ in the Jacobian determinant $[H,K]$, one immediately \VS{-5pt}sees
\equa{J1--}{\tr([H,K])=0\mbox{ for any $H,K\in\CC$}.} Suppose $F,G\in\BB$ such that $[F,G]=1$.
Then $(F\!-\!c)^a=[\frac{(F-c)^a}{a},(F\!-\!c)^{2-a}G]$ (assume $a\ne0$), we obtain the necessity by \eqref{J1--}.

``\,$\Longleftarrow$\,'': First assume $\deg_yF\ne0$.
Replacing
$F$ by $F^{-1}$ if necessary, we can suppose $m:=\deg_y F\in\N$. Express $y$ as in (\ref{equa2.7}), then in fact $\bar b_i$ is
 $c_{1,i}$ in (\ref{AQAQ}). Thus (\ref{COEFF-c}) and the
sufficiency condition show that $x^{-1}$ does not appear in $\bar
b_i$ for all $i\ne1$. Therefore there exist $b_i\in\AA$ such that
$\ptl_xb_i=-\frac{1-i}{m}\bar b_i$ for all $i\ge0$ (cf.~(\ref{PTLBBB})). Set $G=\sum_{i=0}^\infty b_iF^{\frac{1-m-i}{m}}$. Then the proof
of Lemma \ref{lemm2.3} shows $[F,G]=1$.

Now assume $\deg_yF=0$.
Let $f_i:=\Coeff(F,y^i)$ for $i\le0$.
By replacing $F$ by $F-c$ for some $c\in\C$ if necessary, we can suppose $f_0\ne0$ does not contain the constant term.
We can write $f_{-1}=\b\ptl_xf_0+h$ for some $\b\in\C$ and $h\in\AA$ such that
either $h=0$ or else $i:=\deg_xf_0-1\ne j:=\deg_xh$. If $h\ne0$, by taking $a:=\frac{i-j}{i+1}\ne0$ ($i+1=\deg_xf_0\ne0$ since $\Coeff(f_0,x^0)=0$), we have $F^{a}=(f_0+f_{-1}y^{-1}+\cdots)^{a}=f_0^{a}\big(1+af_0^{-1}(\b\ptl_xf_0+h)y^{-1}+\cdots\big)
=f_0^{a}+(\b\ptl_x(f_0^{a})+af_0^{a-1}h)y^{-1}+\cdots$, and $(-1,-1)\in\Supp\,F^a$ by noting that $\ptl_x(f_0^{a})$ does not contain the term $x^{-1}$ but $\deg_xf_0^{a-1}h=-1$, a contradiction. Thus $h=0$.
Now one can uniquely determine $g_1=(\ptl_xf_0)^{-1}$ and $g_{-i}$ for $i\ge1$ inductively such that $G=g_1y+\sum_{i=1}^\infty g_{-i}y^{-i}$ satisfying $[F,G]=1$, by comparing the coefficients of $y^{-i}$ for $i\ge0$.

(3) Suppose $H=[K,F]$. Then by \eqref{Traaa}, we see $HF^a=[KF^a,F]$ has the  vanishing trace property for all $a\in\Q$. Conversely suppose $HF^a$ has  the  vanishing trace property for all $a\in\Q$.
Let $G$ be as constructed in the proof of (1) such that $[F,G]=1$ and $m+n=1$, $m_0+n_0=1$ (cf.~notation \eqref{MAMSMSMSMS}).
Using \eqref{equa2.7} and the fact that $FG$ has the highest term $xy$, we see that any element $x^iy^j$ with $i,j\in\Q$
can be written as $F^kG^\ell+F_1$, where $k,\ell\in\Q$ satisfy $km_0+\ell n_0=i,$ $km+\ell n=j$, and some $F_1\in\CC$ with $\deg_yF_1<j$. Thus by induction on $\deg_yH$, we see that
$H$ can be expressed as $H=\sum_{i,j\in\Q}c_{ij}F^iG^j$ for some $c_{ij}\in\C$.
Note that $\tr(F^kG^\ell)=0$ if and only if $(k,\ell)\ne(-1,-1)$. Thus
$G^{-1}$ cannot appear in the expression of $H$, otherwise $HF^{-i-1}$ would have nonzero trace if $c_{i,-1}\ne0$.
Taking $K=\sum_{i,j\in\Q}\frac{c_{ij}}{j+1}F^iG^{j+1}$, we obtain $H=[K,F]$.

(4)  one has
\begin{eqnarray}\label{2}
&\tr([H,F]G)\!\!\!&=
\tr(\ptl_xH\ptl_yFG-\ptl_yH\ptl_xFG)\nonumber\\&&=
\tr\Big(-H(\ptl_x\ptl_yFG+\ptl_yF\ptl_xG)+H(\ptl_y\ptl_xFG+\ptl_xF\ptl_yG)\Big)=\tr\, H.
\nonumber\end{eqnarray}This completes the proof of the theorem.
\hfill$\Box$\vskip5pt
As a by-product of the above theorem, one can easily obtain
\begin{coro}\label{CCCC}\begin{itemize}\parskip-3pt\item[\rm(1)]
If $(F,G)$ is a Jacobi pair in $\C[x,y]$ such that $\Supp\,F$ has a vertex $(m_0,m)$ with $m_0,m>0$
$($cf.~\eqref{new-supp-F-region}$)$, then $F,G$ are not generators of $\C[x,y]$ $($thus in particular, the proof of the two-dimensional Jacobian conjecture is equivalent to proving that a Jacobi pair $(F,G)$ in $\C[x,y]$ with $\Supp\,F$ having a vertex $(m_0,m)$ with $m_0,m>0$ does not exists$)$.
\item[\rm(2)] Any $\si$ in \eqref{VAR} regarded as an automorphism of $\C[x^{\pm\frac1N}]((y^{-\frac1N}))$ does not change the vanishing trace property.
\end{itemize}
\end{coro}\noindent{\it Proof.~}~(1) Say $[F,G]=1$ and assume $F,G$ are generators of $\C[x,y]$. Then $H=x^{m_0-1}y^{m-1}=\sum_{i,j\in\Z_+}a_{ij}F^iG^j$ for some $a_{ij}\in\C$, and so $H=[K,F]$ for $K=\sum_{i,j\in\Z_+}\frac{a_{ij}}{j+1}F^iG^{j+1}$. However,
$\tr(HF^{-1})\ne0$, a contradiction with Theorem \ref{Jacobian-el}(3).

(2) An automorphism $\si$ of the form \eqref{VAR} is a product of automorphisms of forms
$\si_1:(x,y)\mapsto (x,y+\l x^{-\frac{p}{q}})$ and
$\si_2:(x,y)\mapsto(x^{\frac{q}{q-p}},x^{\frac{-p}{q-p}}y)$ for some $\l\in\C$ and $p,q\in\Z$ with $p<q$.
Thus we can suppose  $\si=\si_1$ or $\si_2$.
In the first case, $\si$ is simply the exponential operator $e^{\ad_h}$ for $h=\frac{q\l}{q-p}  x^{1-\frac{p}{q}}$ (cf.~\eqref{exppp}), which
does not change the vanishing trace property by \eqref{J1--}. As for the later case, one can check directly that
$\si_2(x^iy^j)=x^{\frac{qi-pj}{q-p}}y^j\ne (xy)^{-1}$ if $(i,j)\ne(-1,-1)$.
\hfill$\Box$\vskip5pt

{}Now let $(F,G)$ be a \QJ  pair in $\BB$ and we use notations in \eqref{equa2.5}--\eqref{equa2.7+}. {}From (\ref{equa2.7}), we \VS{-5pt}obtain
\begin{equation}\label{equa-ptl-f-y}
\frac{1}{\ptl_y F}= \mbox{$\sum\limits_{i=0}^\infty$}
\frac{1-i}{m}\varphi_i F^{\frac{1-i}{m}-1}.
\end{equation}
\begin{lemm}
\label{lemm2.3} We have $m+n\ge1$, and $b_i\!\in\!\C$ if
$i\!<\!n\!+\!m\!-\!1$. Furthermore\VS{-5pt},
\begin{equation}\label{PTLBBB}\mbox{$\dis\ptl_xb_{i+m+n-1}=-\frac{(1-i)J}{m}\bar b_i$ \ \ if
\ $i\ge0$.}\end{equation} In particular\VS{-5pt},
\begin{equation}\label{b=m+n-1=}
\ptl_xb_{m+n-1}=-\frac{J}{m}\bar
b_0=-\frac{J}{m}h^{-\frac{m'_0}{m}}\mbox{ \ $($by
$(\ref{equa2.7+}))$}.
\end{equation}
\end{lemm}
\ni{\it Proof.~}~By (\ref{equa2.5}), $\ptl_xG=A+B\ptl_xF$,
$\ptl_yG=B\ptl_yF$, \VS{-5pt}where
$$\mbox{$A=\sum\limits_{i=0}^\infty (\ptl_xb_i)
F^{\frac{n-i}{m}}$,\ \ \
$B=\sum\limits_{i=0}^\infty{\dis\frac{n-i}{m}}b_iF^{\frac{n-i}{m}-1}$.}$$
We \VS{-5pt}obtain
\begin{equation}\label{compar-coeff}\mbox{$J=-A\ptl_yF$, \ i.e., \ $\sum\limits_{i=0}^\infty
(\ptl_xb_i)\dis
F^{\frac{n-i}{m}}=-\frac{J}{\ptl_yF}$.}\end{equation}  The lemma
follows from (\ref{equa-ptl-f-y}) by comparing the coefficients of
$F^{\frac{n-i}{m}}$ in (\ref{compar-coeff}) for
$i\in\Z_+$.\hfill$\Box$
\begin{rema}\label{m=1}\rm
{}From Lemma
\ref{lemm2.3}, we see that in case $F,G\in\C[x,y]$ and $m=1$ then we
can replace $G$ by $G-\sum_{i=0}^{n-1}b_iF^{n-i}$ to reduce
$\deg_yG$ to zero, thus obtain  $F,G$  are generators of $\F[x,y]$
and $F$ is a monic polynomial of $y$. Since eventually we shall
consider a usual Jacobi pair, from now on we suppose $m\ge2$.
Furthermore, we can always suppose $n\ge m$ and $m\mbox{$\not|$}\,n$.
\end{rema}

\def\aaa{\si}Using (\ref{equa2.7+}), we obtain $\deg_x\bar b_0=-\frac{m_0}{m}\ne-1$
(by Theorem \ref{Jacobian-el} since $(-\frac{m_0}{m},-1)\notin\Supp\,F^{-\frac1m}$), and for $i>0$,
\begin{eqnarray}\label{equa2.7+-+0}
\deg_x\varphi_i\!\!\!&\le\!\!\!&-\frac{m_0(1-i)}{m}+\rb{-3pt}{$^{\,\dis\max}_{_{\sc0\le
j<i}}$}\{\deg_x\bar b_j+\frac{m_0}{m}(1-j)+(i-j)p\} \nonumber\\
&\le\!\!\!&\si_i,\mbox{ \ \ where \
}\aaa_i:=i\big(p+\frac{m_0}{m}\big)-\frac{m_0}{m},\end{eqnarray} where,
the last inequality follows from induction on $i$. Thus
\begin{equation}
\label{equa2.7+-+} 1+\si_0\le\deg_xb_{m+n-1}\le\max\{0,1+\si_0\},\ \
\ \deg_xb_{m+n-1+i}\le \max\{0,1+\aaa_i\}\mbox{ \ for }i\ge1.
\end{equation}
We denote
\begin{equation}\label{deno-regard}\mbox{$\eta_i=\deg_xb_{m+n-1+i}$ \ for \ $i\ge0$.}\end{equation}
 We shall use (\ref{eq-lemm-1}) and (\ref{equa2.5}) to
compute $G\comp{r}$. Thus we set $F_i=b_iF^{\frac{n-i}{m}}$. Then
$p(F_i)\le p$ by Lemma \ref{lemm2.02}(1) and (2). Assume $m+n\ge2$.
Then $b_0\in\C$, and we have the data $(\tildem_{i0},\tildem_i)$ in
Lemma \ref{lemm-product-KH}(1) being $\tildem_i=n-i$, and
\begin{eqnarray*}
&&
\tildem_{i0}=\frac{m_0(n-i)}{m}\mbox{ if }i<m+n-1, \mbox{ and
}\\
&&\tildem_{m+n-1+i,0}=\frac{m_0(1-m-i)}{m}+\eta_i\mbox{
if }i\ge0.\end{eqnarray*}
 Then (\ref{eq-lemm-1}) gives
\begin{eqnarray}\label{G's-component}
\!\!\!\!\!\!\!\!\!\!\!\!&\!\!\!\!\!\!\!\!\!\!\!\!&\mbox{$G\comp{r}=
\sum\limits_{i=0}^\infty(F_i)\comp{r+\tildem_{00}-\tildem_{i0}+p(\tildem_0-\tildem_i)}$}\nonumber\\
\!\!\!\!\!\!\!\!\!\!\!\!&\!\!\!\!\!\!\!\!\!\!\!\!
&\mbox{$\phantom{G\comp{r}}=\sum\limits_{i=0}^{m+n-2}
(F_i)\comp{r+(p+\frac{m_0}{m})i}
+\sum\limits_{i=0}^\infty(F_{m+n-1+i})\comp
{\a_{r,i}}, \mbox{ where }$}\\
\label{G's-component-r}
\!\!\!\!\!\!\!\!\!\!\!\!&\!\!\!\!\!\!\!\!\!\!\!\!
&\a_{r,i}:=r+(p+\frac{m_0}{m})(m+n-1)+\frac{m_0}{m}+\si_i-\eta_i
.
\end{eqnarray}

\begin{lemm}\label{b-i=0-if}Assume $m+n\ge2$. If $p<-\frac{m_0}{m}$, then $b_i=0$ for $1\le i<m+n-1$.
\end{lemm}\ni{\it Proof.~}~Suppose there exists the smallest $i_0$ with
$1\le i_0<m+n-1$ such that $b_{i_0}\ne0$ (i.e., $F_{i_0}\ne0$). Then
setting $r=-(p+\frac{m_0}{m})i_0>0$ in (\ref{G's-component}) gives a
contradiction (cf.~statements after (\ref{eq-comp})):
\begin{equation}\label{0=Gnot=0}0=G\comp{r}=b_{i_0}(F^{\frac{n-i_0}{m}})\comp{0}+\cdots\ne0,\end{equation}
where the last inequality follows from the fact that $y^{n-i_0}$
appears in $(F^{\frac{n-i_0}{m}})\comp{0}$ but not in any omitted
terms.\hfill$\Box$\vskip7pt

We denote \begin{equation}\label{denote-b'-i}b'_i=\Coeff(b_i,x^0)
\mbox{ \ \ (which is \ $b_i$ \ if \ $i<m+n-1$).}\end{equation} We
always suppose that our \QJ  pair $(F,G)$ satisfies the
condition (by Lemma \ref{C[x,y]}, this condition will be
automatically satisfied by Jacobi pairs in $\C[x,y]$)
\begin{equation}\label{condition-p}p\ne-\frac{m_0}{m}.\end{equation}
 We always denote $\mu$ to be (until the end of this section)
\begin{equation}\label{mu-set}\mu=m+n-\frac{1+p}{\frac{m_0}{m}+p}\in\Q.\end{equation}
\begin{lemm}\label{b'-i=0}
Assume $m+n\ge2$ and $p<-\frac{m_0}{m}$. Then $b'_i=0$ for all
$i\ge1$ $($thus in particular $\eta_i=1+\deg_x\bar
b_i\le1+\si_i$ for all $i\ge0)$.
\end{lemm}
\ni{\it Proof.~}~Suppose $b'_{i_0}\ne0$ for smallest $i_0\ge1$. Then
$i_0\ge m+n-1$ by Lemma \ref{b-i=0-if}. Setting
$r=-(\frac{m_0}{m}+p)i_0>0$ in (\ref{G's-component}) gives a
contradiction as in (\ref{0=Gnot=0}). \hfill$\Box$

\begin{lemm}\label{p-ge---}
Assume $m\!+\!n\!\ge\!2$. Then $p\!\ge\!
-\frac{m_0}{m}\!+\!\frac{m-m_0}{m(m+n-1)}\!=\!-\frac{m_0+n_0-1}{m+n-1}$, where $n_0\!:=\!\deg_xg_0$.\end{lemm}\ni{\it
Proof.~}~First by letting $i=0$ in
(\ref{equa2.6-}), we have $n_0=\deg_xg_0=\frac{m_0}{m}n$, from which we have $-\frac{m_0}{m}+\frac{m-m_0}{m(m+n-1)}=-\frac{m_0+n_0-1}{m+n-1}$.

Assume conversely $p<-\frac{m_0}{m}+\frac{m-m_0}{m(m+n-1)}$.
\VS{-5pt}Then
$$\begin{array}{ll}r\!\!\!&\dis:=-\Big((p+\frac{m_0}{m})(m+n-1)+\frac{m_0}{m}+\si_0-\eta_0\Big)
>-(1+\si_0-\eta_0)=-1+\frac{m_0}{m}+\eta_0=0,\end{array}$$
where the last equality follows by noting that either $m_0>m$ (so in this case $p<-\frac{m_0}{m}+\frac{m-m_0}{m(m+n-1)}<-\frac{m_0}{m}$ and  $b'_{m+n-1}=0$ by Lemma \ref{b'-i=0}, thus
$\eta_0=1-\frac{m_0}{m}$ by \eqref{b=m+n-1=}), or else
$m_0<m$ (in this case $\eta_0=1-\frac{m_0}{m}$ again by \eqref{b=m+n-1=}).
By Lemma
\ref{b-i=0-if}, we have either $p+\frac{m_0}{m}\ge0$ or $b_i=0$ for
$0<i<m+n-1$. Thus
(\ref{G's-component}) gives a contradiction:
$$0\!=\!G\comp{r}=F\comp{r}\!+\!(F_{m+n-1})\comp{\a_{r,0}}\!+\!\cdots=(F_{m+n-1})\comp{0}\!+\!\cdots\ne0.
\eqno\Box$$\vskip7pt

Note from Lemma \ref{p-ge---} that $(\frac{m_0}{m}+p)(m+n)-(1+p)\ge0$. Thus we obtain
%
\begin{equation}\label{mu<0}\mu=m+n-\frac{1+p}{\frac{m_0}{m}+p}\le0\mbox{ \ \ if \
}p<-\frac{m_0}{m}.\end{equation}
\begin{rema}\label{m+n>1-im}\rm \begin{itemize}\parskip-4pt\item[(1)] It is very important to assume
$m+n\ge2$ (otherwise $b_0\notin\C$) and assume $p(G)\le p$
(otherwise $G\comp{r}$ can be nonzero for $r>0$). \item[(2)] Note
that if $p$ is sufficiently small, one cannot replace $G$ by $\check
G:=F^\ell+G$ for $\ell>0$ without changing $p(G)$; for instance, if
$p<-\frac{m_0}{m}$ then $p(\check G)=-\frac{m_0}{m}>p$.\end{itemize}
\end{rema}
\begin{lemm}\label{C[x,y]}
If $F,G\in\C[x,y]$ with $m+n\ge2$, then $p> -\frac{m_0}{m}$ and
$p(G)=p$ $($thus the positions of $F$ and $G$ are symmetric$)$.
\end{lemm}\ni{\it Proof.~}~If $m_0<m$ the result  follows from Lemma \ref{p-ge---}. Thus
suppose $m_0>m$. Note that in order for $[F,G]\in \C\bs\{0\}$, $x$
must appear as a term in $F$ or $G$. If $F$ contain the term $x$,
then $\deg_xf_m\ge1$ and
$$p\ge\frac{\deg_xf_m-\deg_xf_0}{m}\ge-{\frac{m_0}{m}}+\frac{1}{m}>-\frac{m_0}{m}.$$
If $G$ contains the term $x$, then $$p\ge
p(G)\ge\frac{\deg_xg_n-\deg_xg_0}{n}\ge\frac{1-n_0}{n}=\frac{1}{n}-\frac{m_0}{m}>-\frac{m_0}{m},$$
where $n_0=\deg_xg_0=\frac{m_0}{m}n$. 
Thus $p>-\frac{m_0}{m}$. Now letting $r=0$ in
(\ref{G's-component}) shows $G\comp{0}=(F^{\frac{m}{n}})\comp{0}$
which is not a monomial by (\ref{0-th-com}). Thus $p(G)=p$ by Lemma
\ref{lemm2.02}(3).
 \hfill$\Box$\vskip7pt
\begin{lemm}\label{bi-lemm} Suppose $F,G\in\C[x](y)$ with $m+n\ge2$. If $b_i\ne0$ and $i<m+n-1$,
then $m|m'_0(n-i)$. In particular, $m|m'_0n$ and $m|m'_0i$.
\end{lemm}\ni{\it Proof.~}~The lemma is trivial if $m'_0=0$. Suppose $i_0<m+n-1$ is minimal such that
$b_{i_0}\ne0$ and $m\mbox{$\not|$}\,m'_0(n-i_0)$. Multiplying
(\ref{equa2.6}) by $h^{\frac{m'_0(n-i_0)}{m}}$, we see that
$h^{\frac{m'_0(n-i_0)}{m}}$ is a rational function of $x$, which
contradicts Lemma \ref{ration-function} and (\ref{h---}). Thus
$m|m'_0(n-i_0)$. Since $b_0\ne0$ by (\ref{equa2.6}), we have
$m|m'_0n$.
\hfill$\Box$
\begin{lemm}\label{h=x-lemm} Suppose $F,G\in\C[x](y)$ and
$0\!<\!m_0\!\le\! m$ or $m\mbox{$\not|$}\,m_0$.
 Then $h\!=\!x,$ i.e.,\,$d_0\!=\!1,\, m'_0\!=\!m_0$.
\end{lemm}\ni{\it
Proof.~}~{}From (\ref{equa2.6}) and Lemma \ref{bi-lemm}, we see
$b_{m+n-1}$ has the form $b_{m+n-1}=h^{-\frac{m'_0}{m}-a}g$ for some
$a\in\Z$ and $0\ne g\in\C[x]$ (since $b_{m+n-1}\ne0$ by
(\ref{b=m+n-1=})) such that $h\mbox{$\not|$}\,g$. Since $\deg_x\bar
b_0=-\frac{m_0}{m}$, we have $\deg_xb_{m+n-1}=1-\frac{m_0}{m}$ or
$0$. If $\deg_xb_{m+n-1}=1-\frac{m_0}{m}$, then
\begin{equation}\label{derb-m+n-1-bar-b-0}
d_g:=\deg_xg=\deg_xb_{m+n-1}+\big(\frac{m'_0}{m}+a\big)\deg_x h=1+a
d_0 \mbox{ \ (by (\ref{b=m+n-1=}))}.
\end{equation}
If $\deg_xb_{m+n-1}=0$ and $m_0\ne m$ (which can happen only when
$m<m_0$ and $d_g=\frac{m_0}{m}+ad_0$), then $m|m_0$, a
contradiction with our assumption in the lemma.

 Thus we have (\ref{derb-m+n-1-bar-b-0}), which implies either $a\ge0$ or $a=-1,\,d_0=1$. In the
latter case, we have $h=x$ (since $h$ is monic without the constant
term). Thus suppose $a\ge0$. Now (\ref{b=m+n-1=}) shows
\begin{equation}\label{b-m+n-1-bar-b-0}
\a
h^{-a-1}g\ptl_xh+h^{-a}\ptl_xg=h^{\frac{m'_0}{m}}\ptl_yb_{m+n-1}=-\frac{J}{m}h^{-\frac{m'_0}{m}}\bar
b_0=- \frac{J}{m}\mbox{, \ where  \ $\a=-\frac{m'_0}{m}-a$}.
\end{equation}
\def\Aa{\eta}\def\Ab{\l}Factorize $h,g$ as
products of irreducible polynomials of $x$:
\begin{eqnarray}\label{equa-fac1-hg}
h=h_1^{i_1}\cdots h_{\ell}^{i_\ell},\ \ g=g_0h_1^{j_1}\cdots
h_r^{j_r}g_{r+1}^{j_{r+1}}\cdots g_s^{j_s},
\end{eqnarray}
for some $\ell,s,i_1,...,i_\ell,j_1,...,j_s\in\N$ and $0\le r\le
{\rm min}\{s,\ell\}$, $0\ne g_0\in\F$, where
$$g_1\!:=\!h_1,\ \,...,\ \,
g_r\!:=\!h_r,\ \ h_{r+1},\ \,...,\ \ h_\ell,\ \,g_{r+1},\,...,\,g_s
\in\C[x],$$ are different irreducible monic polynomials of $x$
(thus, of degree $1$). Multiplying (\ref{b-m+n-1-bar-b-0}) by
$h^{a+1}$, using (\ref{equa-fac1-hg}), and canceling the common
factor
$$h_1^{i_1+j_1-1}\cdots h_r^{i_r+j_r-1}h_{r+1}^{i_{r+1}-1}\cdots
h_\ell^{i_\ell-1},$$ noting that $\ptl_x h_\Aa=\ptl_x g_\Aa=1$ for
all $\Aa$, we obtain
\begin{eqnarray}\label{equa-fac2-hg}
\!\!\!\!\!\!&\!\!\!\!\!\!& \Big(-(m'_0+am)g_{r+1}\cdots
g_s\mbox{$\sum\limits_{\Aa=1}^\ell$} i_\Aa\frac{h_1\cdots
h_\ell}{h_\Aa} +m h_{r+1}\cdots h_\ell
\mbox{$\sum\limits_{\Ab=1}^s$} j_\Ab\frac{g_1\cdots
g_s}{g_\Ab}\Big)g_{r+1}^{j_{r+1}-1}\cdots g_s^{j_s-1}
\nonumber\\
\!\!\!\!\!\!&\!\!\!\!\!\!& =-J h_1^{i_1a+1-j_1}\cdots
h_r^{i_ra''+1-j_r} h_{r+1}^{i_{r+1}a+1}\cdots h_\ell^{i_\ell a+1}.
\end{eqnarray}
If $\ell>r$, then $h_\ell$ divides all terms except one term
corresponding to $\Aa=\ell$ in (\ref{equa-fac2-hg}), a
contradiction. Thus $\ell=r$. Since $g_{r+1},...,g_s$ do not appear
in the right-hand side of (\ref{equa-fac2-hg}), we must have
$j_{r+1}=...=j_s=1$, and since the left-hand side is a polynomial,
we have \begin{equation}\label{ccc-hg}i_ka+1-j_k\ge0\mbox{ \ \ for \
}k=1,...,\ell.\end{equation} If $i_ka+1-j_k>0$ for some $k$, then
$h_k$ divides all terms except two terms corresponding to $\Aa=k$
and $\Ab=k$ in (\ref{equa-fac2-hg}), and the sum of these two terms
is a term (not divided by $h_k$) with coefficient
$-(m'_0+am)i_k+mj_k$. This proves
\begin{equation}\label{equa-all-k-hg}
i_ka+1-j_k=0\mbox{ \ or \ }-(m'_0+am)i_k+mj_k=0 \mbox{ \ for \
}k=1,...,\ell.
\end{equation}
If $a\ge1$, then either case of (\ref{equa-all-k-hg}) shows $i_k\le
j_k$, and thus, $h|g$, a contradiction with our choice of $g$. Hence
$a=0$.

Now (\ref{derb-m+n-1-bar-b-0}) shows $d_g=1$. Write $g=g_0x+g_1$ for
some $g_1\in\C$. Then (\ref{b-m+n-1-bar-b-0}) gives
\begin{equation*}
\int\frac{d
h}{h}=\frac{m}{m'_0}\big(\frac{J}{m}+g_0\big)\int\frac{dx}{g_0x+g_1},\mbox{ \
thus \ }
\mbox{$h=c(g_0x+g_1)^{\frac{m}{m'_0g_0}(g_0+\frac{J}{m})}$ \ for
some $c\ne0$.}\end{equation*} Since $h\in\C[x]$, we must have
$\b:=\frac{m}{m'_0g_0}(g_0+\frac{J}{m})\in\N$, then (\ref{h---})
shows $\b=1$ (and so $m\!\ne\! m'_0$ since $J\!\ne\!0$), and
$h\!=\!x$ (since $h$ is monic without the constant term). 
\hfill$\Box$ \begin{rema}\label{m0-m}\rm If
$F,G\in\C[x,y]$ and $p\le0$, then by (\ref{p-degree})
(cf.~(\ref{p-line})), we have $\deg_xF=m_0>0$, and by exchanging $x$
and $y$ if necessary, we can sometimes suppose $m_0>m$ and $p<0$
(i.e., $\Coeff(F,x^m)=y^{m_0}$ by Lemma \ref{h=x-lemm}) or sometimes
suppose $m_0<m$ and $p\le0$.
\end{rema}
%
%
\subsection{Jacobi pairs in $\C((x^{-1}))[y]$}
{}From now on we shall suppose $F,G\in\C((x^{-1}))[y]$ with
$m,m+n\ge2$. Then each component $F\comp{r}$ is in $\C[x^{\pm1},y].$
Let $\pri{F}$ be the primary polynomial as in
Definition \ref{comp}(2)(ii) and $\mu$ be as in \eqref{mu-set}.

\begin{lemm}\label{d|n--}
If $\mu<0$ $($so $p<-\frac{m_0}{m}$ by \eqref{mu-set}$)$ 
then $d|n$ $($recall notation
$d$ in $(\ref{denote-d}))$,
 and $p(G)=p$ $($thus
the positions of $F$ and $G$ are symmetric in this case$)$.
\end{lemm}
\ni{\it Proof.~}~Note that if the equality holds in Lemma \ref{p-ge---} then one can obtain $\mu=0$, thus inequality holds in Lemma \ref{p-ge---}.
Setting $r=0$ in (\ref{G's-component}) gives
(cf.~(\ref{0-th-com})) \begin{equation}\label{G0===}
G\comp{0}=(F^{\frac{n}{m}})\comp{0}=x^{\frac{m_0n}{m}}\pri{F}^{\frac{n}{d}} ,\end{equation}
by noting the following facts: \begin{itemize}\parskip-3pt
\item[(i)]  $b_i\!=\!0$ for $1\!\le\! i\!<\!m\!+\!n\!-\!1$ by Lemma \ref{b-i=0-if};
\item[(ii)]  $H\comp{a}\!=\!0$ for $H\!\in\!\C((x^{-1}))[y]$ and all $a\!>\!0$;
\item[(iii)] by \eqref{G's-component-r},
$\a_{0,i}=(p+\frac{m_0}{m})(m+n-1)+\frac{m_0}{m}+\si_i-\eta_i>\frac{m-m_0}{m}
+\frac{m_0}{m}+\si_i-\eta_i=1+\si_i-\eta_i\ge0$ for $i\ge0$, where the part ``\,$>$\,'' is obtained by Lemma \ref{p-ge---}, and
the part ``\,$\ge$\,'' is by Lemma \ref{b'-i=0}.\end{itemize}
Since
$G\comp{0}$ is a polynomial on $y$, we have $d|n$ by (\ref{pri-pro}), \eqref{G0===} and Lemma
\ref{ration-function}. By Lemma
\ref{lemm2.02}(3), $p(G)=p$.\hfill$\Box$

\begin{lemm}\label{d|b}
If $0\le i<\mu$ and $b'_i\ne0$, then $d|(n-i)$.
\end{lemm}
\ni{\it Proof.~}~By Lemma \ref{b'-i=0}, we can suppose
$p>-\frac{m_0}{m}$. Assume $i_0\in\Z_+$ is smallest such that
$i_0<\mu$ with $b'_{i_0}\ne0$ and $d\mbox{$\not|$}{\,}(n-i_0)$. We shall use
(\ref{G's-component}) to compute $G\comp{r}$ for
$
r=-(p+\frac{m_0}{m})i_0.$ 

First suppose $i_0\ge m+n-1$. Then $\mu>m+n-1$. 
Thus we must have $m_0>m$ by \eqref{mu-set}. 
Set $i'_0=i_0-(m+n-1)$.
Then using $i_0<\mu=m+n-\frac{1+p}{\frac{m_0}{m}+p}$, we obtain 
$$i'_0<1-\frac{1+p}{\frac{m_0}{m}+p}=\frac{\frac{m_0}{m}-1}{\frac{m_0}{m}+p},\mbox{ \ so \ }
1+i'_0(\frac{m_0}{m}+p)-\frac{m_0}{m}<0,$$ i.e., $1+\si_{i'_0}<0$
(cf.~(\ref{equa2.7+-+0})). Thus by (\ref{equa2.7+-+}), we have
$\eta_{i'_0}\le0$. Since $b'_{m+n-1+i'_0}\ne0$ (recall notation $b'_i$ in \eqref{denote-b'-i}), we have
$\eta_{i'_0}=0$.
Note from
\eqref{equa2.7+-+0} that when $i\gg0$, we have $1+\si_i>0$. Thus by
(\ref{equa2.7+-+}), $\eta_i\le\si_i+1$, and so by \eqref{G's-component-r},
$$\begin{array}{ll}\a_{r,i}&\!\!\!=
\dis(p+\frac{m_0}{m})(m+n-1-i_0)+\frac{m_0}{m}+\si_i-\eta_i\\[6pt]&\!\!\!=
\dis(p+\frac{m_0}{m})(m+n-1-i_0)+\frac{m_0}{m}+\si_{i'_0}-\eta_{i'_0}
+\si_i-\eta_i-(\si_{i'_0}-\eta_{i'_0})\\[6pt]
&\!\!\!=\dis\si_i-\eta_i-(\si_{i'_0}-\eta_{i'_0})\ge-1-\si_{i'_0}>0\mbox{ \
\ if \ }i\gg0.\end{array}$$
Say the above holds when $i>i_1$.
Then (\ref{G's-component}) gives (using $(F_i)\comp{s}=0$ for all
$s>0$)
\begin{eqnarray}\label{G's-component===}
\!\!\!\!\!\!\!\!\!\!\!\!&\!\!\!\!\!\!\!\!\!\!\!\!&
\mbox{$G\comp{-(p+\frac{m_0}{m})i_0}=\sum\limits_{i=0}^{m+n-2}b_i(F^{\frac{n-i}{m}})
\comp{(i-i_0)(p+\frac{m_0}{m})}+b'_{i_0}(F^{\frac{n-i_0}{m}})\comp{0}
+\sum\limits_{^{^{\sc\ \ 0\le
i<i_1}}_{_{\sc i\ne
i'_0,\,\a_{r,i}\le0}}}(F_{m+n-1+i})\comp{\a_{r,i}}
.$}
\end{eqnarray}
Let $a\ge0$ be any rational number. Similar to \eqref{G's-component===},
we also have
\begin{eqnarray}\label{G's-component===1}
\!\!\!\!\!\!\!\!\!\!\!\!&\!\!\!\!\!\!\!\!\!\!\!\!&
\mbox{$G\comp{-(p+\frac{m_0}{m})i_0+a}=\sum\limits_{i=0}^{m+n-2}
b_i(F^{\frac{n-i}{m}})\comp{(i-i_0)(p+\frac{m_0}{m})+a}
+\sum\limits_{0\le i\le i_1,\,\a_i=-a}b'_{m+n-1+i}(F^{\frac{1-m-i}{m}})\comp{0}
$}\nonumber\\ \!\!\!\!\!\!\!\!\!\!\!\!&\!\!\!\!\!\!\!\!\!\!\!\!&
\mbox{$\phantom{G\comp{-(p+\frac{m_0}{m})i_0+a}=}
+\sum\limits_{0\le
i\le i_1,\,\a_{r,i}+a<0}(F_{m+n-1+i})\comp{\a_{r,i}+a}
.$}
\end{eqnarray}
Note 
that
$G\comp{-(p+\frac{m_0}{m})i_0+a}$ is a polynomial on $y$.
Set $$A:=\max\{-(i-i_0)(p+\frac{m_0}{m}),\,-\a_{r,j}\,|\,0\le i\le m+n-2,\,0\le j\le i_1,\a_{r,j}\le0\}.$$
If we first take $a=A$ in \eqref{G's-component===1}, then
the last summand vanishes and the first summand is summed over those $i$'s such that
$(i-i_0)(p+\frac{m_0}{m})+a=0$, so we can use
%
Lemma \ref{ration-function} to obtain that all nonzero terms
in  \eqref{G's-component===1} are rational.
Now by taking $a<A$ (and $a\ge0$) in \eqref{G's-component===1}
 and by induction on $A-a$ (note that there are only finitely many $a$'s
with $0\le a\le A$ such that there exist nonzero terms in \eqref{G's-component===1}), we can see that all
nonzero terms in \eqref{G's-component===1} (thus in \eqref{G's-component===}) are rational. In particular, $\pri{F}^{\frac{n-i_0}{d}}
=x^{-\frac{m_0(n-i_0)}{m}}(F^{\frac{n-i_0}{m}})\comp0$ (cf.~\eqref{0-th-com}) is rational,
a contradiction with \eqref{pri-pro}.

Now suppose $0\le
i_0<m+n-1$. 
As above,
noting from
\eqref{equa2.7+-+0} that when $i\gg0$, we have $1+\si_i>0$ and
$\si_i-\eta_i\ge-1$ by (\ref{equa2.7+-+}). Thus
$$\a_{r,i}=r+(p+\frac{m_0}{m})(m+n-1)+\frac{m_0}{m}+\si_i-\eta_i\ge
(p+\frac{m_0}{m})(m+n-1-i_0)+\frac{m_0}{m}-1>0,$$
where the part ``$>$'' is obtained by using $i_0<\mu=m+n-\frac{1+p}{\frac{m_0}{m}+p}$.
As in \eqref{G's-component===}, say the above
inequality holds for $i>i_1$. Then
(\ref{G's-component}) shows (in this case $b'_{i_0}=b_{i_0}$ by \eqref{denote-b'-i})
\begin{equation*}
\mbox{$G\comp{-(p+\frac{m_0}{m})i_0}=\sum\limits_{0\le
i<i_0}b_i(F^{\frac{n-i}{m}})\comp{(i-i_0)(p+\frac{m_0}{m})}
+b_{i_0}(F^{\frac{n-i_0}{m}})\comp{0}+\sum\limits_{0\le i\le i_1,\,\a_{r,i}\le0}(F_{m+n-1+i})\comp
{\a_{r,i}}.$}\end{equation*}
Thus as in \eqref{G's-component===}, we obtain a contradiction.
%
%
%
%
\hfill$\Box$\vskip7pt

\def\tildeb{\tilde b_\mu}We always set\begin{equation}\label{tildeB}
\tildeb=\left\{\begin{array}{ll}
b_\mu&\mbox{if \ }\mu\in\Z_+\mbox{ and }\mu<m+n-1,\\[4pt]
0&\mbox{otherwise}.\end{array}\right.\end{equation}
 Let (cf.~notation $\si_i$ in
(\ref{equa2.7+-+}))
\begin{eqnarray}\label{l-r0}\!\!\!\!\!\!\!\!\!\!\!\!\!\!\!
R_0\!\!\!&:=\!\!\!&\mbox{$\tildeb
(F^{\frac{n-\mu}{m}})\comp{0}+\sum\limits_{i=m+n-1}^\infty {\Coeff}(b_{i},x^{1+\aaa_{i-(m+n-1)}})x^{1+\aaa_{i-(m+n-1)}}(F^{\frac{n-i}{m}})\comp{0}$}
\nonumber\\
\!\!\!\!\!\!\!\!\!\!\!\!\!\!\!&=\!\!\!&\mbox{$ \tildeb
x^{\frac{m_0(n-\mu)}{m}}\pri{F}^{\frac{n-\mu}{d}}-{\dis\frac{J}{m}}
\sum\limits_{i=0}^\infty\dis\frac{1\!-\!i}{1\!+\!\aaa_i}\bar
b_{i,0}x^{1-m_0+ip}\pri{F}^{\frac{1-i-m}{d}}$, where
$\bar b_{i,0}\!=\!\Coeff(\bar b_i,x^{\aaa_i})\!\in\!\F$,}\end{eqnarray}  where the last
equality follows from (\ref{0-th-com}) and Lemma \ref{lemm2.3} (if
$1+\si_i=0$, i.e., $i=\mu-(m+n-1)\ge0$, then
 $\tilde b_{\mu}=0$ by \eqref{tildeB}, in this case
we use the convention that
$-\frac{J(1-i)}{m(1+\aaa_i)}\bar b_{i,0}$ is regarded as the limit
$\lim_{i\to \mu-(m+n-1)}-\frac{J(1-i)}{m(1+\aaa_i)}\bar b_{i,0}$ which is
{\bf defined} to be $b'_\mu$, cf.~\eqref{denote-b'-i}).
%
We claim \VS{-7pt}that\begin{equation}\label{R-0not=0} R_0\ne0,\mbox{ \ \ and
\ }\deg_yR_0=\frac{1+p}{\frac{m_0}{m}+p}-m,\mbox{ \ or \
}1-m.\end{equation} This is because: if $\tildeb\ne0$ (then
$\mu\in\Z_+$) then $\Coeff(R_0,y^{n-\mu})=\tildeb\ne0$; if
$\tildeb=0$ then $\Coeff(R_0,y^{1-m})=\frac{1}{1+\aaa_0}\bar
b_{0,0}x^{1-m_0}\ne0$ ($\bar
b_{0,0}=\Coeff(f_0,x^{m_0})=1$ by \eqref{equa2.7+}).
\begin{lemm}\label{R===0} We have
\begin{equation}\label{ra-L0}
\mbox{$R_0=G\comp{-\mu (p+\frac{m_0}{m})}-\d
b_0(F^{\frac{n}{m}})\comp{-\mu (p+\frac{m_0}{m})}-\sum\limits_{0\le
i<\mu}b'_i(F^{\frac{n-i}{m}})\comp{(i-\mu)(p+\frac{m_0}{m})}=\pri{F}^{-a-m'}P
$},\end{equation}  for some $a\in\Z,\,P\in\C[x^{\pm1}][y]$, such
that $\pri{F}\mbox{$\not|$}\,P$, where $\d=1$ if $\mu<0$ or $0$
otherwise.\end{lemm}\ni{\it Proof.~}~The last equality follows from
Lemmas \ref{lemm-product-KH}(5), 
\ref{bi-lemm} and \ref{d|b}. To prove the
first equality of (\ref{ra-L0}), set $r=-\mu(\frac{m_0}{m}+p)$ in
(\ref{G's-component}), we \VS{-7pt}obtain
\begin{equation}\label{TO-pro}\mbox{$
G\comp{-\mu(\frac{m_0}{m}+p)}=\sum\limits_{j=0}^{m+n-2}b_j(F^{\frac{n-j}{m}})\comp{(j-\mu)(p+\frac{m_0}{m})}
+\sum\limits_{i=0}^\infty(F_{m+n-1+i})\comp{1+\si_i-\eta_i}\VS{-7pt}.
$}\end{equation} First assume $p>-\frac{m_0}{m}$. Compare the $i$-th
term of (\ref{TO-pro}) with corresponding terms of (\ref{l-r0}) and
(\ref{ra-L0}):\begin{itemize}\parskip-3pt\item[(1)] If $1+\si_i\ge0$
(so $m+n-1+i\ge\mu$ and (\ref{ra-L0}) does not have such a term), then
either \begin{itemize}\parskip-1pt\item[(i)] $\eta_i=1+\si_i$: the $i$-th term
of (\ref{TO-pro}) corresponds to the $i$-th term of (\ref{l-r0}), or
\item[(ii)] $\eta_i<1+\si_i$: the $i$-th terms of
(\ref{TO-pro}) and (\ref{l-r0}) are both
zero.\end{itemize}\item[(2)] If $1+\si_i<0$ (so $m+n-1+i<\mu$ and
(\ref{l-r0}) does not have such a term), then either
\begin{itemize}\parskip-1pt\item[(i)] $\eta_i=0$:  the $i$-th term of (\ref{TO-pro})
corresponds to the $i$-th term of (\ref{ra-L0}), or
\item[(ii)] $\eta_i<0$: the $i$-th terms of (\ref{TO-pro})
and (\ref{ra-L0}) are both zero.\end{itemize}\end{itemize} This
proves the lemma in this case. Assume $p<-\frac{m_0}{m}$. Then by
Lemmas \ref{b'-i=0} and \ref{d|n--}, the summand in (\ref{ra-L0}) is
empty and the first summand in (\ref{TO-pro}) has only one term
corresponding to $j=0$. We again have the lemma.
\hfill$\Box$\vskip7pt

By (\ref{R-0not=0}) and (\ref{ra-L0})\VS{-7pt},
\begin{equation}\label{degree-P}
d_P:=\deg_yP=da+\frac{1+p}{\frac{m_0}{m}+p}\mbox{ (in case \
$\tildeb\ne0$) \ \ or \ \ }da+1\mbox{ (in case \ $\tildeb=0$)}\VS{-7pt}.
\end{equation}
Computing the zero-th component of (\ref{equa2.7}), using
(\ref{0-th-com}), similar as in (\ref{G's-component}), we obtain
\begin{equation*}
y=\mbox{$\sum\limits_{i=0}^\infty(\bar
b_iF^{\frac{1-i}{m}})\comp{0}=\sum\limits_{i=0}^\infty{\rm
Coeff}(\bar
b_i,x^{\aaa_i})x^{\aaa_i}(F^{\frac{1-i}{m}})\comp{0}=\sum\limits_{i=0}^\infty\bar
b_{i,0}x^{ip}\pri{F}^{\frac{1-i}{d}}.$}
\end{equation*}
Taking $\ptl_y$ \VS{-7pt}gives
\begin{equation}\label{ptl-zero-th-com-y}
1=\mbox{$\sum\limits_{i=0}^\infty\dis \frac{1-i}{d}\bar
b_{i,0}x^{ip}\pri{F}^{\frac{1-i}{d}-1}\ptl_y\pri{F}.$}
\end{equation}
\begin{lemm}\label{quasi-J}$(F\comp{0},R_0)=(x^{m_0}\pri{F}^{m'},\pri{F}^{-a-m'}P)$
is a \QJ  pair.\end{lemm}\ni{\it Proof.~}~First note that
$R_0$ is a $p$-type q.h.e., thus $p(R_0)=p$ (or $-\infty$ if it is a
monomial). By (\ref{ra-L0}), $R_0\in\C(x,y)$ has the form
(\ref{monic-F-G}). {}From (\ref{l-r0}) and 
(\ref{ptl-zero-th-com-y}), we see
$[F\comp{0},R_0]=J\in\C\bs\{0\}$ (cf.~proof of Lemma
\ref{lemm2.3}), which can be also proved as follows:
Denote $G_1 = \sum_{i=0}^\infty b'_i F^{\frac{n-i}{m}}$
and $H = G - G_1.$
Then $[F,H] = J$. Note that $[F\comp0,H\comp0]$ has the highest term $1$ (up to a scalar) since $F$ has the highest term $x^{m_0}y^m$ and
$H$ has the highest term $x^{1-m_0} y^{1-m}$.
Thus $[F_{[0]},H_{[0]}]$ is a $p$-type q.h.e. whose support lies
in a line passing the point $(0,0)$.
Comparing the $0$-th components in $[F,H] = J$,
we obtain $ [F_{[0]},H_{[0]}]=J.$
Since $H_{[0]}$ and $R_0$ only differ by some $\C$-combinations of rational powers of $F$,
we have $[F_{[0]},R_0] =[F_{[0]},H_{[0]}] = J.$
\hfill$\Box$\vskip7pt

 Multiplying (\ref{l-r0}) by
$\pri{F}^{\frac{\mu-n}{d}}$, taking $\ptl_y$ and using
(\ref{ra-L0}), (\ref{mu-set}) and (\ref{ptl-zero-th-com-y}), we have (cf.~the remark after \eqref{asas})
\begin{eqnarray}\label{diff-1}
\!\!\!\!\!\!\!\!\!\!\!\!\!\!\!\!\!\!\!\!&
\!\!\!\!\!\!\!\!\!\!\!\!\!\!\!\!\!\!\!\!&
\big(-a-m'+\frac{\mu-n}{d}\big)\pri{F}^{-a-m'+\frac{\mu-n}{d}-1}P\ptl_y\pri{F}+\pri{F}^{-a-m'+\frac{\mu-n}{d}}\ptl_yP
\nonumber\\
\!\!\!\!\!\!\!\!\!\!\!\!\!\!\!\!\!\!\!\!&
\!\!\!\!\!\!\!\!\!\!\!\!\!\!\!\!\!\!\!\!&
=\ptl_y(R_0\pri{F}^{\frac{\mu-n}{d}})=-\frac{J}{m}\mbox{$\sum\limits_{i=0}^\infty$}{\dis\frac{1-i}{1+\aaa_i}}\bar
b_{i,0}x^{1-m_0+ip}\ptl_y\big(\pri{F}^{-\frac{1+\aaa_i}{d(p+\frac{m_0}{m})}}\big)\nonumber\\
\!\!\!\!\!\!\!\!\!\!\!\!\!\!\!\!\!\!\!\!&
\!\!\!\!\!\!\!\!\!\!\!\!\!\!\!\!\!\!\!\!&={\dis\frac{Jx^{1-m_0}}{mp+m_0}}
\mbox{$\sum\limits_{i=0}^\infty$}\frac{1-i}{d}\bar
b_{i,0}x^{ip}\pri{F}^{\frac{1-i}{d}-1}\ptl_y\pri{F}\cdot\pri{F}^{-\frac{1+p}{d(p+\frac{m_0}{m})}}
 =
{\dis\frac{Jx^{1-m_0}}{mp+m_0}}\pri{F}^{-\frac{1+p}{d(p+\frac{m_0}{m})}}.\end{eqnarray}
We write $p=\frac{p'}{q}$ for some coprime   integers $p',\,q$ such
that $q>0$. Note that $\pri{F}$ being a monic $p$-type q.h.e., has
the form
\begin{equation}\label{primary-F}\mbox{$\pri{F}=y^d+\sum\limits_{i=1}^d \pri c_ix^{ip}y^{d-i}\in\C[x^{\pm1}][y]$ \
for some \ $\pri c_i\in\F$}.\end{equation} Noting that
$x^{m_0}\pri{F}^{m'}=F\comp{0}\ne x^{m_0}y^m$ since $p\ne-\infty$, we have $\pri c_i\ne0$ for some $i$. Hence at least one of
$p,2p,...,dp$ is an integer. Thus
\begin{equation}\label{q=}1\le q\le d.\end{equation}
Multiplying (\ref{diff-1}) by
$\pri{F}^{\frac{1+p}{d(p+\frac{m_0}{m})}}$, we obtain the following
differential equation on $\pri{F}$ and $P$,
\begin{equation}\label{equa-int1}
\a_1\pri{F}^{-a-1}P\ptl_y\pri{F} +\a_2\pri{F}^{-a}\ptl_yP=\a_3,
\mbox{ \ or \ } \a_1P\ptl_y\pri{F}
+\a_2\pri{F}\ptl_yP=\a_3\pri{F}^{a+1},
\end{equation}
where
\begin{equation}\label{asas}
\a_1=-a\a_2-m'(p'+q),\, \ \ \ \a_2=p'm+m_0q,\, \ \ \
\a_3=
%
qJx^{1-m_0}.
\end{equation}
We remark that the only purpose of \eqref{diff-1} is to prove \eqref{equa-int1}, which can be also directly proved by the following formal arguments:
Noting that $\pri F,\,P$ are $p$-type q.h.e., so
 $x\ptl_x\pri F,\,x\ptl_x\pri P$ are combinations of $F,\,y\ptl_y\pri F$ and
$P,\,y\ptl_y\pri P$,
using this in $[x^{m_0}{\pri F}^{m'},{\pri F}^{-a-m'}P]\!=\!J$, one can easily deduce
\eqref{equa-int1} for some $\a_1,\a_2\!\in\!\C$ and some $0\!\ne\!\a_3\!\in\!\C[x^{\pm1}]$.
Note that the first term of $R_0$ in \eqref{l-r0} does not contribute to
\eqref{equa-int1}, and whether or not $P$ is a polynomial on $y$ does not affect \eqref{equa-int1} either, thus for the purpose of proving \eqref{equa-int1} with $\a_i$ satisfying
\eqref{asas},
 if necessary
by changing the coefficient $\tilde b_\mu$, we may suppose that the $y$-degree of $P$ is $\deg_yP\!=\!ad\!+\!\frac{1+p}{\frac{m_0}{m}+p}\!=\!ad\!+\!\frac{(p'+q)m}{m_0q+p'm}$\VS{-4pt}. Then by
comparing the coefficients of the highest degree of $y$ in \eqref{equa-int1}, we immediately obtain
$\a_1 d\!+\!\a_2(ad\!+\!\frac{(p'+q)m}{m_0q+p'm})\!=\!0.$
Thus up to a scalar, we can suppose $\a_1,\a_2$ have the forms as in\,\eqref{asas}\,(what is $\a_3$ is not important since later on we only need the fact that $\a_3\!\ne\!0$).
\begin{lemm}\label{Factor-q}
\begin{itemize}\parskip-3pt\item[{\rm(1)}] For every irreducible factor $Q\in\C[x^{\pm1}][y]$ of $\pri{F}$
or $P$, we have $\deg_yQ=1$ or $q|{\ssc\,}\deg_yQ$. If $Q$ is monic and
$q\mbox{$\not|$}\,\deg_yQ$, then $Q=y$.
\item[{\rm(2)}] If $d=q>1$, then $\pri{F}$ is
irreducible and $\pri{F}=y^q+x^{p'}$ $($up to rescaling  $x)$.
\end{itemize}
\end{lemm}\ni{\it Proof.}~~(1) Note that $R_0$ is a $p$-type q.h.e., and $\pri{F}$ is a $p$-type q.h.e.~(cf.~the right-hand side
of (\ref{primary-F})). By Lemma \ref{lemm-product-KH}(4), $P=R_0\pri{F}^{a+m'}$ must be a $p$-type
q.h.e. By Lemma \ref{lemm-product-KH}(4) again,
every irreducible factor $Q$ of $\pri{F}$ or $P$ must be q.h.e.~of
the form $\sum_{i=0}^{\g}u_ix^{pi}y^{\g-i}$, where $\g=\deg_yQ$ and
$u_i\in\C$. If $q\mbox{$\not|$}\,\g$, then $p\g=\frac{p'\g}{q}$ cannot
be an integer, thus $u_{\g}=0$, and $Q$ contains the factor $y$.

(2) It follows from (1) and (\ref{pri-pro}) (cf.~(\ref{primary-F})).
\hfill$\Box$ \vskip 7pt

\subsection{Some examples of \QJ  pairs with $p\le0$}


\begin{exam}\label{1ex3.2.1}\rm
Below we obtain some \QJ  pairs starting from a simple one
by adding some powers of $F$ (or $G$) to $G$ (or $F$)\VS{-3pt}.
\begin{itemize}\parskip-3pt\item[(1)]
$(xy^2,y^{-1})\to(xy^2,xy^2+y^{-1})\to(x^2y^4+xy^2+2xy+y^{-2},xy^2+y^{-1}):=(F,G).
$ In this case,
$p=-\frac{1}{3},\,m_0=2,\,m=4,\,F\comp{0}=(xy^2+y^{-1})^2$ and
$n=2$.
\item[(2)]
$(x^{\frac{1}{2}}y,x^{\frac{1}{2}})\to(x^{\frac{1}{2}}y,xy^2+x^{\frac{1}{2}})\to
(x^2y^4+2x^{\frac{3}{2}}y^2+x+x^{\frac{1}{2}}y,xy^2+x^{\frac{1}{2}}):=(F,G)$.
In this case,
$p=-\frac{1}{4},\,m_0=2,\,m=4,\,F\comp{0}=(xy^2+x^{\frac{1}{2}})^2$
and $n=2$.
\end{itemize}
\end{exam}\vskip-5pt

In Example \ref{1ex3.2.1}, we have $m|n$ or $n|m$. If this holds in
general, then the two-dimensional Jacobian conjecture can be proved by induction on $m$.
However, the following example shows that this may not be the case. 
\vskip-22pt\

\begin{exam}\label{2ex3.2.1}\rm
\begin{itemize}\parskip-3pt\item[(1)]
$F=xy^2+2x^{\frac{5}{8}}y=F\comp{0},\,p=-\frac{3}{8},\,m_0=1,\,m=2,$
and $
G=x^{\frac{3}{2}}y^3+3x^{\frac{9}{8}}y^2+\frac{3}{2}x^{\frac{3}{4}}y
-\frac{1}{2}x^{\frac{3}{8}}=R_0,\,n=3$.
\item[(2)]
$F=x^2y^{10}+2xy^4=F\comp{0},\,p=-\frac{1}6,\,m_0=2,\,m=10$, and
$G=x^3y^{15}+3x^2y^9+\frac{1}2xy^3-\frac{1}2y^{-3}=R_0,\,n=15$.
\end{itemize}
\end{exam}\vskip-5pt

{}For the above examples, we observe that if $m_0<m$ then either
$F,G$ must contain some negative power of $y$ or contain some
non-integral power of $x$.

In case $F,G\in\C[x,y]$, by Remark \ref{m0-m}, we can suppose $m_0<m$ if $p\le0$. To better understand
Jacobi pairs, we first
 suppose $m_0>m$ and obtain the following two lemmas.
 We rewrite $F,G$ as (where $f_{00}=1$\VS{-9pt})
\begin{eqnarray}
\label{F____}\!\!\!\!\!\!\!&&F=\mbox{$\sum\limits_{i=0}^mf_iy^{m-i}=\sum\limits_{(m_0-i,m-j)\in\Supp\,F}f_{ij}x^{m_0-i}y^{m-j},$}\\
\label{G____}\!\!\!\!\!\!\!&&G=\mbox{$G_1+R,\ \
G_1=\sum\limits_{0\le i<\mu}b'_iF^{\frac{n-i}{m}},\ \
R=\sum\limits_{\min(m+n-1,\mu)\le j<\infty}b''_jF^{\frac{n-j}{m}},$}
\end{eqnarray}
where $b'_i$ is defined in (\ref{denote-b'-i}), and $b''_j=b_j-b'_j$
if $j<\mu$ and $b''_j=b_j$ if $j\ge\mu$ (recall notation $\mu$ in \eqref{mu-set}). Then $[F,G]=[F,R]$.
\begin{theo}\label{C[x,y]'}\begin{itemize}\parskip-3pt\item[\rm(1)]
Suppose $F,G$ is a Jacobi pair in $\C[x,y]$. Then $p\ne-1$.
\item[\rm(2)]
Suppose $F,G\in\C[x^{\pm1},y]$ with $m_0>m$ and
$p<0$. Then $p\le-1$. Furthermore, if $p=-1$ then $\res_x(F\ptl_xG)\ne0$ and $F\comp0=x^{m_0}(y-\b x^{-1})^m$ for some $\b\in\C\bs\{0\}$. Thus by replacing $y$ by $y+\b x^{-1}$ for some $\b\in\C$ if necessary, we can suppose $p<-1$.
\end{itemize}\end{theo}
\ni{\it Proof.~}~(1) If $m_0<m$, then $p>-\frac{m_0}{m}>-1$ by Lemma \ref{p-ge---}. Thus assume $m_0>m>0$. Suppose $p=-1$. Then $\mu=m+n$ and $\tildeb=0$ (cf.~\eqref{tildeB}). So (\ref{R-0not=0}) shows
$\deg_yR_0=1-m$, and (\ref{l-r0}) in fact shows $R_0=R\comp{0}$.
Thus we can \VS{-7pt}write \equa{R===AA}{\mbox{$R=\sum\limits_{i=0}^\infty
r_iy^{1-m-i}=\sum\limits_{(1-m_0-i,1-m-j)\in\Supp\,R}r_{ij}x^{1-m_0-i}y^{1-m-j},\ \ r_i\in\C((x^{-1})),\ r_{i,j}\in\C.$}}
By Lemma \ref{clai-b+m+n} and (\ref{G____}), we see $F^{-1}$ does
not appear in $R$, thus we can \VS{-7pt}write
\equan{eqno(1)}{R=r'_0 F^{\frac{1-m}{m}}+\sum\limits_{i=1}^\infty r'_i F^{\frac{-m-i}{m}}\mbox{
for some $r'_i\in\AA$.}}
Since $F^{\frac{1-m}{m}}=f_0^{\frac{1-m}{m}}(y^{1-m}+\frac{1-m}{m}f_0^{-1}f_1y^{-m}+\cdots)$,
comparing the coefficients of $y^{1-m}$ and $y^{-m}$ in the above equation on $R$,
we obtain $r_0=r'_0f_0^{\frac{1-m}{m}}$ \VS{-7pt}and \equa{R1==}{\dis r_1=\frac{1-m}{m}r'_0f_0^{\frac{1-m}{m}}f_0^{-1}f_1
=\frac{1-m}{m}r_0f_0^{-1}f_1.}
 If necessary, by
replacing $x$ by $x+\b$ for some $\b\in\C$ (which does not affect
$p$), we can suppose $f_{10}=0$. \VS{-5pt}Thus
\equa{F1001=0}{\mbox{$0=f_{10}=f_{01}$}\VS{-5pt},} where the second equality
follows from $p<0$. Computing the coefficients of $x^{-1}y^0$ and $x^0y^{-1}$
in $J=[F,R]$, using (\ref{G____}), (\ref{R===AA}) and
(\ref{F1001=0}), we obtain $r_{10}=r_{01}=0$.
This and (\ref{R1==}) \VS{-7pt}give
\equa{r_11=0}{r_{11}=\dis\frac{1-m}{m}r_{00}f_{11}.} Note that we obtain
(\ref{r_11=0}) only by (\ref{F1001=0}) (and is independent of $p$).
Thus by symmetry (exchanging $x$ and $y$), we also have $r_{11}=\frac{1-m_0}{m_0}r_{00}f_{11}$.
\VS{-7pt}Therefore  \equa{rrr=0}{r_{11}=f_{11}=0,}
which can be also obtained by using $\Res_x(F\ptl_xG)=\Res_y(F\ptl_yG)=0$.
 Since
$F\comp{0}$ and $R_0$ are $p$-type q.h.e with $p=-1$, of the form $F=\sum_{i=0}^m f_{ii}x^{m_0-i}y^{m-i}$ and $R_0=\sum_{i=0}^\infty
r_{ii}x^{1-m-i}y^{1-m-i}$, we can easily deduce
\equan{F-R-0}{
x\ptl_xF\comp{0}-y\ptl_yF\comp{0}=(m_0-m)F\comp{0},\,\ \ \
x\ptl_xR_0-y\ptl_yR_0=-(m_0-m)R_0.}
Using this in
$J=(\ptl_xF\comp{0})(\ptl_yR_0)-(\ptl_yF\comp{0})(\ptl_xR_0)$, we
obtain
$
{xJ=(m_0-m)\ptl_y(F\comp{0}R_0).}
$ For convenience,
we can suppose $J=m_0-m$. Therefore, \equa{SAAA}{F\comp{0}R_0=xy+a_0\mbox{ for some $a_0\in\AA$. }}
Since $F\comp{0}R_0$ is a $p$-type q.h.e., we have $a_0\!\in\!\C$. By
(\ref{rrr=0}), $a_0\!=\!0$. Using (\ref{G____}), we have (here and below,
all equalities associated with integration mean that they hold up to
some elements in $\AA$\VS{-7pt}) \equan{integr}{{\dis\int}
G(\ptl_yF)dy=\sum\limits_{i=0}^{m+n-1}{\dis\frac{m}{m+n-i}}b'_iF^{\frac{m+n-i}{m}}+Q,\mbox{
\ where \ }Q={\dis\int} R(\ptl_yF)dy.} As in (\ref{l-r0}) and
(\ref{ra-L0}),  this shows $Q\comp{0}=\int
R_0(\ptl_yF\comp{0})dy\in\AA(y)$ by Lemma \ref{d|b}. Factorize
$F\comp{0}=a'_0\prod_{i=1}^k(y-a_i)^{\l_i}$ in $\AA(y)$, where
$a_i\in\AA,\,a'_0\ne0,\,\l_i>0$ and $a_1,...,a_k$ are distinct. Then
(\ref{rrr=0}) and the fact that $F\comp{0}$ is not a monomial show
$k>1$ and some $a_i\ne0$. Then (\ref{SAAA}) \VS{-7pt}gives
\equan{NOT=ra}{Q\comp{0}=\sum\limits_{i=1}^k\l_ix{\dis\int\frac{ydy}{y-a_i}}
=\sum\limits_{i=1}^k\l_ix\Big(y+a_i\,{\rm ln}(y-a_i)\Big),} which
cannot be a rational function on $y$. A contradiction.

(2) Consider (\ref{equa-int1}) and (\ref{asas}). If $d=1$ then (\ref{primary-F})
shows $p\in\Z$. Thus $p\le-1$. So suppose $d\ge2$. If $a\le-2$, then
$\pri{F}$ divides the left-hand side of the first equation of
(\ref{equa-int1}), a contradiction.
 Thus $a\ge-1$.

\def\hp{h}If $a=-1$, then the second case of (\ref{degree-P}) cannot occur
since $d_P\ge0$ and $d\ge2$. Thus $\tildeb\ne0$, i.e., $\mu<m+n-1$,
which implies $m_0<m$, a contradiction with our assumption.

Now suppose $a\ge0$. Similar to the proof of Lemma \ref{h=x-lemm},
we factorize $\pri{F},P$ as products of irreducible polynomials of
$y$ in the ring $\AA[y]$\VS{-7pt}: 
\begin{eqnarray}\label{equa-fac1===}
\pri{F}=\pri{f}_1^{i_1}\cdots \pri{f}_{\ell}^{i_\ell},\ \
P=\hp_0\pri{f}_1^{j_1}\cdots
\pri{f}_r^{j_r}\hp_{r+1}^{j_{r+1}}\cdots \hp_s^{j_s},
\end{eqnarray}
for some $\ell,s,i_1,...,i_\ell,j_1,...,j_s\in\N$ and $0\le r\le
{\rm min}\{s,\ell\}$, $0\ne \hp_0\in\AA$, 
\VS{-7pt}where
$$\hp_1\!:=\!\pri{f}_1,\ \,...,\ \,
\hp_r\!:=\!\pri{f}_r,\ \ \pri{f}_{r+1},\ \,...,\ \ \pri{f}_\ell,\
\,\hp_{r+1},\,...,\,\hp_s \in\AA[y]
\VS{-4pt},$$ are different
irreducible monic polynomials of $y$  of degree $1$. As in
(\ref{equa-fac2-hg}), we obtain (where $\a_1,\a_2,\a_3$ are as in
(\ref{asas})\VS{-7pt})
\begin{eqnarray}\label{equa-fac2}
\!\!\!\!\!\!&\!\!\!\!\!\!& h_0\Big(\a_1\hp_{r+1}\cdots
\hp_s\mbox{$\sum\limits_{\Aa=1}^\ell$} i_\Aa\frac{\pri{f}_1\cdots
\pri{f}_\ell}{\pri{f}_\Aa} +\a_2 \pri{f}_{r+1}\cdots \pri{f}_\ell
\mbox{$\sum\limits_{\Ab=1}^s$} j_\Ab\frac{\hp_1\cdots
\hp_s}{\hp_\Ab}\Big)\hp_{r+1}^{j_{r+1}-1}\cdots \hp_s^{j_s-1}
\nonumber\\
\!\!\!\!\!\!&\!\!\!\!\!\!& =\a_3 \pri{f}_1^{i_1a+1-j_1}\cdots
\pri{f}_r^{i_ra+1-j_r} \pri{f}_{r+1}^{i_{r+1}a+1}\cdots
\pri{f}_\ell^{i_\ell a+1}.
\end{eqnarray}
Similar to the arguments after (\ref{equa-fac2-hg}), we have
$r=\ell$ \VS{-5pt}and
\begin{equation}\label{j-k=}j_{r+1}=...=j_s=1.\end{equation} Also for
all $k\le\ell$\VS{-7pt},
\begin{eqnarray}\label{ccc}
\!\!\!\! \!\!\!\! \!\!\!\! &\!\!\!\! \!\!\!\! \!\!\!\!
&i_ka+1-j_k\ge0, \mbox{ and}\\
\label{equa-all-k} \!\!\!\! \!\!\!\! \!\!\!\! &\!\!\!\! \!\!\!\!
\!\!\!\! & i_ka+1-j_k=0\mbox{ \ or \ }\a_1i_k+\a_2j_k=0.
\end{eqnarray}
First assume $a=0$. Then (\ref{ccc}) shows $j_k=1$ for all $k$. Thus
$s=d_P$, and (\ref{equa-fac2}) is simplified \VS{-9pt}to
\begin{equation}\label{equa-all-k+} \hp_{\ell+1}\cdots
\hp_s\mbox{$\sum\limits_{\Aa=1}^\ell$}(\a_1i_\Aa+\a_2)\frac{\pri{f}_1\cdots
\pri{f}_\ell}{\pri{f}_\Aa} +\a_2\pri{f}_1\cdots
\pri{f}_\ell\mbox{$\sum\limits_{\Ab=\ell+1}^s$}\frac{\hp_{\ell+1}\cdots
\hp_s}{\hp_\Ab}=\a_3\VS{-9pt}.
\end{equation}
If $\ell=1$,  then (\ref{equa-fac1===}) shows that $i_1=d$ and
$\pri{F}=\pri{f}_1^{d}$. Write $\pri{f}_1=y+\pri{f}_{11}$ for some
$0\ne\pri{f}_{11}\in\AA$ 
since $\pri{F}\ne y^d$. Comparing
the coefficients of $y^{d-1}$ in both sides of (\ref{primary-F})
gives $d\pri{f}_{11}=\pri c_1x^p$, we see $\pri c_1\ne0$, thus $p\in\Z$ (since
$\pri{F}\in\C[x^{\pm1},y]$), i.e., $p\!\le\!-1$. Thus suppose $\ell\!\ge\!2$.
Comparing the coefficients of the terms with highest $y$-degree (i.e.,
degree $s\!-\!1$, which is $\ge\!\ell\!-\!1\!\ge\!1$) in (\ref{equa-all-k+}) shows (using (\ref{asas}) and
$a\!=\!0$\VS{-11pt})
\begin{equation*}
0=\mbox{$\sum\limits_{\Aa=1}^\ell(\a_1i_\Aa+\a_2)
+\sum\limits_{\Ab=\ell+1}^s\a_2$}=\a_1d+s\a_2
=-m(p'+q)+s(p'm+m_0q)\VS{-5pt},
\end{equation*} i.e.,
$d_P=s=\frac{m(p'+q)}{m_0q+p'm}$, and the first case of
(\ref{degree-P}) occurs, a contradiction with $m_0>m$.

Now suppose $a\ge1$. If $p'+q>0$, then either case of (\ref{equa-all-k})
shows $i_k<j_k$ for all $k$ (note from (\ref{asas}) that
$-\a_1>\a_2$), i.e., $\pri{F}|P$, a contradiction with our choice of
$P$. Thus $p'+q\le0$, i.e., $p=\frac{p'}{q}\le-1$.
Now assume $p=-1.$ If necessary, by replacing $y$ by $y+\b x^{-1}$ for some $\b \in\C$, we can suppose
$f_{11}=0$ (cf.~notations in \eqref{F____}). Note that we still have \eqref{r_11=0} since its proof  does not require $F,G$ to be polynomials on $x$. Thus we have \eqref{rrr=0}. Then as in the proof of
Theorem \ref{C[x,y]'}, we have $p\ne-1$. This completes the proof of the theorem.
\hfill$\Box$
\subsection{Reducing the problem to the case $p\le0$}\label{S3.4}
In this subsection, we want to show that the proof of the two-dimensional Jacobian conjecture  can be reduced to the case $p\!\le\!0$\VS{-7pt}.

\begin{lemm}\label{lemm-g-d2-i-KEY}Suppose $(F,G)$ is a \QJ  pair in $\C[x^{\pm1},y]$ with $p=\frac{p'}{q}>0$ and $m_0\ge0$.
 Then the primary
polynomial $\pri{F}$ satisfies one of the following $($up to
re-scaling variable $x)$\VS{-5pt}:\begin{itemize}\parskip-3pt
\item[{\rm(i)}] $\pri{F}=y+x^{p'}$ and $q=1$.
\item[{\rm(ii)}] $\pri{F}=y^{q}+x$ and $q>1,\,m_0=0$.
\item[{\rm (iii)}] $\pri{F}\!=\!(y^{q}\!+\!x)^{i_1}(y\!+\!\d_{q,1}\a x)^{i_2}$ for some $1\!\ne\!\a\!\in\!\C$,
$i_1,i_2\!\ne\!0,\,{\rm gcd}(i_1,i_2)\!=\!1,\,(i_1,i_2)\!\ne\!(1,1)$, $m_0\!=\!0$.
\end{itemize}
\end{lemm}
\noindent{\it Proof.~}~If $d=1$, then we have (i)
(cf.~(\ref{primary-F})). Thus suppose $d\ge2$. If $a\le-2$, then
$\pri{F}$ divides the left-hand side of the first equation of
(\ref{equa-int1}), a contradiction.
 Thus $a\ge-1$.

Suppose $a=-1$. Then the second case of
(\ref{degree-P}) cannot occur since $d_P\ge0$ and $d\ge2$.  \VS{-5pt}Thus
\begin{equation}\label{d+dp<=2}
d+d_P=\frac{m(p'+q)}{p'm+m_0q}\le
\frac{m(p'+q)}{mp'}=1+\frac{q}{p'}\le 1+q\le 1+d.
\end{equation}
So $d_P\le1$. If $d_P=0$, i.e., $\ptl_yP=0$, then the second
equation of (\ref{equa-int1}) shows $d=1$, a contradiction. If
$d_P=1$, then all equalities must hold in (\ref{d+dp<=2}), i.e.,
$m_0=0,\,p'=1,\,d=q$. Thus we have (ii) by Lemma \ref{Factor-q}(2).

Now suppose $a\ge0$. As in the proof of Theorem \ref{C[x,y]'}, we have
\eqref{equa-fac1===}--\eqref{equa-all-k+}.
If $a>0$, then either case of (\ref{equa-all-k}) shows $i_k<j_k$ for
all $k$ (note from (\ref{asas}) that $-\a_1>\a_2$ since $p'+q>0$),
i.e., $\pri{F}|P$, a contradiction with our choice of $P$.
Thus $a=0$. If $\ell=1$, as in the arguments after \eqref{equa-all-k+}, we
obtain $d=1$ (since $\pri F$ is power free, cf.~\eqref{pri-pro}). Thus $\ell\ge2$ (and so $s\ge\ell\ge2$). Again, as in the arguments after
\eqref{equa-all-k+}, we  \VS{-9pt}have
\begin{eqnarray}\label{Final-snow}
0<d_P=s=\frac{m(p'+q)}{m_0q+p'm}\le 1+\frac{q}{p'}.
\end{eqnarray}
Thus we have the first case of (\ref{degree-P}). Note that any irreducible polynomial $Q$
in the ring $\F[x^{\pm1}][y]$ does not contain an irreducible factor
of multiplicity $\ge2$ in $\AA[y]$. 
 If $d_P<q$ (then $q>d_p=s\ge\ell\ge2$), by Lemma \ref{Factor-q}(1), $P$  has only one
different irreducible factor, i.e., $s=1$, a contradiction.
Thus $d_P\ge q$. Assume
$q=d_P$, which is equal to $s\ge2$. Then $p'=1$ by (\ref{Final-snow}) because: if $p'\ge3$ then $1+\frac{q}{p'}\le1+\frac{q}{3}<q$; if $p'=2$ then
$q>2$ (since $p',q$ are coprime) and $1+\frac{q}{p'}=1+\frac{q}{2}<q$.
Hence \eqref{Final-snow} gives $m=q^2m_0$.
However, since we have the first case of (\ref{degree-P}), we see $b_\mu\ne0$ by \eqref{tildeB}, and $\frac{m_0}{m}\mu$
is integral by Lemma \ref{bi-lemm}. Then
by \eqref{mu-set},
$\frac{m_0}{m}\cdot\frac{m(p'+q)}{m_0q+p'm}$ (which, by
\eqref{Final-snow}, is equal to
$\frac{m_0}{m}d_P\!=\!\frac{m_0}{m}q\!=\!\frac1q$)
is integral, a
contradiction with the fact that $q\!=\!d_P\!\ge\!2$. Thus $d_P\!>\!q$ and  by \VS{-5pt}\eqref{Final-snow},
\begin{equation}\label{j=}d_P=q+1,\ \ \ \ p'=1,\ \ \ m_0=0\VS{-5pt}.
\end{equation}
First suppose $d_P\ge3$. Let $H\in\F[x^{\pm1}][y]$, a monic
polynomial of $y$ of degree $k$, be an irreducible factor of $\pri{F}$ (in the
ring $\F[x^{\pm1}][y]$). Then Lemma \ref{Factor-q} shows that either $\deg_yH=1$
(and $H=y$ by noting that $q\ge2$), or $q|k$ (in this case $H$ has
$k$ different irreducible factors in $\AA[y]$). 
Since
$\pri{F}$ has only $\ell$ different irreducible factors in
$\AA[y]$ 
and $\ell\le d_P=q+1$, we see that $\pri{F}$ has to
have the form
 (up to  re-scaling  $x$\VS{-5pt})
\begin{equation}\label{F1'sform}
\pri{F}=(y^{q}+x)^{i_1}y^{i_2}\mbox{ \ for some }i_1,i_2\in\Z_+\mbox{ \ \ with \ }qi_1+i_2=d
\VS{-5pt}.\end{equation}
 If
$d_P\!=\!s\!=\!2$, then $q\!=\!d_P\!-\!1\!=\!1$, $p\!=\!\frac{p'}{q}\!=\!1$, and $\ell\!=\!2$ since
$2\!\le\!\ell\!\le\! s$. \VS{-1pt}Thus
$\pri{F}\!=\!(y\!+\!x)^{i_1}(y\!+\!\a x)^{i_2}$ (up to  re-scaling  $x$)
for some $1\!\ne\!\a\!\in\!\C,\,i_1,i_2\!\in\!\Z_+$ by noting that each irreducible factor (in $\AA[y]$)
of the $p$-type q.h.e.~$\pri{F}$ is a $p$-type q.h.e., hence, of the form $y+\b x^p=y+\b x$ for some
$\b\in\C$.

In any case,  $i_1,i_2$ are coprime by
(\ref{pri-pro}). If $i_1\!=\!0$ (then $i_2\!=\!1$) or $i_2\!=\!0$ (then
$i_1\!=\!1$), or $i_1\!=\!i_2\!=\!1$, we see that $\pri{F}|P$, a contradiction. Thus
we have (iii). (We can also prove $i_1\ne i_2$ as  follows:
 If $i_1\!=\!i_2$, then  using  $m_0=0$ and Definition \ref{comp}(2)(ii), we have $F\comp0\!=\!\pri{F}^{m'}\!=\!(y^q\!+\!x)^{i_1m'}(y\!+\!\d_{q,1}\a x)^{i_1m'}$, which
is a Jacobian element by Lemma \ref{quasi-J}.
By re-denoting $y^q\!+\!x,\,y\!+\!\d_{q,1}\a x$ to \VS{-3pt}be $x,\,y$ respectively, we obtain from $F\comp0$ that
$H\!:=\!x^{i_1m'}y^{i_1m'}$ is a Jacobian element, but $(-1,-1)\!\in\!\Supp{\sc\,}H^{-\frac1{i_1m'}}$, a contradiction with Theorem \ref{Jacobian-el}.)
This proves
the lemma. \hfill$\Box$\vskip7pt

Now for any \QJ  pair $(F,G)$ in $\C[x^{\pm1},y]$, if we have Lemma \ref{lemm-g-d2-i-KEY}(i), then by
applying the automorphism $(x,y)\!\mapsto\!(x,y\!-\!x^{p'})$, $F$ becomes an element with prime degree $<\!\frac{p'}{q}$. If we have Lemma \ref{lemm-g-d2-i-KEY}(ii) or (iii), then by applying the automorphism $(x,y)\!\mapsto\!(x\!-\!y^q,y)$, $F$ becomes  an element with a lower $y$-degree.
Thus by induction on $\deg_yF$ \VS{-2pt}and on the prime degree $p$, after a finite steps
(note that we have only finite possible choices of $p\!=\!\frac{p'}{q}$ since $1\!\le\! q\!\le\! \deg_yF$ by \eqref{q=}), we can suppose
either $F$ has $y$-degree $\!\le\!1$, or else,
$p\!\le\!0$ (this implies $m_0\!>\!0$, otherwise, by the definition of the prime degree $p$
in Definition \ref{defi2.1}, we would obtain $F,G\!\in\!\C[x^{-1},y]$ and $(F,G)$ cannot be a \QJ  pair). Thus from now on, we can suppose $m_0\!>\!0$ and $p\!\le\!0$.
Note that if the original $F,G$ are in $\C[x,y]$, then after the above process, the resulting $F,G$ are still in $\C[x,y]$. In this case, by Remark \ref{m0-m},
we can suppose $0\!<\!m_0\!<\!m$ and $p\!\le\!0$. The 
above have in fact \VS{-3pt}proved
\begin{lemm}\label{voit}
If there exists a Jacobi pair in $\C[x,y]$ which violates the two-dimensional Jacobian conjecture, then there exists a Jacobi pair $(F,G)$ in $\C[x,y]$ satisfying $0<m_0<m$ and $p\le0$.
\end{lemm}

{}For our purpose of discussions, from now on we  \VS{-5pt}assume\equa{ASSMM}{\mbox{ $F,G\in\C[x^{\pm1},y]$ with $0<m_0<m$, and $p\le0$.}}
\VS{-17pt}
In case $p\!<\!0$, we always write
$p\!=\!-\frac{p''}{q}$ for some coprime positive integers $p'',q$.
We
can suppose $d\!\ge\!2$ because (cf.~\eqref{primary-F}): If $p\!=\!0$, by replacing $y$ by $y\!+\!\l$, where $y\!=\!\l$ is root of $\pri{F}$, we can suppose $\pri c_d\!=\!0$ (in case $\pri F$ becomes a monomial after the replacement, this means that $p$ becomes $<0$); if $p\!<\!0$, then since $p\!>\!-\frac{m_0}{m}\!>\!-1$ by Lemma \ref{p-ge---} and \eqref{ASSMM}, we have $p\!\notin\!\Z$, and $\pri c_1$ has to be
zero.

We can use similar arguments as those in
the proof of Lemma \ref{lemm-g-d2-i-KEY} to obtain
$a=-1,0$ in \eqref{equa-int1} (even though in our case here $p\le0$, some arguments in the proof of Lemma \ref{lemm-g-d2-i-KEY} does not depend on $p$; for instance, we again have
\eqref{ccc} and \eqref{equa-all-k}, so if $a\ge1$, then either of \eqref{equa-all-k} again shows $i_k<j_k$
because by Lemma \ref{p-ge---}, $p>-\frac{m_0}{m}>-1$, i.e., $p'+q>0$ still holds). If
we drop the condition $\pri{F}\mbox{$\not|$}\,P$, we can always suppose
$a=0$ because if $a=-1$ we can re-denote $P$ to be $\pri{F}P$ so that
$a$ become zero. Furthermore,
if we consider $p\le0$ and use similar arguments as in the proof of Lemma \ref{lemm-g-d2-i-KEY}
(cf.~(\ref{equa-fac2})--\eqref{equa-all-k}), we can  \VS{-4pt}obtain
 \begin{lemm}\label{irr-factor}
Every irreducible factor of $\pri{F}$ in $\C[x^{\pm1},y]$ is a factor of
$P$, and every irreducible factor of $P$ has multiplicity $1$.
So
$2\!\le\! d\!\le\!d_P\!:=\!\deg_yP\!=\!\frac{1+p}{\frac{m_0}{m}+p}\!\in\!\Z$ $($thus $b_\mu\!\ne\!0$,
cf.~$(\ref{mu-set}))$\VS{-2pt}. Furthermore, if $\pri{F}\mbox{$\not|$}\,P$, then every irreducible factor of $P$ in $\C[x^{\pm1},y]$ is also a factor of $\pri{F}\VS{-5pt}.$\end{lemm}

\begin{lemm}\label{y|P}We can suppose $y|P$\VS{-7pt}.\end{lemm}

\ni{\it
Proof.~}~If $p=0$, then since we can suppose $\pri c_d=0$, i.e., $y|P$.
So assume
$p=-\frac{p''}{q}<0$\VS{-2pt}. Since $b_\mu\ne0$, we see $\mu,\frac{m_0}{m}\mu$
are integers (cf.~Lemma \ref{bi-lemm})\VS{-3pt}. By \eqref{mu-set},
$\frac{m(q-p'')}{m_0q-mp''},\,\frac{m_0(q-p'')}{m_0q-mp''}$ are
integers. Subtracting them by $1$, we see
$\frac{(m-m_0)q}{m_0q-mp''},\,\frac{(m-m_0)p''}{m_0q-mp''}$ are
integers. Since $p'',q$ are coprime, we see
$r=\frac{m-m_0}{m_0q-mp''}$ is an integer. Since $d_P=qr+1$ (cf.~\eqref{degree-P}), by
Lemma \ref{Factor-q}(1), we see $P$ must have a monic irreducible factor
of degree $1$, which has to be $y$.
 \hfill\VS{-7pt}$\Box$

\subsection{Newton polygons}\label{S-NewtonP}

In this subsection, we start from a \QJ  pair  $(F,G)$ satisfying \eqref{ASSMM},
\VS{-2pt}and  $2\!\le\! m\!<\!n$, $m\mbox{$\not|$}\,n$, and regard $F,G$ as in $\C[x^{\pm\frac1N},y]$ for some sufficient large $N\!\in\!\N$ (at the beginning we can take $N\!=\!1$, later on we may need to enlarge $N$ after we apply some variable changes).
We define the {\it Newton polygon} of $F$ to be the minimal polygon such that the region surrounded by it contains
$(\Supp\,F)\cup\{0\}$. Thus the Newton polygon of $F$ may look as in \eqref{new-supp-F-region}, where, $m_0,\dot m_0,\Dhat m_0\!\in\!\Q_+,\,m,\dot m,\Dhat m\!\in\!\N$, and points such as $\tau,\dot \tau,\Dhat \tau,...$ are called {\it vertices} of $\Supp\,F$, \VS{-1pt}and line segments such as $L,\,\dot L,\,\Dhat L$,... are called {\it edges} of $\Supp\,F$; we shall always use $(n_0,n),$ $(\dot n_0,\dot n),$ $(\Dhat n_0,\Dhat n),...$
to denote corresponding vertices of $G$\VS{-8pt}.
{\small\equa{new-supp-F-region}{\mbox{$\ \ \ \ \ \ $}%
 \put(-5,3){$\vector(1,0){100}$}\put(-2,-5){$\sc0$}\put(95,-5){$\sc
 x$}
 \put(5,-3){$\vector(0,1){80}$}\put(-2,74){$\sc y$}
\put(5,70){\put(0,-67){\line(-1,0){100}}
\put(-22,-1){$\sc\circ$}
\put(-20,0){\line(-1,-1){44}}
\put(-67,-47){$\sc\circ$}\put(-64,-47){\line(3,-2){30}}\put(-34,-69){$\sc\circ$}
$\put(-20,0){$\line(1,0){79}$}$}\put(65,70){$\line(0,-1){15}$}
\put(64,68){$\sc\circ$}\put(65,74){$\sc\tau=(m_0,m)\in\Supp\,F$}
\put(64,52){$\sc\circ$}\put(69,48){$\sc\dot\tau=(\dot m_0,\dot m)
\in\Supp\,F$}
\put(65,55){$\line(-1,-2){16}$}\put(47,20){$\sc\circ\,\Dhat \tau=(\Dhat m_0,\Dhat m)\in\Supp\,F$}
\put(48,21){\line(-1,-1){17}}\put(28,1){$\sc\circ$}
%
%
%
%
\put(15,40){$\sc\Supp\,F$} \put(43,8){$\sc \Dhat L:\mbox{ \footnotesize slope } \frac{\Dhat q}{\Dhat p}$}
\put(60,35){$\sc \dot L:\mbox{ \footnotesize slope } \frac{\dot q}{\dot p}$}\put(69,59){$\sc L$}
 \ \hspace*{80pt}\
}}
We always assume $p,q$ together with their different versions (e.g., the ``\ $\dot{ }$\ '' version,
the ``\ $\Dhat{ }$\ '' version, etc.) are in $\Z_+$, and $q$ together with  its different versions are nonzero.

The next lemma says that starting from any (not necessarily the top most)
vertex $\dot\tau=(\dot m_0,\dot m)$, we can find a lower vertex $(\Dhat m_0,\Dhat m)$ as long as the edge $L$ satisfies some condition.
\begin{theo}\label{Nerw}
Suppose $\dot\tau\!=\!(\dot m_0,\dot m)\!\in\!\Supp{\sc\,}F$ is a vertex of $\Supp{\sc\,}F$ such that $0\!<\!\dot m_0\!<\!\dot m$ with $\dot m_0\!\in\!\Q_+,\,\dot m\!\in\!\N$, and
$\dot L$ is the edge of $\Supp{\sc\,}F$ as shown in \eqref{new-supp-F-region} such that
its slope is either $\frac{\dot q }{\dot p }\!>\!0$ $($cf.~the statement before \eqref{p-line}$)$ or $\infty$ $($in this case, we write its slope as $\frac{\dot q }{\dot p }$ with $\dot p \!=\!0,\,\dot q \!=\!1)$\VS{-2pt}.
Suppose correspondingly $(\dot n_0,\dot n)$ is a vertex \VS{-5pt}of $\Supp{\sc\,}G$ such that $\frac{\dot n_0}{\dot m_0}\!=\!\frac{\dot n}{\dot m}\!=\!\frac{n}{m}$\VS{-2pt}, and the edge of $\Supp{\sc\,}G$ whose top vertex is $(\dot n_0,\dot n)$ also has \VS{-5pt}slope $\frac{\dot q }{\dot p }$.
Assume $\dot m\!\ge\!2$ and $($cf.~Lemma $\ref{p-ge---}\VS{-9pt})$
\equa{p0-q0}{\dis-\frac{\dot p }{\dot q }>
-\frac{\dot m_0}{\dot m}+\frac{\dot m-\dot m_0}{\dot m(\dot m+\dot n-1)}=-\frac{\dot m_0+\dot n_0-1}{\dot m+\dot n-1}.}
Then either $\dot F\Comp0$ $($the
part of $F$ with support being the edge $\dot L)$ has \VS{-4pt}only one irreducible factor, or else,
if necessary by
changing $(x,y)$ to  $(x,y+\a_0 x^{-\frac{\dot p }{\dot q }})$ for some $\a_0\in\C$, there is a \VS{-2pt}vertex $\Dhat \tau\!=\!(\Dhat m_0,\Dhat m)\!\in\!\Supp\,F$ and a corresponding vertex $(\Dhat n_0,\Dhat n)$ \VS{-5pt}of $\Supp\,G$ such that $\frac{\Dhat n_0}{\Dhat m_0}\!=\!\frac{\Dhat n}{\Dhat m}\!=\!\frac{n}{m}$ an\VS{-4pt}d \begin{equation}\label{Vertex-2}0<\Dhat m_0<\Dhat m<\dot m\mbox{ \ with \ $\Dhat m_0\in\Q_+,$ \ \ and \ \ $2\le \Dhat m\in\N.$}\end{equation}
\end{theo}\VS{-10pt}\PAR\NOindent
\ni{\it Proof.~}~Suppose \VS{-2pt}$F=\sum_{i,j}f_{ij}x^iy^j\in\C[x^{\pm\frac1N},y]$ for some $f_{ij}\in\C$, and  $\dot \tau,\,\dot L$ are as  shown in
\eqref{new-supp-F-region}. Then $\frac{\dot m_0-\Dhat m_0}{\dot p }=\frac{\dot m-\Dhat m}{\dot q }$,
and 
$\dot F\Comp{0}=\sum_{(i,j)\in \dot L}f_{ij}x^iy^j,$ which can be re-written \VS{-9pt}as
\equa{F====0000}{
\mbox{$\dot F\Comp0=\sum\limits_{j=0}^{\dot m}f'_{j}x^{\dot m_0-j\frac{\dot p }{\dot q }}y^{\dot m-j}$, \ where }f'_{j}=f_{\dot m_0-j\frac{\dot p }{\dot q },\dot m-j}.
}
Thus $\dot F\Comp{0}$ is
a $-\frac{\dot p }{\dot q }{\sc\,}$-type q.h.p (cf.~Definition \ref{comp}(3)). Now regard
$-\frac{\dot p }{\dot q }$ as the ``prime degree'' of $F$ and $\dot F\Comp0$ as the leading
polynomial of $F$, and
 define $-\frac{\dot p }{\dot q }{\sc\,}$-type $r$-th
component $\dot F\Comp r$ of $F$ as in (\ref{eq-comp}) (with the data $(p,m_0,m)$ being $(-\frac{\dot p }{\dot q },\dot m_0,\dot m)$).
We can write $F$  and define
$F^a\in\CC$ for $a\in\Q$ as
(where $\dot F\Comp{<0}=\sum_{r<0}\dot F\Comp r$ is the
ignored polynomial, cf.~Definition \ref{comp}(2)(i)\VS{-9pt})
\equa{f-a-------}{F=\dot F\Comp0+\dot F\Comp{<0},\ \ \ \ F^a=\sum\limits_{s=0}^\infty\dbinom{a}{s}\dot F\Comp0^{^{\sc
a-s}}\dot F\Comp{<0}^{^{\sc s}},}
i.e., we expand $F^a$ according to its
$-\frac{\dot p }{\dot q }{\sc\,}$-type components. This is well-defined.
{\def\tau{z}%
In fact, if we \VS{-4pt}let
$\tau\in\C\bs\{0\}$ be an indeterminate, and replace $x,y$ by $\tau
x,\tau^{-\frac{\dot p }{\dot q }} y$, and regard elements as in the field
(cf.~(\ref{CC-TAU})\VS{-12pt}) \equan{SSSSSSSSSS}{
\CC_\tau=\{F=\sum\limits_{j\in\Q}F_j\tau^j\,|\,F_j\in\CC,\,\Supp_\tau
F\subset a-\dis\frac{1}{b}\Z_+\mbox{ \ for some }a,b\in\Z,\,b>0\}, }
then the $-\frac{\dot p }{\dot q}{\sc\,}$-type $r$-th component $\dot F\Comp{r}$ is simply a
$\tau$-homogenous element of $\CC_\tau$, and (\ref{f-a-------})
simply means that we expand $F^a$ as an element in $\CC_\tau$. Then
we can use all arguments before and  define $\dot R_0$ as in
(\ref{l-r0}) such that $(\dot F\Comp{0},\dot R_0)$ is a \QJ
pair (cf.~Lemma \ref{quasi-J}), and in case $\dot F\Comp0$ has at least two irreducible factors, we have (cf.~\eqref{R-0not=0}, Lemmas
\ref{R===0} and \ref{irr-factor}\VS{-7pt})\equa{NOWOO-R}{\mbox{$\dot R_0=\dot F\Comp0^{-1}\dot P$ for
some $\dot P\in\C[x^{\pm\frac1N},y]$ of $y$-degree
$\dis d_{\dot P}=\frac{(\dot q -\dot p )\dot m}{\dot m_0\dot q -\dot m \dot p }\ge2$.}}
}\VS{-19pt}\PAR\NOindent%
To be more precise, we give detailed arguments below. Write $G=\dot G\Comp0+\dot G\Comp{<0}$ as in
\eqref{f-a-------}, and suppose $\dot G\Comp0$ has the highest term $x^{\dot n_0}y^{\dot n}$ for some $\dot n_0\in\Q_+,$ $\dot n\in\N$ such that $\frac{\dot m_0}{\dot n_0}=\frac{\dot m}{\dot n}$.
\VS{-7pt}Denote\begin{eqnarray}\label{tilde-F-G}\!\!\!\!\!\!&\!\!\!\!\!\!\!\!\!\!\!\!\!\!\!\!\!\!\!\!&
\widetilde F=z^{\frac{\dot p }{\dot q }\dot m-\dot m_0}F(zx,z^{-\frac{\dot p }{\dot q }}y),\
\ \widetilde G=z^{\frac{\dot p }{\dot q }\dot n-\dot n_0}G(zx,z^{-\frac{\dot p }{\dot q }}y),\ \
 \widetilde J:=[\widetilde F,\widetilde G]=z^{-\dot \mu}J,
\end{eqnarray}where $\dot \mu=\dot m_0+\dot n_0-1-
\frac{\dot p }{\dot q }(\dot m+\dot n-1)>0$ (cf.~\eqref{p0-q0}), $J=[F,G].$
Then clearly, for any $r\in\Q_+$,
the $-\frac{\dot p }{\dot q }\,$-type $r$-th components $\dot F\Comp{r}$ of $F$ and $\dot G\Comp{r}$ of $G$ are simply the coefficients \VS{-4pt}of $z^r$ in $\widetilde F$ and $\widetilde G$ respectively.
As in
\eqref{equa2.5}, we can write $\dot G\Comp0=\sum_{i=0}^\infty \tilde b_{i0}\dot F\Comp0^{\frac{\dot n-i}{\dot m}}$ \VS{-3pt}for some $\tilde b_{i0}\in\AA$ with $\tilde b_{00}=1$.
Denote $\widetilde G{}^1=\widetilde G- \sum_{i=0}^\infty \tilde b_{i0}\widetilde F^{\frac{\dot n-i}{\dot m}}$. Then $\deg_z\widetilde G{}^1\in\Q$ is $<0$, denoted by $-r_1$.
Write $\widetilde G{}^1=z^{-r_1}(\dot G{}^1\Comp0+\dot  G{}^1\Comp{<0})$ as in
\eqref{f-a-------}, and suppose $\deg_y\dot G{}^1\Comp0=k_{r_1}$. Then again we can \VS{-7pt}write \equa{tilde-G-111}{\mbox{$\dot G{}^1\Comp0=\sum\limits_{i=0}^\infty \tilde b_{i,r_1}\dot F\Comp0^{\frac{k_{r_1}-i}{\dot m}}$ \ for some $\tilde b_{i,r_1}\in\AA$ \ with $\tilde b_{0,r_1}\ne0$.}}\VS{-19pt}\PAR\NOindent Continuing this way, we can \VS{-7pt}write (note that if $\dot \mu=0$, the following still holds with the first summand vanishing\VS{-5pt})
\equa{widetilde-G==}{\widetilde G=\mbox{$\sum\limits_{r\in\Q_+,\, r<\dot \mu}
z^{-r}\sum\limits_{i=0}^\infty$}\tilde b_{ir}\widetilde F^{\frac{k_r-i}{\dot m}}+z^{-\dot \mu}\widetilde R,
}\VS{-19pt}\PAR\NOindent for some $k_r\in\Z$, $\tilde b_{ir}\in\AA$ with $k_0=\dot n$, and some $\widetilde R\in\CC_z$ with $\deg_z\widetilde R\le0$.
As in Lemma \ref{lemm-product-KH}(6), there exists some $\a\in\N$ such that $\tilde b_{ir}=0$ if $r\notin\frac1{\a}\Z_+$. Assume there exists $r_0<\dot \mu$ being smallest such that $\ptl_x\tilde b_{i_0,r_0}\ne0$ for some $i_0\in\Z_+$ (and we take $i_0$ to be smallest). Comparing the coefficients of $z^{-r_0}$ in the last equation of \eqref{tilde-F-G}, we easily obtain a contradiction (we have an equation similar to the first equation of \eqref{compar-coeff}). Thus $\tilde b_{ir}\in\C$.
Since $\dot G{}^1\Comp0$ in \eqref{tilde-G-111} is a $-\frac{\dot p }{\dot q }\,$-\VS{-5pt}type q.h.e.~(cf.~Definition \ref{comp}(3)) and $\Supp\,\dot F\Comp0^{\frac{k_{r_1}-i}{\dot m}}$ lies in a
different line with slope $\frac{\dot q }{\dot p }$ for different $i$, there is at most one $i$ such that  $\tilde b_{i,r_1}\ne0$. Thus
we can assume $\tilde b_{i,r_1}=0$ for $i>0$, and 
\eqref{widetilde-G==} can be rewritten as
(although we do not need the actual value of $k_r$, one can compute that
$k_r=\dot n-\frac{\dot q \dot mr}{\dot m_0\dot q -\dot m\dot p }$ \VS{-5pt}if $\tilde b_r\ne0$ by using the fact that $\Supp\,\dot F\Comp0^{\frac{k_r}{\dot m}}$ lies in the line passing through $(\dot n_0-r,\dot n)$ with slope $\frac{\dot q }{\dot p }$, i.e., $(\frac{\dot m_0}{\dot m}k_r,k_r)=(\dot n_0-r-i_0\frac{\dot p }{\dot q },\dot n-i_0)$ for some $i_0\in\Z$\VS{-7pt})
\equa{Qwidetilde-G==}{\ \ \ \ \ \ \ \widetilde G=\mbox{$\sum\limits_{r\in\frac1{\a}\Z_+,\,r<\dot \mu}
z^{-r}$}\tilde b_{r}\widetilde F^{\frac{k_r}{\dot m}}+z^{-\dot \mu}\widetilde R\mbox{ \ for some \ }\tilde b_{r}\in\C.
}
Hence, $[\widetilde F,\widetilde R]\!=\!z^{\dot \mu}[\widetilde F,\widetilde G]\!=\!z^{\dot \mu}\widetilde J\!=\!J$.
Write $\widetilde R\!=\!\dot R\Comp0\!+\!\dot R\Comp{<0}$ as before. Then we obtain $[\dot F\Comp0,\dot R\Comp0]\!=\!J$.
As before, \VS{-3pt}we use $\dot{\pri{F}}$ to denote the primary polynomial of $\dot F\Comp0$. Then as in the proof of \VS{-4pt}Lemma \ref{d|b}, we see from \eqref{Qwidetilde-G==} that  $\dot F\Comp0^{\frac{k_r}{\dot m}}$ is an integral power of $\dot{\pri{F}}$ if $\tilde b_r\ne0$. Hence
 $\dot R\Comp0$ is a rational function on $y$ of the form $\dot F\Comp0^{-b}\dot P$ for some $b\in\Z$, $\dot P\in\C[x^{\pm\frac1N},y]$ (as in Lemma \ref{R===0}).
Noting that $(\deg_x\dot R\Comp0,\deg_y\dot R\Comp0)\!=\!(\dot n_0\!-\!\dot \mu,\dot n)\!+\!\b(\frac{\dot p }{\dot q },1)$ for some $\b\!\in\!\Q$ (by considering $\Supp\,\dot R\Comp0$), which must be also equal to
either $(1\!-\!\dot m_0,1\!-\!\dot m)$ or $\g(\dot m_0,\dot m)$  for some $\g\!\in\!\Q$ (by the fact that $[\dot F\Comp0,\dot R\Comp0]\!=\!J$), we can solve that $\deg_y\dot R\Comp0=1-\dot m$ or
 $\frac{(\dot q -\dot p )\dot m}{\dot m_0\dot q -\dot m\dot p }-\dot m$ (as in \eqref{R-0not=0}).
Thus, if $\dot{\pri{F}}$ has at least two irreducible factors (in $\AA[y]$),
we can obtain \eqref{NOWOO-R} by taking $\dot R_0=\dot R\Comp0$
 (cf.~Lemmas  \ref{quasi-J} and \ref{irr-factor}).
%
%

We may suppose $\dot{\pri{F}}$ has at least 2 irreducible
factors  (in $\AA[y]$). \VS{-10pt}Otherwise \equa{OTHWWW}{\mbox{$\dot F\Comp{0}=x^{\dot m_0}(y- \a_0 x^{-\frac{\dot p }{\dot q }})^{\dot m}$ for some $0\ne\a_0\in\C$.}}
By \VS{-3pt}replacing $y$ by $y+ \a_0 x^{-\frac{\dot p }{\dot q }}$ (in this case since $x^{-\frac{\dot p }{\dot q }}\in\C[x^{\pm\frac1N},y]$ we must have $\dot q|N$ and we do not need to enlarge $N$), \VS{-3pt}$F$ becomes an element with
the ``prime degree'', denoted $-\frac{p'}{q'}$, being smaller than \VS{-3pt}$-\frac{\dot p }{\dot q }$.
Since $-\frac{p'}{q'}\ge -\frac{\dot m_0}{\dot m}$ (as in the proof of Lemma \ref{p-ge---}) and $1\le q'< N\dot m$
(cf.~\eqref{q=}, noting that $F$ is now in $\C[x^{\pm\frac1N},y]$), we only have finite possible choices of $-\frac{p'}{q'}$. Thus $F$ can eventually become an element such that either the ``\,$>$\,'' in \eqref{p0-q0} becomes equality  (i.e., $\dot \mu$ becomes zero),
or else \eqref{p0-q0} still holds but
$\dot{\pri{F}}$ has at least 2 irreducible
factors. Let us assume the later case happens
({\bf remark}: if $F,G\in\C[x^{\pm1},y]$ and \eqref{OTHWWW} occurs, then $\frac{\dot p }{\dot q }\in\Z$ and so $\dot p=0$ since the right-hand side of \eqref{p0-q0} is $>-1$, and thus after the above replacement, $\dot p$ becomes nonzero, and we still have the following three facts:
(i) $F,G\in\C[x^{\pm1},y]$\VS{-3pt};
(ii) $\dot F\Comp0\in x{\ssc\,}\C[x,y]$ since in \eqref{F====0000},
 $\dot m_0-j\frac{\dot p }{\dot q }\ge \dot m_0-\dot m\frac{\dot p }{\dot q }>0$ \VS{-3pt}(whether or not the ``\,$>$\,'' in \eqref{p0-q0} becomes equality, we always have $-\frac{\dot p }{\dot q }>-\frac{\dot m_0}{\dot m}$), analogously, $G\Comp0\in x{\ssc\,}\C[x,y]$\VS{-2pt}; (iii) $\dot \mu\!>\!0$ (which is equivalent to \eqref{p0-q0}), otherwise
\eqref{widetilde-G==} with $\dot \mu=0$ shows $\widetilde G\!=\!\widetilde R$ and so $(F\Comp0,G\Comp0)\!=\!(\dot F\Comp0,\dot R\Comp0)\!=\!J$, a contradiction with fact (ii), thus fact (iii) shows that we always have the later case, and in particular after the variable change \eqref{x-yyyyy}, we {\bf always} have a side $\dot L$ with negative slope in \eqref{new-supp-FCOM-region}).

If $\dot p =0$, then Lemma \ref{irr-factor} shows $\dot m_0|\dot m$ (in general, for $a,b\in\Q$ with $a\ne0$,
notation $a|b$ means $\frac ba\in\Z$),
$d_{\dot P}=\frac{\dot m}{\dot m_0}$\VS{-3pt}, and the leading polynomial
$\dot F\Comp{0}$ has at most $d_{\dot P}$ irreducible
factors, thus at least an irreducible factor of $\dot F\Comp{0}$ has
multiplicity $\ge\frac{\dot m}{d_{\dot P}}=\dot m_0$\VS{-3pt}. If every irreducible factor of
$\dot F\Comp{0}$ has multiplicity $\dot m_0$, then $\dot{\pri{F}}$ must be a power
of $\dot P$, which contradicts the dot version of  (\ref{equa-int1}). Thus at least one
irreducible factor, say, \equa{f0-p==0}{\dot{\pri{f}}_1=y-\a_0\mbox{ for some $\a_0\in\C$, has maximal multiplicity, say, $\Dhat m>\dot m_0$.}}
 Now replacing $y$ by $y+\a_0$, we obtain that $(\Dhat m_0,\Dhat m)$ (with $\Dhat m_0=\dot m_0$) is a vertex of $F$ satisfying
$0\!<\!\Dhat m_0\!<\!\Dhat m\!<\!\dot m$ (we {\bf remark} here that if our starting pair $(F,G)$ is in $\C[x,y]$, then the resulting pair after the variable change
is still in $\C[x,y]$). 
Now assume $\dot p \!>\!0$. As above, suppose the following irreducible factor   (in $\AA[y]$) of $\dot{\pri{F}}$\VS{-13pt},
\equa{SF1===}{\ \ \ \ \ \ \ \ \ \ \mbox{$\dot{\pri{f}}_1\!=\!y\!-\!\a_0 x^{-\frac{\dot p }{\dot q }}$
has maximal
multiplicity, denoted $\Dhat m$ (thus $\Dhat m<\dot m$).}}
Then \VS{-1pt}$\Dhat m>\frac{\dot m}{d_{\dot P}}=\frac{\dot q \dot m_0-\dot p {\ssc\,}\dot m}{\dot q -\dot p }$ (cf.~\eqref{degree-P})\VS{-3pt}.
Applying the variable change
 (if necessary we 
enlarge $N$ by  $\dot q$ times  to ensure that after the variable change, all elements under consideration are in $\C[x^{\pm\frac1N},y]$)
\equa{Var1111}{\si:(x,y)\mapsto(x,y+\a_0 x^{-\frac{\dot p }{\dot q }}),}
we obtain a
vertex $(\Dhat m_0,\Dhat m)$, where $0\!<\!\Dhat m_0\!=\!\dot m_0\!-\!(\dot m\!-\!\Dhat m_0)\frac{\dot p }{\dot q }\!<\!\Dhat m\!<\!\dot m$.

In any case, we have \eqref{Vertex-2} except that we have not yet proved $2\le \Dhat m$.
\VS{-2pt}Using \eqref{G's-component}, \eqref{p0-q0} and \eqref{Qwidetilde-G==}, as in the proof of Lemma \ref{p-ge---}, we have
$\dot G\Comp0=\dot F\Comp0^{\frac{n}{m}}$, and so $\frac{\Dhat n_0}{\Dhat m_0}\!=\!\frac{\Dhat n}{\Dhat m}\!=\!
\frac{\dot n}{\dot m}\!=\!\frac{n}{m}$\VS{-3pt}.
Finally, if $\Dhat m=1$, we would obtain
$n=m\frac{\Dhat n}{\Dhat m}=m \Dhat n$, a contradiction with the assumption that $m\mbox{$\not|$}\,n$. \hfill$\Box$\vskip5pt

The proof of Theorem \ref{Nerw} shows that
 we can eventually obtain a vertex, denoted $(\Dhat m_0,\Dhat m)$\VS{-2pt},
 such that (where $\frac{\Dhat q}{\Dhat p}$ is the slope of the edge $\Dhat L$ with top vertex $(\Dhat m_0,\Dhat m)$\VS{-11pt})
\equa{AMAM}{\dis \!\!\!\!\!\!\!\!\!\!\Dhat m\ge2,\ \ \ \ \frac{\Dhat p}{\Dhat q}=\frac{\Dhat m_0}{\Dhat m}-\frac{\Dhat m-\Dhat m_0}{\Dhat m(\Dhat m+\Dhat n-1)},}
and if $\Dhat F$ is the
part of $F$ with support being the edge $\Dhat L$ and $\Dhat G$ is
the part of $G$ analog to $\Dhat F$, then $(\Dhat F,\Dhat G)$ is a \QJ  pair (cf.~\eqref{Qwidetilde-G==} with $\dot\mu$ replaced by $\Dhat\mu=0$ and statements after \eqref{Qwidetilde-G==}).

Now we apply the variable \VS{-9pt}change (as before, if necessary we 
enlarge $N$ by $\Dhat q\!-\!\Dhat p$ times) \equa{x-yyyyy}{(x,y)\mapsto(x^{\frac{\Dhat q}{\Dhat q-\Dhat p}},x^{-\frac{\Dhat p}{\Dhat q-\Dhat p}}y).}
Note that any line with slope $\frac{\Dhat q}{\Dhat p}$ \VS{-3pt}is mapped under \eqref{x-yyyyy} to a line
parallel to the $y$-axis, so the above pair $(\Dhat F,\Dhat G)$ determined by the edge $\Dhat L$ in \eqref{new-supp-F-region} are mapped to a  pair $(x^af,x^bg)$ for some $a,b\in\Q$ and $f,g\in\C[y]$ with $\deg_yf=\Dhat m,\,\deg_yg=\Dhat n$. Using $[x^af,x^bg]\in\C\bs\{0\}$, \VS{-3pt}we obtain
$a+b=1,\,\frac{a}{\Dhat m}=\frac{b}{\Dhat n}$, i.e., $a=\frac{\Dhat m}{\Dhat m+\Dhat n},\,b=\frac{\Dhat n}{\Dhat m+\Dhat n}$\VS{-3pt}.
Note also that any line with slope $>\frac{\Dhat q}{\Dhat  p}$ is mapped to a line with a negative slope, thus
the edge $\dot L$ with slope $\frac{\dot q}{\dot p}>\frac{\Dhat q}{\Dhat  p}$ \VS{-3pt}in \eqref{new-supp-F-region} is mapped to an edge $\dot L$ with a negative slope, denoted $-\frac1{\dot\a}$ (cf.~\eqref{new-supp-FCOM-region}). Thus,
by exchanging symbols $\Dhat L,\Dhat \tau,\Dhat m,\Dhat n$ and $L,\tau,m,n$, we see that
$F$ and $G$ become (up to nonzero scalars) elements of the forms in \eqref{F-now-is} and
\eqref{G-now-is} (so from now on, $m$ and $n$ {\textbf{do not}} denote the $y$-degrees of $F$ and $G$), such that the Newton polygon of $F$ looks as in \eqref{new-supp-FCOM-region},
 where the existence of the edge $\dot L$ with negative slope $-\frac{1}{\dot\a}$ follows from the remark in a few lines after \eqref{OTHWWW}\VS{-7pt}.
\begin{eqnarray}\label{F-now-is}F\!\!\!&=&\!\!\!
\mbox{$x^{\frac{m }{m +n }}\Big(f+\sum\limits_{i=1}^{M_1}x^{-\frac{i}{N}}f _i\Big)$,  where,
$f=f (y)=y^{m }+\sum\limits_{i=1}^{m }c_iy^{m -i}$ with $c_1=0$,}\\[-3pt]
\label{G-now-is}G\!\!\!&=&\!\!\!\mbox{$x^{\frac{n }{m +n }}\Big(g+\sum\limits_{i=1}^{M_2}x^{-\frac{i}{N}}g _i\Big),$ where, $g=g(y)=y^{n }+\sum\limits_{i=1}^{n }d_iy^{n -i},$}\end{eqnarray}
\VS{-19pt}\PAR\NOindent
for some $m,n\in\N$ (with $2\le m<n$), 
$M_1,M_2\in\Z_+,\,c_i,d_i\in\C$ (we can always assume $c_1=0$ by replacing
$y$ by $y-c$ for some $c\in\C$ if necessary) and $f _i,g_i\in\C[y]$
(note that $x^{\frac{m }{m +n }}f $ is the part of $F$ whose support is the edge $L$, and that
the point $(\frac{m}{m+n},0)$ may not belong to the edge $L$, i.e., $c_m$ can be zero, cf.~the
last statement of this section)\VS{-7pt}.
%
{\small
\equa{new-supp-FCOM-region}{%
 \put(25,3){$\vector(1,0){100}$}\put(78,-5){$\sc0$}\put(125,-5){$\sc
 x$}
 \put(84,-6){$\vector(0,1){80}$}\put(88,68){$\sc y$}
\put(5,70){\put(30,-67){\line(-1,0){110}}
\put(34,0){$\sc\circ$}
%
\put(35,0){\line(-1,-1){67}}\put(28,-1){$\sc\Dhat\tau$}
\put(-100,-49){}\put(-76,-47){}
\put(-34,-69){$\sc \circ$}
$\put(35,1){$\line(3,-1){34}$}$}\put(75,60){$\line(1,-1){25}$}
\put(74,58){$\sc\circ$}\put(71,65){$\sc\,\dot \tau$}\put(60,67){$\sc\Dhat L$}
\put(99,32){$\sc\circ$}\put(104,30){$\sc\,\tau =(\frac{m }{m +n },m )$}
\put(102,13){$\sc\,L:\mbox{ \footnotesize slope } \infty $}
\put(100,35){$\line(0,-1){32}$}
%
%
%
%
\put(35,30){$\sc\Supp\,F$} \put(91,47){$\sc \dot L:\mbox{ \footnotesize slope } -\frac{1}{\dot \a}$}
 \ \hspace*{80pt}\
}}\VS{-19pt}\PAR\NOindent
Furthermore, the proof of Theorem \ref{Nerw} shows that the part $\dot f^0$ of $F$ whose support is the edge $\dot L$ with slope $-\frac1{\dot\a}$ has the \VS{-9pt}form
\equa{f0===}{\dot f^0=x^{\frac{m}{m+n}}y^m\prod\limits_{i=1}^{\de}(a_ix^{-\dot\a}y+1)^{m_i},}
for some $\de,m_i\in\N,\,a_i\in\C\bs\{0\}$
and $m_i\le m$ with inequality holds for at least some $i$
(since $m$ is the maximal multiplicity among all irreducible factors of $\dot f_0$, cf.~\eqref{SF1===}), furthermore, $m,m_i,i=1,...,\de$, have at least a common divisor (otherwise, $n$ must be a multiple of $m$, and so, for the original $F$ and $G$, we also have $m|n$, a contradiction with the assumption).

{}From \eqref{f0-p==0} (and the remark after it), \eqref{Var1111}, \eqref{x-yyyyy}, and
proof of  Theorem \ref{Nerw}, and
Corollary \ref{CCCC}(2), we obtain
\begin{coro}\label{Coeoe}\begin{itemize}\parskip-3pt\item[\rm(1)]
The pair $(F,G)$ is in fact obtained from a Jacobi pair, denoted by $(F,G)$, in $\C[x,y]$ by applying
an automorphism of the form
\equa{VAR}{\si: (x,y)\mapsto (x^{\frac{q}{q-p}},
x^{-\frac{p}{q-p}}y+\l_1x^{-\frac{p_1q}{q_1(q-p)}}+\cdots+\l_sx^{-\frac{p_s q}{q_s(q-p)}}),}
for some $0\ne \l_i\in\C$, $ p,q,p_i,q_i\in\N$
with $\frac{p_i}{q_i}<\frac{p_{i+1}}{q_{i+1}}<\frac{p}{q}<1$ for $1\le i<s$.
\item[\rm(2)]
For any $i,j\in\Q$, we have
 $\tr\,(x^iy^j\bar F)=0$ if and only if
\equa{AMSM}{\tr\,H=0,\mbox{ where }H=
x^{\frac{qi-pj}{q-p}}(y+\l_1x^{\frac{q_1p-p_1q}{q_1(q-p)}}+\cdots+\l_sx^{\frac{q_sp-p_s q}{q_s(q-p)}})^jF.}
Note that since $\bar F,\bar G\in\C[x,y]$, we clearly have
 $\tr\,(x^iy^j\bar F)=0$ if $(i,j)\notin\Z_-^2$
$($we wish that one may obtain some condition on $F$ from \eqref{AMSM}$)$.
\end{itemize}\end{coro}

We can decompose $F$ and $G$ as sums of $\dot\a$-type q.h.e.~(cf.~Definition \ref{comp}(3)\VS{-7pt})
\begin{eqnarray}\label{fi-defoii}&&
F=\mbox{$\sum\limits_{i\in\Q_+}$} \dot f^i,\mbox{ where
 }\dot f^i=x^{\frac{m}{m+n}}y^{m-i}\mbox{$\sum\limits_{j={\rm max}\{0,i-m\}}^{\infty}$} f_{ij} x^{-\dot\a j}y^{j}\mbox{ for some $f_{ij}\in\C$},\\[-3pt]
\label{gi-defoii}&&
G=\mbox{$\sum\limits_{i\in\Q_+}$} \dot g^i,\mbox{ where
 }\dot g^i=x^{\frac{n}{m+n}}y^{n-i}\mbox{$\sum\limits_{j={\rm max}\{0,i-n\}}^{\infty}$} g_{ij} x^{-\dot\a j}y^{j}\mbox{ for some $g_{ij}\in\C$},
\end{eqnarray}
where all sums are finite. We  call $\dot\a$ the {\it leading degree} of $F$ and $G$.

{}From \eqref{F-now-is} and \eqref{G-now-is}, we see that $(x^{\frac{m}{m+n}}f,x^{\frac{n}{m+n}}g)$, being obtained from the pair $(\Dhat F,\Dhat G)$ under the mapping \eqref{x-yyyyy}, is a \QJ  pair,
 \VS{-5pt}i.e.,
\equa{f-g(y)--}{\dis\frac{m}{m+n}fg'-\frac{n}{m+n}f'g=J,
}
for some $J 
\in\C\bs\{0\}$, where the prime stands for the derivative $\frac{d}{dy}$.
Thus $f$ and $g$ do not have common irreducible factors and all irreducible factors of $f$ and $g$ have multiplicity $1$.

\section{Poisson algebras}\label{sect-remark}
In this section, we use the natural Poisson algebra structure on $\C[y]((x^{-\frac1N}))$ to discuss Jacobi pairs.
The main results of this section are Theorem \ref{MMM-lemm} and
Corollary \ref{Rema4.1}. 

%
%
\subsection{Poisson  algebra $\C[y]((x^{-\frac1N}))$ and exponential operator $e^{\ad_H}$}\label{Subs4.1}
Let $\PP:=\C[y]((x^{-\frac1N}))$ (cf.~notation in \eqref{ring-ABC1--BB}), where $N$ is some fixed sufficient large integer (such that all elements considered below are in $\PP$). By Definition \ref{gen-JP}(1), $\PP$ is a poisson algebra.
For $H\in\PP$, we use $\ad_H$ to denote the operator on $\PP$ such that \equa{ad-H}{\mbox{$\ad_H(P)=[H,P]$ \ for \ $P\in\PP$.}} If $H$ has the form \equa{===H===}{\mbox{$H= x(a_0y+a_1)+\widetilde{H}$, where $\deg_x\widetilde{H}<1,\,a_0,a_1\in\C$,}} then we can
define the following exponential operator on $\PP$\VS{-6pt}:
\equa{exppp}{e^{\ad_H}:=\sum\limits_{i=0}^\infty{\dis\frac{1}{i!}\ad^i_H,}}
which is well defined in $\PP$, and is in fact an automorphism of the Poisson algebra $(\PP,[\cdot,\cdot],\cdot)$ with inverse $e^{\ad_{-H}}$, i.e.\VS{-5pt}, \equa{Ead-H---}{\mbox{$(e^{\ad_H})^{-1}\!=\!e^{\ad_{-H}},\ e^{\ad_H}(PQ)\!=\!e^{\ad_H}(P)e^{\ad_H}(Q),\
e^{\ad_H}([P,Q])\!=\![e^{\ad_H}(P),e^{\ad_H}(Q)]$,}}  for $P,Q\in\PP$. For any
 $H_1,H_2\in\PP$ having the form \eqref{===H===}, one can \VS{-3pt}verify
\equa{Exp-law}{e^{\ad_{H_1}}\cdot e^{\ad_{H_2}}=e^{\ad_{K}}\cdot e^{\ad_{H_1}},\mbox{ \ where $K=e^{\ad_{H_1}}(H_2)$.
}}
Note that for any $H_i\in x^{1-\frac{i}{N}}\C[y][[x^{-\frac1N}]]$, $i=1,2,...$, both \VS{-3pt}operators
\equa{Both-op}{\Plar e^{\ad_{H_i}}:=\cdots e^{\ad_{H_2}}e^{\ad_{H_1}},\ \ \ \
\Prar e^{\ad_{H_i}}:=e^{\ad_{H_1}}e^{\ad_{H_2}}\cdots, }
are well defined on $\PP$: For any giving $P\in x^{\frac{a}{N}}\C[y][[x^{-\frac1N}]]$ for some $a\in\Z$, to compute, say,   $Q=\Prar e^{\ad_{H_i}}(P)$\VS{-3pt}, we can write $Q$ as $Q=\sum_{j=0}^\infty x^{\frac{a-j}{N}}Q_j$ with $Q_j\in\C[y]$. Then the computation of $Q_j$ for each $j$ only involves finite numbers of operators: $e^{\ad_{H_i}},\,i\le j$. This together with \eqref{Exp-law} also shows that for any $H_i$ as above, we can
find some $K_i\in x^{1-\frac{i}{N}}\C[y][[x^{-\frac1N}]]$ such \VS{-3pt}that \equa{lar-rar}{\Plar e^{\ad_{H_i}}
=\Prar e^{\ad_{K_i}}.}
\begin{theo}\label{MMM-lemm} Suppose $(F,G)$ is a \QJ  pair in $\PP$ satisfying 
\eqref{F-now-is} and \eqref{G-now-is}. There \VS{-8pt}exist \equa{HK=======}{H=x+\sum\limits_{i=1}^\infty x^{1-\frac{i}{N}}h_i,\ \
\ K=y+\sum\limits_{i=1}^\infty x^{-\frac{i}{N}}k_i\in\PP\mbox{ \ \ \ with \ \ }h_i,k_i\in\C[y],} such that $[H,K]=1$ \VS{-7pt}and \equa{F------G}{F=H^{\frac{m}{m+n}}f(K),\ \ \ \ \ G=H^{\frac{n}{m+n}}g(K).}
\end{theo}
{\it Proof.~}~Let $F,G\in\PP$ be as in
\eqref{F-now-is} and \eqref{G-now-is}. Let $i$ be the smallest positive integer such that $(f_{i},g_{i})\ne(0,0)$. By computing the terms with $x$-degree $-\frac{i}{N}$ in $[F,G]=J$, we \VS{-7pt}obtain
\equa{G-F-i0}{\dis\frac{m}{m+n}fg'_{i}-(\frac{n}{m+n}-\frac{i}{N})f'g_i+
(\frac{m}{m+n}-\frac{i}{N})f_ig'-\frac{n}{m+n}f'_ig=0,}
where the prime stands for the derivative $\frac d{dy}$. First assume $i\ne N$. Taking $\bar h_i=\frac1J (f_ig'-f'g_i)$, $\bar k_i=\frac{1}{(m+n)J}(mfg_i-nf_ig)$, using \eqref{f-g(y)--} and \eqref{G-F-i0}, we \VS{-7pt}have
\equa{f-i-g-i=}{\dis\Big(1-\frac{i}{N}\Big)\bar h_i+\bar k'_i=0,\ \ \ f_i=\frac{m}{m+n}\bar h_if+f'\bar k_i,\ \ \ g_i=\frac{n}{m+n}\bar h_ig+g'\bar k_i.}
Take $Q_i=\frac{-N}{N-i}x^{1-\frac{i}{N}}\bar k_i$. \VS{-7pt}Then \equa{FFFGGG-iii}{\mbox{ $x^{\frac{m}{m+n}-\frac{i}{N}}f_i=-[Q_i,x^{\frac{m}{m+n}}f],$ \ \ \ \ \ $x^{\frac{n}{m+n}-\frac{i}{N}}g_i=-[Q_i,x^{\frac{n}{m+n}}g].$}}
Thus if we apply the automorphism $e^{\ad_{Q_i}}$ to $F$ and $G$, we can suppose $f_i=g_i=0$.

Now assume $i=N$. Then \eqref{G-F-i0} gives $mfg_i-nf_ig=(m+n)cJ$ for some $c\in\C$. This together with \eqref{f-g(y)--} implies
$m(g_i-cg')f=n(f_i-cf')g$. Since $f$ and $g$ are coprime (cf.~the statement after \eqref{f-g(y)--}),
we have $g|(g_i-cg')$ in $\C[y]$, i.e.,
there exists $\tilde k_i\in\C[y]$ such \VS{-7pt}that\equa{i====N}{\mbox{
$\dis g_i-cg'=\frac{n}{m+n}g \tilde k_i,\mbox{ \ \ thus, \ }f_i-cf'=\frac{m}{m+n}f\tilde k_i$.}}
Now if we first apply $e^{\ad_{Q_i}}$ to $F,G$ (with $Q_i=\int\tilde k_idy$), and
then applying the \VS{-7pt}automorphism \equa{tau-c}{\mbox{$\tau_c:(x,y)\mapsto(x,y-cx^{-1})$,}} we can suppose $f_i=g_i=0$.

The above shows that by applying infinite many automorphisms $e^{\ad_{Q_i}},\tau_c$, $i\ge1$, the pair $(F,G)$ becomes
$(x^{\frac{m}{m+n}}f,x^{\frac{n}{m+n}}g)$. Thus $F=\si(x^{\frac{m}{m+n}}f)=H^{\frac{m}{m+n}}f(K)$, $G=\si(x^{\frac{n}{m+n}}g)=H^{\frac{n}{m+n}}g(K)$ for some automorphism $\si$ of the form
(cf.~\eqref{Ead-H---} and \eqref{lar-rar})
\equa{PAAAAAA}{\si:=\ \cdots e^{\ad_{P_j}}\cdots e^{\ad_{P_{N+1}}}\tau_{-c\,}e^{\ad_{P_{N}}}\cdots e^{\ad_{P_1}}=(\cdots e^{\ad_{Q_j}}\cdots e^{\ad_{Q_{N+1}}}\tau_{c\,} e^{\ad_{Q_{N}}}\cdots e^{\ad_{Q_1}})^{-1},} for some $P_i\in x^{1-\frac{i}{N}}\C[y][[x^{-\frac1N}]]$, where $H=\si(x),$ $K=\si(y)$ have the forms as in the theorem since $\deg_xQ_i<1$ for all $i\ge1$.\hfill$\Box$\vskip5pt

Theorem \ref{MMM-lemm} can be generalized as follows (one may wish to obtain some information on $F,G$ from this result).
\begin{coro}\label{Rema4.1}
For any \QJ pair  $(F,G)$ in $\C[x^{\pm\frac1N},y]$, 
and any line $L$ which meet the boarder of $\Supp\,F$ $($either meet an edge of $\Supp\,F$ or a vertex of $\Supp\,F)$ such that $L$ does not pass the origin. Regarding $L$ as the prime line, one can start from $F\comp0=F_L$, the part of $F$ corresponding to $L$,  to obtain that there exists an automorphism $\si$ of $\DD=\C[y]((x^{-\frac1N}))$ of the form \eqref{PAAAAAA}
 $($we need to change $\DD$ to $\DD=\C[y]((x^{\frac1N}))$ if $L$ is an edge at the left side of $\Supp\,F)$
with $P_i,Q_i$ in the localized ring $\DD[F_{[0]}^{-1}]$
such that $\si(F)=F\comp0$ and $\si(G)=\phi(F\comp0)+R_0$, where $R_0$ is as in Lemma $\ref{quasi-J}$, and $\phi(F\comp0)$ is a function of $F\comp0$ of the form $\phi(F\comp0)=\sum_{i\in\Q_+}a_iF_{[0]}^i$ with $a_i\in\C$ and $F_{[0]}^i$ is a rational function if $a_i\ne0$. 
\end{coro}

%
\subsection{Jacobi pairs in $\C[y]((x^{-\frac1N}))$ satisfying
\eqref{F-now-is}--\eqref{f0===}}
Now let $F,G\!\in\!\C[x^{\pm\frac1N},y]$ be a Jacobi pair satisfying
\eqref{F-now-is}--\eqref{f0===}.
Let $H,K\!\in\!\PP$ be as in Theorem \ref{MMM-lemm} (note that $H,K$ are not necessarily in $\C[x^{\pm\frac1N},y]$).
\begin{lemm}\label{lemm-H-K===}
There exist a unique $\a$, 
called the \textbf{leading degree} of $H$ and $K$ $($as in \eqref{fi-defoii}$)$,
such that $H=\sum_{i\in\frac1{\b}\Z_+} H\UP{i},$ $K=\sum_{\in\frac1{\b}\Z_+} K\UP{i}$ for some $\b\in\N$, where $H\UP{i},K\UP{i}$ are respectively the $\a$-th $\,-i{\ssc\,}$-th components of $H,K$ of the following forms $($cf.~Definition $\ref{comp}(1))$\VS{-7pt},
\equa{H-K=====}{
H\UP{i}=xy^{-i}\sum\limits_{j\ge i}h_{ij} x^{-j\a}y^j,\ \ \ K\UP{i}=y^{1-i}\sum\limits_{j\ge{\rm max}\{0,i-1\}}k_{ij}x^{-j\a}y^j,}
for some $h_{ij},k_{ij}\in\C$ with $k_{10}=0$ $($note that $K$ does not contain the constant term by \eqref{HK=======}$)$ and at least some $h_{0j}$ or $k_{0j}$ is nonzero for some $j\ge1$.
\end{lemm}

\begin{rema}\rm\label{rema4.1}
We remark here that when some negative power of $y$ appears in an expression (as in \eqref{H-K=====}), we always regard the element as in a proper space which is a subspace of the \VS{-5pt}space
\equa{LARGE-s}{
\widetilde{\PP}:=\C((y^{-1}))((x^{-\frac1N}))\mbox{ (cf.~notation in \eqref{ring-ABC1--BB})}.}
\end{rema}\noindent{\it Proof of Lemma \ref{lemm-H-K===}.~}~For $i\ge0$, we inductively define $Q_i$,
 $F^{(i)}\!=\!x^{\frac{m}{m+n}}(f\!+\!\sum_{j=1}^\infty x^{-\frac{j}{N}}f^{(i)}_j)$ and
 $G^{(i)}\!=\!x^{\frac{m}{m+n}}(g\!+\!\sum_{j=1}^\infty x^{-\frac{j}{N}}g^{(i)}_j)$ as follows: $Q_0\!=\!0,$ $F^{(0)}\!=\!F$,
 $G^{(0)}\!=\!G$ (in particular, $f^{(0)}_j\!=\!f_j$, $g^{(0)}_j\!=\!g_j$). For $i\!\ge\!1$, $Q_{i}$ is
 defined as in the proof of Theorem \ref{MMM-lemm} with $(F,G)$ replaced by
  $(F^{(i-1)},G^{(i-1)})$,
 and \VS{-7pt}set
 \equan{NewF-i===}{\mbox{
$F^{(i)}=e^{\ad_{Q_i}}(F^{(i-1)}),$
\ \ $G^{(i)}=e^{\ad_{Q_i}}(G^{(i-1)})$.
}}
We remark that since our purpose here is to prove 
\eqref{H-K=====}, from the discussions below, we see it does not matter whether or not $\tau_c$ defined in
\eqref{tau-c} is involved in \eqref{PAAAAAA}. Thus for convenience, we may assume $\tau_c$ is not involved. Alternatively, one can also
formally regard the automorphism $\tau_c$ as $e^{\ad_{\ol{Q}_N}}$ with $\ol{Q}_N={\rm ln}\,x$ and regard $\ol{Q}_N$ as an element
with $x$-degree and $y$-degree being zero, and $\ptl_x\ol{Q}_N=x^{-1}$, $\ptl_y\ol{Q}_N=0$. In this way,
 $\ol{Q}_N$ can be regarded as another $Q_N$, and we define
 $\ol{F}^{(N)}=e^{\ad_{\ol{Q}_N}}(F^{(N)}),$
 $\ol{G}^{(N)}=e^{\ad_{\ol{Q}_N}}(G^{(N)})$ and
 $F^{(N+1)}=e^{\ad_{Q_{N+1}}}(\ol{F}^{(N)}),$
 $G^{(N+1)}=e^{\ad_{Q_{N+1}}}(\ol{G}^{(N)})$, etc.

 %

Now we choose $\gamma$ to be the minimal rational number such that (thus at least one equality holds below\VS{-7pt})
 \equa{b-ge-0}{\deg_yQ_{i}\le1+i\gamma,\ \deg_yf_i\le m+i\gamma ,\ \deg_yg_i\le n+i\gamma \mbox{ \ \ for }0\le i\le N+M_1+M_2,}
 where $M_1,M_2$ are as in \eqref{F-now-is}, \eqref{G-now-is}. Note from the edge $\dot L$ in \eqref{new-supp-FCOM-region}, we see that $\gamma\ge\frac{1}{N\dot\a}>0$. We claim that for all $i\ge0$, we \VS{-7pt}have
\equa{WEHAVE}{\mbox{$\deg_yQ_i\le1+i\gamma$, \ \ $\deg_yf^{(i)}_j\le m+j\gamma,\ \ \ \deg_yg^{(i)}_j\le n+j\gamma$ \ for all $j\ge0$.}}
By \eqref{b-ge-0}, the claim holds for $i=0$ (note that $f^{(0)}_j=g^{(0)}_j=0$ if $j>M_1+M_2$ by \eqref{F-now-is} and \eqref{G-now-is}).
 Inductively assume that \eqref{WEHAVE} holds for some $i=i_0-1\ge0$. Now assume $i=i_0$.
We want to \VS{-7pt}prove \equa{Q-i0==1+}{\deg_yQ_{i_0}\le1+i_0\gamma.}
If $i_0\le N$, we already have \eqref{Q-i0==1+} by \eqref{b-ge-0}. Assume $i_0>N$. By the inductive assumption, $\deg_yf^{(i_0-1)}_{i_0}\le m+i_0\gamma$. Then the first equation of
  \eqref{FFFGGG-iii} with $i=i_0-1$ 
  shows that either
  $\deg_yf^{(i_0-1)}_{i_0}=\deg_yQ_{i_0}+\deg_yf-1$ or else
  $(1-\frac{i_0}{N},\deg_yQ_{i_0})=\b(\frac{m}{m+n},\deg_yf)$
   for some $\b\in\Q$. The later case cannot occur since $i_0>N$ (i.e., $1-\frac{i_0}{N}<0$). The first case
implies \eqref{Q-i0==1+}. This proves \eqref{Q-i0==1+} in any case.
Now by definition of $e^{\ad_{Q_{i_0}}}$ (cf.~\eqref{exppp}),
we \VS{-7pt}know \equa{f-i00000}{\mbox{$x^{\frac{m}{m+n}-\frac{j}{N}}f^{(i_0)}_j=\ $ a combination of elements of forms $\ad^{j_1}_{Q_{i_0}}(x^{\frac{m}{m+n}-\frac{j_2}{N}}f^{(i_0-1)}_{j_2})$, }}
for $j_1,j_2\in\Z_+$ satisfying $j_1i_0+j_2=j$.
Thus for all $j$\VS{-7pt},
\equan{f-i0-j}{\deg_yf^{(i_0)}_{j}\le \max\big\{j_1(\deg_yQ_{i_0}-1)+\deg_yf^{(i_0-1)}_{j_2}\,
\big|\,j_1,j_2\in\Z_+,\,j_1i_0+j_2=j\big\}\le m+j\g,
}
where the last inequality is obtained by the inductive assumption.
 Analogously,   $\deg_yg^{(i_0)}_{j}\le n+j\g$.
%
This completes the proof of \eqref{WEHAVE}.
Now let $i_1\ge1$ be minimal such that when $i=i_1$, at least one equality holds in \eqref{b-ge-0}. We \VS{-10pt}claim \equa{Q-i-0}{\deg_yQ_{i_1}=1+i_1\gamma.} Otherwise $\deg_yQ_{i_1}<1+i_1\gamma$ and, say,
$\deg_yf_{i_1}= m+i_1\gamma$. We want to prove  by induction on $\ell$\VS{-7pt},
 \equa{j-i-1=====}{\mbox{$\deg_yf^{(\ell{\ssc\,})}_{i_1}= m+i_1\gamma$ \
for $0\le \ell\le i_1-1$.}}
By definition, \eqref{j-i-1=====} holds for $\ell=0$. Inductively assume \eqref{j-i-1=====} holds for $\ell-1<i_1-1$. Similar to \eqref{f-i00000},
$x^{\frac{m}{m+n}-\frac{i_1}{N}}f^{(\ell{\ssc\,})}_{i_1}$\VS{-3pt} is a combination of $\ad^{j_1}_{Q_{\ell}}(x^{\frac{m}{m+n}-\frac{j_2}{N}}f^{(\ell-1)}_{j_2})$ (for $j_1,j_2\in\Z_+$ with $j_1\ell+j_2=i_1$), whose $y$-degree is either $<m+i_1\gamma$ if $j_1\ne0$, or else $=m+i_1\g$ if $j_1=0$ (note that \VS{-3pt}$x^{\frac{m}{m+n}-\frac{i_1}{N}}f^{(\ell-1)}_{i_1}$ does indeed appear as a term in $x^{\frac{m}{m+n}-\frac{i_1}{N}}f^{(\ell{\ssc\,})}_{i_1}$). Thus \eqref{j-i-1=====} holds. However, \eqref{FFFGGG-iii} with $(f_i,Q_i)$ replaced by $(f^{(i_1-1)}_{i_1},Q_{i_1})$ implies that $\deg_yf^{(i_1-1)}_{i_1}\le \deg_yQ_{i_1}-1+\deg_yf<m+i_1\gamma$, a contradiction with \eqref{j-i-1=====} (with $\ell=i_1-1$). This proves \eqref{Q-i-0}.

Using \eqref{Ead-H---} and \eqref{Exp-law}, we can explicitly determine $P_i$ in terms of $Q_1,...,Q_i,\tau_c$ from \eqref{PAAAAAA}; for instance\VS{-5pt},$$\mbox{
$P_1=-\tau_c(Q_1)$, \ $P_2=-\tau_c(e^{\ad_{P_1}}(Q_2))$, \  $P_3=-\tau_c(e^{\ad_{P_2}}e^{\ad_{P_1}}(Q_3))$, ...}$$ In particular, if we write $P_i=\sum_{j=i}^\infty x^{1-\frac{j}{N}}p_{ij}$ for some $p_{ij}\in\C[y]$, using \eqref{WEHAVE} and \eqref{Q-i-0},
one can show that
$\deg_yp_{ij}\le1+j\gamma$ for all $i,j$ and $\deg_yp_{ij}<1+j\gamma$ if $j<i_1$, and furthermore,
$i_1$ is the minimal integer such that the equality $\deg_yp_{i,i_1}=1+i_1\gamma$ holds when $i=i_1$.
This implies (using $H=\si(x),$ $K=\si(y)$ and definitions of $h_i,k_i$ in \eqref{HK=======}, as discussions above\VS{-7pt})
$$\deg_yh_i\le i\gamma,\ \deg_yk_i\le1+i\gamma\mbox{ \ for }i\ge1,\mbox{ \ and the equalities hold when $i=i_1$}\VS{-5pt}.$$Thus we have \eqref{H-K=====} by choosing $\a=\frac1{\g N}$.\hfill$\Box$\vskip5pt

Clearly, we have $\a\le\dot\a$ (otherwise from \eqref{F------G}, we would obtain $\dot f^0=x^{\frac{m}{m+n}}y^m$, a contradiction with \eqref{f0===}).
We rewrite $F=\sum_{i=0}^\infty F\UP{i}$ according to the leading degree $\a$ of $H$ and $K$ (but not according to the leading degree $\dot \a$ of $F$ and $G$; in particular if $\a=\dot\a$ then $F\UP{i}=\dot f^i$, cf.~\eqref{fi-defoii}), and call $F\UP{i}$ the {\it $\a$-type $\,-i{\ssc\,}$-th component} of $F$. Then \VS{-7pt}(cf.~\eqref{f0===})
\equa{FUP==}{H^{\frac{m}{m+n}}\UPo0 K^m\UPo0=F\UPo0=\left\{\begin{array}{ll}x^{\frac{m}{m+n}}y^m\prod\limits_{i=1}^\de(a_ix^{-\dot\a}y+1)^{m_i}\!\!&\mbox{if \ }\a=\dot\a,\\[5pt]
x^{\frac{m}{m+n}}y^m&\mbox{if \ }\a<\dot\a.\end{array}\right.}
Note that $(H^{\frac{m}{m+n}}\UPo0 K^m\UPo0,H^{\frac{n}{m+n}}\UPo0 K^{1-m}\UPo0)$ is a \QJ  pair, also
$(F\UPo0,F^{-1}\UPo0P\UPo0)$ is a \QJ  pair, \VS{-7pt}where $P\UPo0$ has the form
 (cf.~\eqref{NOWOO-R}, \eqref{f0===}, Lemmas \ref{quasi-J} and \ref{irr-factor})
\equa{PUP==}{P\UPo0=\left\{\begin{array}{ll}xy\prod\limits_{i=1}^{\de+\de'}(a_ix^{-\dot\a}y+1)&\mbox{if \ }\a=\dot\a,\\[2pt]
xy&\mbox{if \ }\a<\dot\a,\end{array}\right.}
for some $\de'\ge0$ and $a_i\in\C$.
Thus up to a nonzero scalar, $H^{\frac{n}{m+n}}\UPo0 K^{1-m}\UPo0$ \VS{-3pt}must have the form $F^{-1}\UPo0P\UPo0+\psi$ with $[F\UPo0,\psi]=0$, and so $\psi$ is a function on $F\UPo0$ (noting from \eqref{equa2.5} that in case $[G,F]=0$, we can obtain that each $b_i$ does not depend on $x$, i.e., $G$ is a function on $F$). Since $H^{\frac{n}{m+n}}\UPo0 K^{1-m}\UPo0$ is an $\a$-type q.h.e.~(cf.~Definition \ref{comp}(3)), $\psi$ must be of the form $\theta_0F\UPo0^\mu$ for $\theta_0\in\C$ and $\mu\in\Q$ (we do not need to know the exact value of $\mu$, however by noting
 that $\Supp(F^{-1}\UPo0P\UPo0+\psi)$ is on a line with slope $-\frac1{\a}$, the term $x^{\frac{\mu m}{m+n}}y^{\mu m}$ in $F\UPo0^\mu$  must be of the form $x^{\frac{n}{m+n}-i_0\a}y^{1-m+i_0}$ for some $i_0\in\N$
since $F^{-1}\UPo0P\UPo0+\psi$ has a term $x^{\frac{n}{m+n}}y^{1-m}$, one can compute $\mu=\frac{n+(m+n)(1-m)\a}{m(1+(m+n)\a)}$).
Thus we can solve (up to nonzero scalars\VS{-7pt})
\begin{eqnarray}\label{H-K-UP0=}&&
H\UPo0=F\UPo0^{-\frac{m+n}{m(m+n-1)}}(P\UPo0+\theta_0F\UPo0^{\mu+1})^{\frac{m+n}{m+n-1}},\nonumber\\[-3pt]
&&K\UPo0=F\UPo0^{\frac{m+n}{m(m+n-1)}}(P\UPo0+\theta_0F\UPo0^{\mu+1})^{-\frac{1}{m+n-1}},\  \ \ \ \ H\UPo0K\UPo0=P\UPo0+\theta_0F\UPo0^{\mu+1}.\end{eqnarray}
In particular, if $\a<\dot\a$, then $\theta_0\ne0$ (otherwise $H\UPo0=x,$ $K\UPo0=y$, a contradiction with the definition of $\a$).
Applying $\ptl_x,\ptl_y$ to \eqref{F------G}, using \eqref{f-g(y)--}, we \VS{-8pt}have
\equa{-ptl-x-K}{\dis H\ptl_xK=\frac{1}{(m+n)J}(mF\ptl_xG-nG\ptl_xF),\ \ \ H\ptl_yK=\frac{1}{(m+n)J}(mF\ptl_yG-nG\ptl_yF),}
which imply that they are {\it polynomials} (i.e., elements in $\C[x^{\pm\frac1N},y]$).
In particular, by \eqref{H-K=====}, $H\UPo0K\UPo0= yH\UPo0\ptl_yK\UPo0+\frac1{\a}xH\UPo0\ptl_xK\UPo0$ is a polynomial. Thus the last equation of \eqref{H-K-UP0=} shows that if $\theta_0\ne0$, then $F\UPo0^{\mu+1}$ is a polynomial (thus $(\mu+1) m'\in\Z_+$, where $m'$ is the greatest common divisor of $m, m_i,\,i=1,...,\de$, cf.~\eqref{FUP==}). Note that in any case (either $\a=\dot\a$ or $\a<\dot\a$), each irreducible factor (in $\C[x^{\pm\frac1N},y]$) of $F\UPo0$ is an irreducible factor of $H\UPo0K\UPo0$ with multiplicity $1$ by \eqref{FUP==} and \eqref{PUP==}.
Also note that since $H\UP{i},K\UP{i}$ are $\a$-type q.h.e.~(cf.~Definition \ref{comp}(3)) of the form \eqref{H-K=====}, we always \VS{-7pt}have
\equa{-p-x-py==}{\!
x\ptl_xH\UP{i}\!+\!\a y\ptl_yH\UP{i}\!=\!(1\!-\!i\a)H\UP{i},\
x\ptl_xK\UP{i}\!+\!\a y\ptl_yK\UP{i}\!=\!(1\!-\!i)\a K\UP{i}\mbox{ for all }i\!\in\!\frac1\b\Z_+.\!}
In particular, we can \VS{-7pt}obtain\equa{H0-K0Ja}{\mbox{$\dis1=[H\UPo0,K\UPo0]
=\frac1{\a y}(\a K\UPo0\ptl_xH\UPo0-H\UPo0\ptl_xK\UPo0)$,}}
where the first equality follows by comparing the $\a$-type $0$-th components in $1=[H,K]$. We always \VS{-7pt}denote
\equa{nu-0===}{\dis\nu_0=1+\frac1{\a}.}

Now \VS{-2pt}let \equa{nu====0}{\mbox{$\nu>0$ be smallest such that $(H\UP{\nu},K\UP{\nu})\ne(0,0)$.}}
Then \eqref{-ptl-x-K} gives that (using \eqref{product-HK} and \eqref{-p-x-py==}\VS{-7pt})
\begin{eqnarray}\label{HKUP-HK}\!\!\!\!\!\!\!\!\!R_1\!\!\!&:=&\!\!\!(1-\nu)H\UPo0K\UP{\nu}+H\UP{\nu}K\UPo0\nonumber\\
\!\!\!\!\!\!\!\!\!\!\!\!&=&\!\!\!y(H\UPo0\ptl_yK\UP{\nu}\!+\!H\UP{\nu}\ptl_yK\UPo0)
\!+\!\frac{x}{\a}(H\UPo0\ptl_xK\UP{\nu}\!+\!H\UP{\nu}\ptl_xK\UPo0)\mbox{ is a polynomial}\end{eqnarray}
\begin{lemm}\label{11111Qnu==}
There exists a unique \VS{-7pt}element \equa{Q-nu===}{
Q_\nu=\sum\limits_{j\ge\nu-1}q_{\nu j}x^{1-j\a}y^{1+j-\nu}\in\PP\mbox{ \ for some $q_{\nu j}\in\C$,
 and }\deg_xQ_\nu<1,}
 with $q_{10}=0$ $($if $\nu=1)$ and $q_{\nu_0,\nu_0-1}=0$ $($if $\nu_0\in\N)$, such \VS{-5pt}that
\begin{eqnarray}\label{H-K-UP-nu1}&&
H\UP{\nu}=[Q_\nu,H\UPo0]=\frac{1}{\a y}\big(H\UPo0\ptl_xQ_\nu-(1+\a(1-\nu))Q_\nu \ptl_xH\UPo0\big),
\\ \label{H-K-UP-nu2}&&
K\UP{\nu}=[Q_\nu,K\UPo0]=\frac{1}{\a y}\big(\a K\UPo0\ptl_xQ_\nu-(1+\a(1-\nu))Q_\nu \ptl_xK\UPo0\big).\end{eqnarray}
%
\end{lemm}\vskip-9pt

\noindent{\it Proof.~}~
%
First we set $q_{\nu,\nu-1}=0$ and inductively choose unique $q_{\nu,j}\in\C$ for $j\ge\nu$ to \VS{-4pt}satisfy
$$(1+j-\nu)q_{\nu j}+\mbox{$\sum\limits_{i\ge1}$}\big((1-i\a)(1+j-i-\nu)-i(1-(j-i)\a)\big)h_{0i}q_{\nu,j-i}
=h_{\nu,j},\ \ j\ge\nu\VS{-7pt},$$
which implies $[Q_\nu,H\UPo0]=H\UP{\nu}$.
Set $\ol{K}_\nu=K\UP{\nu}-[Q_\nu,K\UPo0]$. Then by \eqref{nu====0}, we \VS{-5pt}have \equa{1====[H,K]}
{1=[H,K]=[H\UPo0,K\UPo0]+[[Q_\nu,H\UPo0],K\UPo0]+[H\UPo0,[Q_\nu,K\UPo0]+\ol{K}_\nu]+\cdots,}
which implies $[H\UPo0,\ol{K}_\nu]=0$ since the omitted terms  are those components whose component indices are $<-\nu$.
Thus as in the arguments after \eqref{PUP==}, $\ol{K}_\nu=\l_\nu H\UPo0^{a_\nu}$ for some
$\l_\nu\!\in\!\C$ and $a_\nu\!\in\!\Q$.
Since $\Supp\,\ol{K}_\nu$ and $\Supp\,K\UP{\nu}$ are on the same line, we have $x^{a_\nu}=x^{-i_0\a}y^{1-\nu+i_0}$ for some $i_0$, thus, $a_\nu=(1-\nu)\a$. First assume $\nu\ne\nu_0$ (cf.~\eqref{nu-0===}). Then
$K\UP{\nu}=[Q_\nu+\frac{\l_\nu}{1+\a(1-\nu)}H\UPo0^{1+\a(1-\nu)},K\UPo0]$. Thus by re-denoting $Q_\nu+\frac{\l_\nu}{1+\a(1-\nu)}H\UPo0^{1+\a(1-\nu)}$ to be $Q_\nu$, we can suppose $\l_\nu=0$
(note that we still have $\deg_xQ_\nu<1$ since if $\nu>1$ then $\deg_xH\UPo0^{1+\a(1-\nu)}<1$, and if $\nu=1$ then $H\UPo0$ contains the term $x$, but since $K$ does not contain the constant term, i.e., $\l_\nu=0$).

Now assume $\nu=\nu_0$. 
Then $a_{\nu_0}=-1$. Noting from \eqref{H-K=====} that
among all components of $K$, only $K\UP{\nu_0}$ can possibly contain the term $x^{-1}$, also we can deduce from
\eqref{H-K=====} and \eqref{Q-nu===} that $[Q_{\nu_0},K\UPo0]=\frac1y K\UPo0\ptl_xQ_{\nu_0}$ cannot contain the term $x^{-1}$. Thus $\Coeff(K,x^{-1})=\Coeff(\ol{K}_\nu,x^{-1})=\Coeff(\l_\nu H\UPo0^{-1},x^{-1})=\l_\nu$.
If necessary by replacing $y$ by $y-\l x^{-1}$ for some $\l\in\C$, we can always suppose $\Coeff(K,x^{-1})=0,$
i.e., $\l_\nu=0$ (noting that since the original pair $(F,G)$ satisfies \eqref{Res-x-y=0}, and \eqref{Res-x-y=0} is also satisfied by the \QJ  pair
$(x^{\frac{m}{m+n}}f,x^{\frac{n}{m+n}}g)$, by Theorem \ref{Jacobian-el}(1) and \eqref{F------G}, we see $\Res_x(H\ptl_xK)=0$, from this we in fact have
$\Coeff(K,x^{-1})=0$).
Hence in any case, we have
the first equalities of
\eqref{H-K-UP-nu1} and \eqref{H-K-UP-nu2}, and the above proof also shows that such $Q_\nu$ is unique.
Now using \eqref{Q-nu===} we obtain as in \eqref{-p-x-py==}\VS{-7pt},
\equa{ptl-y-Q-nu}{\dis\ptl_yQ_\nu=\frac1{\a y}((1+\a(1-\nu))Q_\nu-x\ptl_xQ_\nu),}
from this and \eqref{-p-x-py==},
we have
the second equalities of
\eqref{H-K-UP-nu1} and \eqref{H-K-UP-nu2}.\hfill$\Box$\vskip5pt
Computing the $\a$-type $\,-\nu{\ssc\,}$-th component $F\UP{\nu}$ of $F$ in \eqref{F------G} gives that
(cf.~\eqref{F-now-is} and \eqref{nu====0}\VS{-4pt})
\equa{AHComp}{\dis R_2:=H\UPo0^{\frac{m}{m+n}}K\UPo0^m\Big(\frac{m H\UP{\nu}}{(m+n)H\UPo0}+\frac{mK\UP{\nu}}{K\UPo0}+c_\nu K\UPo0^{-\nu}\Big)\mbox{ is a polynomial.}}\VS{-19pt}\PAR\NOindent
Now we first assume $\nu\ne\nu_0$.
Using \eqref{FUP==}, \eqref{H-K-UP0=}, \eqref{H-K-UP-nu1}, \eqref{H-K-UP-nu2} and \eqref{ptl-y-Q-nu},
we obtain that  $\frac{\a }{1+\a(1-\nu)}R_1$ and
$\frac{(m+n)\a }{m}R_2H\UPo0^{1-\frac{m}{m+n}}K\UPo0^{1-m}$ are \VS{-5pt} respectively equal \VS{-7pt}to
\begin{eqnarray}\label{RATTTTT}
\!\!\!\!\!\!\!\!\!\!\!\!\!\!\!\!\!\!\!\!\!\!\!\!\!\!\!\!\!\!\!\!\!&\!\!\!\!\!\!R_3\!\!\!&:=
\frac1yH\UPo0K\UPo0\ptl_xQ_\nu-\frac1y\big((1-\nu)H\UPo0\ptl_xK\UPo0+K\UPo0\ptl_xH\UPo0\big)Q_\nu,\\
\!\!\!\!\!\!\!\!\!\!\!\!\!\!\!\!\!\!\!\!\!\!\!\!\!\!\!\!\!\!\!\!\!&\!\!\!\!\!\!R_4\!\!\!&:=
\frac1y(1+(m\!+\!n)\a)H\UPo0K\UPo0\ptl_xQ_\nu-\frac1y(1+\a(1\!-\!\nu))
\big((m\!+\!n)H\UPo0\ptl_xK\UPo0+K\UPo0\ptl_xH\UPo0\big)Q_\nu
\nonumber\\&\!\!\!\!\!\!\!\!\!\!\!\!\!\!\!\!&\phantom{:=}+\frac{(m+n)\a }mc_\nu H\UPo0K\UPo0^{1-\nu}.
\label{R4====}
\end{eqnarray}\VS{-19pt}\PAR\NOindent
Multiplying the first equation by $-(1+(m+n)\a)$ then adding it to the second equation,
and using \eqref{H0-K0Ja}, we \VS{-15pt}solve
\begin{eqnarray}\label{Qnu===}&& Q_\nu=q_\nu+\b H\UPo0K\UPo0^{1-\nu},\mbox{ where }\b=-\frac{m + n}{m( m + n-1
+ \nu) }c_\nu,
\end{eqnarray}
and $q_\nu=
\frac{(1 + ( m + n)\a )R_3 - R_4}{\a( m + n-1
+ \nu)}$.
Thus if $c_\nu=0$, we obtain that $Q_\nu$ is {\it rational} (i.e., an element of the form $\frac{P}{Q}$ with $P,Q\in\PP$). Assume $c_\nu\ne0$ (thus $\nu\ge2$ by \eqref{F-now-is}). Similar to \eqref{HKUP-HK}, we obtain that (note from \eqref{nu====0} that the $\a$-type $\,-2\nu{\ssc\,}$-th components of $H\ptl_yK$ and  $H\ptl_xK$
only involve $H\UPo0,$ $K\UPo0,$ $H\UP{\nu},$ $K\UP{\nu},$ $H\UP{2\nu},$ $K\UP{2\nu}$\VS{-5pt})
\begin{eqnarray}
\label{R5==}
\!\!\!\!\!\!\!\!\!\!\!\!&R_5\!\!\!&:=(1-2\nu)H\UPo0K\UP{2\nu}+(1-\nu)H\UP{\nu}K\UP{\nu}+H\UP{2\nu}K\UPo0,\\
\label{R6==}
\!\!\!\!\!\!\!\!\!\!\!\!&R_6\!\!\!&:=H\UPo0\ptl_xK\UP{2\nu}+H\UP{\nu}\ptl_xK\UP{\nu}+H\UP{2\nu}\ptl_xK\UPo0,
\end{eqnarray}
are rational.
Our attempt was to use
\eqref{-ptl-x-K}--\eqref{Qnu===} to prove:\begin{itemize}\parskip-3pt\item[(i)]
 Suppose $\nu\ne\nu_0$. Then $Q_\nu$ is a {rational function}. Furthermore, $K\UPo0^\nu$ is rational if $c_\nu\ne0$; \item[(ii)]Suppose $\nu=\nu_0$ $($then $\frac1{\a}\in\N)$. Then $\ptl_yQ_{\nu_0}=-\frac{x}{\a y}\ptl_xQ_{\nu_0}$, and $H\UP{\nu_0}=\frac1{\a y}H\UPo0\ptl_xQ_{\nu_0}$,
$K\UP{\nu_0}=\frac1{y}K\UPo0\ptl_xQ_{\nu_0}$. Furthermore,
$\frac{m(1+(m+n)\a)}{(m+n)\a y}\ptl_xQ_{\nu_0}+c_{\nu_0} K\UPo0^{-\nu_0}$ is rational.
\end{itemize}
We claim that if (i) and (ii) hold, then it would imply that
a \QJ  pair $(F,G)$  in $\C[x^{\pm\frac1N},y]$ satisfying
\eqref{F-now-is}--\eqref{f0===} does not exist (which implies the two-dimensional Jacobi conjecture).
Although  (i) or (ii) might not necessarily be true, one may get some information from this.

\section{Weyl algebras}
\setcounter{theo}{0} \setcounter{equation}{0} \vs{0pt}\par In this
section, we first  generalize results of the previous sections. 
All undefined notations
can be found in the previous sections. The main results in this section are Theorems \ref{Dixmier-lemm}, \ref{a-b====}.

We denote
\begin{eqnarray}\label{A-ring-ABC1}
\!\!\!\!&\!\!\!\!&\mbox{$\AAu=\{f=\sum\limits_{i=0}^\infty
f_iu^{\a-\frac{i}{\b}}\,|\,f_i\in\C\mbox{ and some }\a,\b\in\Z,\,\b>0\}, $}
\\
\!\!\!\!&\!\!\!\!&\mbox{$
\BBu=
\{F=\sum\limits_{i=0}^\infty u^{\a-\frac i\b} F_i\,|\,F_i\in\C((v^{-1}))\mbox{ and some }\a,\b\in\Z,\,\b>0\}, $}
\end{eqnarray}
so that $\AAu$ is a field, and $\BBu$ is an associative unital algebra (which is in fact a divisible ring)
such that the product obeys the following law:
\equa{A-law}{v^iu^j=\dis\sum_{s\in\Z_+}s!\binom{i}{s}\binom{j}{s}u^{j-s}v^{i-s}\
\ \mbox{ for }\ \ i\in\Q,\,j\in\Z,} or more generally,
\equa{A-law-1}{v^if=\sum\limits_{s\in\Z_+}{\dis\binom{i}{s}}f_u^{(s)}v^{i-s}, \
gu^j=\sum\limits_{s\in\Z_+}{\dis\binom{j}{s}}u^{j-s}g^{(s)},
\mbox{ \ \ where $\dis
f_u^{(s)}=\ptl^s_uf,\ g^{(s)}=\ptl^s_vg$.}}
for $i\in\Z,\,j\in\Q,\,f\in\AAu, g\in\C((v^{-1}))$.
 Thus, the Weyl
algebra $W_1$ is the subalgebra of $\BBu$ generated by $u,v$. We remark that
an element $F$ of $\BBu$ is a combination of rational powers of $u$ with coefficients in $\C((v^{-1}))$, and we \textbf{always} write
an element $F$ 
in its {\it standard form}, namely, $u$ always appears before $v$ in any term of $F$.
\def\rmap#1{{\overrightarrow{#1}}}\def\lmap#1{{\overleftarrow{#1}}}Then we can define a linear map \equa{bar-map}{\mbox{$\rmap{\phantom{a}}:\BBu\to\BB$ such that $\rmap u=-x$ and $\rmap v=y$.}}
We also denote $\lmap{\phantom{a}}$ the inverse map of $\rmap{\phantom{a}}$. Then clearly, for any $F,G\in\BBu$, we \VS{-7pt}have
\equa{A-B-breack}{FG=\lmap{\rmap F\,\rmap G}+\cdots,\ \ \ [F,G]=\lmap{[\rmap F,\rmap G]}+\cdots,}
where the omitted terms in the first  (resp., second) equation
have $u$-degrees $<\deg_uFG$ (resp., $\deg_u[F,G]$), and
the bracket in the left-hand side is the usual commutator in $\BBu$ defined by \eqref{Brak},
the bracket in the right-hand side is the bracket in $\BB$ defined by the Jacobian determinant \eqref{Lie-b}.
If we define the 
bracket $[\cdot,\cdot]_{\ssc\cal W}$ in $\BB$ as
\equa{--Bra-W}{[F,G]_{\ssc\cal W}=\sum\limits_{i=1}^\infty{\dis\frac1{i!}}\Big((\ptl_{\barx}^i F)(\ptl_{\bary}^iG)-(\ptl_{\bary}^iF)( \ptl_{\barx}^iG)\Big)\mbox{ \ for }F,G\in\C[\barx,\bary],}
%
then clearly, the two Lie algebras $(\BBu,[\cdot,\cdot])$ and $(\BB,[\cdot,\cdot]_{\ssc\cal W})$ are isomorphic under the map
$\rmap{\phantom{a}}$ in \eqref{bar-map}.

Now consider the Weyl algebra $\BBu$.
Any element $F=\sum_{i=0}^\infty u^{\a-\frac{i}{\b}}f_i\in\BBu$ with
$f_0\ne0$ being monic
has the inverse $F^{-1}$, which is defined to be the unique element $H=\sum_{i=0}^\infty
u^{-\a-\frac{i}{\b}}h_i$
with $h_0=f_0^{-1}$ \VS{-7pt}and
\equa{qqq}{1=FH=\sum\limits_{i,j,s\in\Z_+}\dis \binom{-\a-\frac{i}{\b}}{s}u^{-\frac{i+j}{\b}-s}f^{(s)}_ih_j\mbox{ \ (cf.~\eqref{A-law-1}).}
}
Note that $h_i$ is uniquely determined for all $i$.
Further assume $f_0$ has  degree $\deg_vf_0=m>0$, then for any $a,b\in\Z,\,b>0$ with $b|am$, we
can define $F^{\frac{a}{b}}$ to be the unique element
$E=\sum_{i\ge0}u^{\frac{a\a}{b}-\frac{i}{\b}}e_i$ in $\BBu$ such that
$e_0=f_0^{\frac{a}{b}}$ (which is defined as in (\ref{h-a-power}))
 \VS{-7pt}and
\begin{eqnarray}
\label{A-equa2.2}
 \!\!\!\!\!\!\!\!\!\!&\!\!\!\!\!\!\!\!\!\!F^{a}
\!\!\!&=E^b=u^{a\a}e_0^{b}+\sum_{s=0}^{b-1}
(u^{\frac{a\a}{b}}e_0)^s(u^{\frac{a\a}{b}-\frac1{\b}}e_1)(u^{\frac{a\a}{b}}e_0)^{b-1-s}+\cdots.
\end{eqnarray}
Using (\ref{A-law-1}), we see that when writing the right-hand side
of (\ref{A-equa2.2}) as a standard form, the coefficient of
$u^{\frac{a\a}{b}-\frac{s}{\b}}$ is (where we use $*$ to denote some coefficients which can be determined but not needed for our purpose)\VS{-10pt},
\equa{A-law-2}{\ \ \ \ \ \ \ \dis\Coeff(E^b,u^{\frac{a\a}{b}-\frac{s}{\b}})=e_0^{b-1}e_s+
\sum_{\stackrel{\sc 0\le i_0\le i_1\le\cdots\le i_{s}<s}{\sc
i_0+\cdots+i_s+k_0+\cdots+k_s=s}}
*\,e_{i_0}^{(k_0)}e_{i_1}^{(k_1)}\cdots e_{i_{s}}^{(k_{s})},}
Thus for each $s\ge1$, (\ref{A-equa2.2}) has a unique solution for
$e_s$, which has the \VS{-7pt}form \equa{A-h-s-form}{e_s=e_0^k h=f^{\frac{ak}{b}}h\mbox{ \ \
for \ some \ \ }k\in\Z,\,h\in\C((v^{-1})).}
Note that if $F\ne0$, then we have $0=[F,FF^{-1}]=F[F,F^{-1}]$, which implies $[F,F^{-1}]=0$. Then for
any $a,b\in\Q$, we can write $a=\frac{c}{q},b=\frac{d}{q}\in\Q$ with $c,d,q\in\Z$, $q>0$, and if we write $H=F^{\frac1q}$, then
(say $c,d>0$\VS{-7pt})
\equa{AMSSSS}{[F^a,F^b]\!=\![H^{c},H^{d}]=\sum\limits_{i=0}^{c-1}H^i[H,H^{d}]H^{c-i-1}=
\sum\limits_{i=0}^{c-1}\sum\limits_{j=0}^{d-1}H^i(H^j[H,H]H^{d-j-1})H^{c-i-1}=0.}

 Now suppose  $(F,G)$ is a Dixmier pair in $W_1$,
i.e., $F,G\in W_1$ with $[F,G]=1$.
As before, we can express $G$ as
\begin{equation}
\label{A-equa2.5} G=\mbox{$\sum\limits_{i=0}^\infty$}
b_iF^{\frac{n-i}{m}} \mbox{ \ \ for some \ }b_i\in\AAu.
\end{equation}
\def\ign{\mbox{{\it $\lceil$ignored$\rceil$}}}We define
 $p=p(F)$ as before.
First we assume $p>-1$. Then from (\ref{A-law}), we can easily
observe \equa{A-commute}{HK=KH+\ign\mbox{ \ \ for \ any \ \
}H,K\in\BBu,} where, we  use \ign\ to denote terms whose
component lines are  located below the prime line of the proceeding
term, cf.~(\ref{p-line})). Similar as in Definition \ref{comp}, we
define $\pri{F}\in\C[u^{\pm1}][v]$ to be the $p$-type q.h.e.  such
that \equa{A-igno}{F\comp{0}=u^{m_0}\pri{F}^{m'}+\ign,}
 with $m'$
maximal.
 Then
from
$\sum_{i=0}^\infty[F,b_i]F^{\frac{n-i}{m}}=\sum_{i=0}^\infty[F,b_iF^{\frac{n-i}{m}}]=[F,G]\in\C^*$,
we see, \equa{A-wesee}{b_i\in\C\ (i\le m+n-2),\ \ \
\deg_ub_{m+n-1}=1+\si_0,\ \ \ \deg_ub_{m+n-1+i}\le1+\si_i\ (i\ge0),}
where $\si_i$ is defined in (\ref{equa2.7+-+0}), and the last
equation follows by comparing the $p$-type component. Using
(\ref{A-commute}), as in the proof of (\ref{A-h-s-form}), we see
that the analogous result of Lemma \ref{lemm-product-KH}(5) also
holds, i.e., for all $\ell\in\Z$ with $d|\ell$, where
$d=\deg_v\pri{F}$, and all $r\in\Q$, the element
$u^{-\frac{m_o\ell}{m}}(F^{\frac{\ell}{m}})\comp{r}$ of $\BBu$ is of
the form $\pri{F}^aP$ for some $a\in\Z,\,P\in\C[u^{\pm1}][y]$. Using
this and (\ref{A-commute}), we obtain $d|(n-i)$ if $i\le m+n-2$ with
$b_i\ne0$ (cf.~Lemma \ref{d|b}), and we have analogous results of
Lemmas \ref{quasi-J}, \ref{lemm-g-d2-i-KEY} and Theorem \ref{C[x,y]'}.
In particular, we can assume
\equa{p==Dix}{p\le 0,\ \ \ \mbox{ and furthermore,  $p\le -1$ if $m_0>m>0$}.}

 We remark that the following discussions will be similar to Subsection \ref{S-NewtonP}.
We define the Newton polygon of $F$ as in Subsection \ref{S-NewtonP}.
Let $(m_0,m)$ be any (not necessarily the top most) vertex of $\Supp\,F$ such that $m_0\ge0$, and $(n_0,n)$ is the corresponding vertex of $\Supp\,G$ (then $\frac{m_0}{m}=\frac{n_0}{n}$).
As before, we always assume $2\le m<n$ and $m\mbox{$\not|$}\,n$.
We can assume $m_0\le m$  (if necessary by using the automorphism $(u,v)\mapsto(v,-u)$, cf.~Remark  \ref{m0-m}).
Note that in our case here, we cannot assume $m_0\ne m$ (which will be clear later). Also
note that the arguments below do not need to assume $m_0>0$.
As  in Subsection \ref{S-NewtonP}, we regard $F,G$ as elements in
$\C[u^{\pm\frac1N},v]$ for some sufficient large $N$.
We define the {\it prime degree} $p$ as in  Theorem \ref{Nerw}, such that $-\frac1p$ is the slope of the unique edge (denoted by $L$)  of $\Supp\,F$ which is located at the right bottom side of $\Supp\,F$ with  top vertex $(m_0,m)$.
Note that
 if $F$ has the form
 $F=\sum_{i=0}^\infty u^{\a-\frac{i}{\b}}f_i\in\BBu$ with
$f_0\ne0$, we can rescal $F$ so that $f_0$ becomes a monic polynomial of $v$
(however we can only rescal $u$ by $au$, and in the meantime rescal $v$ by $a^{-1}v$ for some $0\ne a\in\C$ in order for $u,v$ to satisfy $[v,u]=1$).
First we need the following.
\begin{lemm}\label{p--large-than}~$p>-\frac{m_0}{m}.$\end{lemm}
\noindent{\it Proof.~}~Suppose conversely, $p\le -\frac{m_0}{m}$. First assume $m_0<m$. Define the automorphism $\si$ of  $\C[u^{\pm\frac1N},v]$ to \VS{-6pt}be \equan{u----->v}{\dis\si:(u,v)\mapsto
\Big(\frac{m-m_0}{m}u^{\frac{m}{m-m_0}},u^{\frac{m_0}{m-m_0}}v\Big).} Note that under $\si$, the edge $L$ is mapped to the $v$-axis, and that $F,G$ are mapped to some elements in $\C[u^{-\frac1N},v]$. However, any pair of elements in $\C[u^{-\frac1N},v]$ cannot form a \QJ  pair, a contradiction. Now assume $m_0=m$.
As in the proof of  Theorem \ref{Nerw}:
Let
$z\in\C\bs\{0\}$ be an indeterminate, and apply the automorphism $(u,v)\mapsto(zu,z^{-1}v)$.
As in \eqref{tilde-F-G},
we have (note that if we use notation $\frac{\dot p}{\dot q}$ as in \eqref{tilde-F-G}, then here we define
$-\frac{\dot p}{\dot q}$ to be $-1=-\frac{m_0}{m}$, but not to be $p$\VS{-7pt})
\begin{eqnarray*}
\!\!\!\!\!\!&\!\!\!\!\!\!\!\!\!\!\!\!\!\!\!\!\!\!\!\!&
\widetilde F=z^{ m- m_0}F(zx,z^{-1}y),\
\ \widetilde G=z^{ n- n_0}G(zx,z^{-1}y),\ \
 \widetilde J:=[\widetilde F,\widetilde G]=z^{\mu}J,
\end{eqnarray*}
where  $\mu=m+n-m_0-n_0>0$. Note that since $p<-1$,
both $\widetilde F$ and $\widetilde G$ only contain non-positive (rational) powers of $z$, thus the last equation cannot hold, a contradiction.
This proves the \linebreak lemma.\hfill$\Box$\vskip5pt
Now by Lemma \ref{p--large-than}, we have $p>-\frac{m_0}{m}\ge-1$, so we can always
apply the automorphism \VS{-7pt}(cf.~\eqref{x-yyyyy}, note that $p$ is the $-\frac{\Dhat p}{\Dhat q}$ there)\equa{MAMSMS}{(u,v)\mapsto ((1+p)u^{\frac1{1+p}},u^{\frac{p}{1+p}}v),} so that the edge $L$ becomes an edge
in the first quadrant, which is parallel to the $y$-axis. Thus we can always assume $0<m_0\le m$ and $p=0$.
Hence we can write $F$ as $F=\sum_{i=0}^{M_1}u^{m_0-\frac{i}{N}}f_i$ for some $f_i\in\C[v]$ such that
$f_0
$ is a  polynomial  (which can be assumed to be monic) of $v$ of degree $m>0$, which contains at least two terms. We can further suppose that $f_0$ has at least two different roots, otherwise by change $v$ to $v+\a$ (where $\a$ is the root of $f_0$),
$f_0$ becomes  $v^m$, i.e., $p$ becomes negative (and we repeat the above to use \eqref{MAMSMS} to change $p$ to zero, this repeating can only last finite times, as in the proof of  Theorem \ref{Nerw}, cf.~statement after \eqref{OTHWWW}). In particular, we can obtain $m_0\ne m$. Analogously, we write $G=\sum_{i=0}^{M_2}u^{n_0-\frac iN}g_i$ with $g_0$ being a polynomial of $v$-degree $n$.
We always regard $F,G$ as in $\BBu$.
Note from $[F,G]=J\ne0$ that $m_0+n_0\ge1$.

First assume $m_0+n_0>1$. Using \eqref{A-B-breack} and Comparing the coefficients of $u^{m_0+n_0-1}$ in $[F,G]=J$,
we \VS{-7pt}obtain \equa{g0==f00000}{g_0=b_0f_0^{\frac{n}{m}} \mbox{ \ for some nonzero $b_0\in\C$.}}
Write $f={\pri{F}}^{m'}$ with $m'$ maximal. Then \eqref{g0==f00000} proves $\pri{F}^{\frac{nm'}{m}}$ is a polynomial.
%
Denote $G^1=G-\sum_{i=0}^\infty b_{i0}F^{\frac{n-i}{m}}$. Discussing as above (with $G$ replaced by $G_1$ and using \eqref{bar-map}) and continuing (similar to the arguments before \eqref{widetilde-G==}),
 we can eventually write $G$ \VS{-7pt}as
%
\equa{Now-G==F-ab}{G=\sum\limits_{s\in\frac1N\Z_+,\,s<m_0-1}b_{s}F^{\frac{k_s}{m}}+u^{1-m_0}R,}
for some $k_s\in\Z$, $ b_{s}\in\C$ with $k_0=n$, and some $R\in\BBu$ with $\deg_uR\le0$, so that $R$ can be written as $R_0=\sum_{i=0}^\infty u^{-\frac{i}{N}}r_i$.
Using \eqref{A-equa2.2},  \eqref{A-h-s-form} and  \eqref{g0==f00000}, as in the proof of Lemma \ref{bi-lemm},
we see from \eqref{Now-G==F-ab} that  $F^{\frac{k_s}{m}}$ is rational if $b_s\ne0$ (here an element $H=\sum_{i=0}^\infty u^{a-\frac{i}{N}}h_i$ is {\it rational} if each $h_i$ is a rational function of $v$), and $r_0$ is a rational function of the form $r_0=\pri F^{-m'-a}P$ for some $a\in\Z,\,P\in\C[v]$. Using
\eqref{bar-map}, we have $J=[F,G]=[F,u^{1-m_0}R]$, which implies that $(x^{m_0}\lmap f,x^{1-m_0}{\lmap {r_0}})$ is a \QJ  pair in $\BB$
by comparing the coefficients of $u^0$ and using \eqref{AMSSSS}. Namely, we have (where the prime stands for $\ptl_v$)
\equa{mmmm0000}{-(m'+am_0)\pri F^{-1-a}P\pri F'+m_0\pri F^{-a}P'=m_0fr'_0-(1-m_0)f'r_0=J\in\C,}
which has exactly the same form of \eqref{equa-int1} (with $p'=0,q=1$). Thus as the discussions after
\eqref{f0-p==0}, we can find a lower vertex of $\Supp\,F$. Continuing the above process (from \eqref{MAMSMS}), $F,G$ can finally become elements such that $m_0+n_0=1$, so we can write $m_0,n_0$ as $m_0=\frac{m}{m+n}$, $n_0=\frac{n}{m+n}$. Thus in fact we have proved the following
\begin{theo}\label{Dixmier-lemm}
Suppose there exists a Dixmier pair $(F,G)$ in $\C[u,v]$ such that the Newton polygon of $F$ has
a vertex $(m_0,m)$ with $m_0, m>0$ $($we can assume $ m_0\le m$ by using the automorphism $(u,v)\mapsto(v,-u)$ if necessary, cf.~Remark  $\ref{m0-m})$,
and $(n_0,n)$ is the corresponding vertex of $\Supp\,G$ satisfying $\frac{m_0}{ m}=\frac{n_0}{n}$ and $2\le m< n$ and $ m\mbox{$\not|\,$} n$.
Then there exists an automorphism $\si$ of $\C[u^{\pm\frac1N},v]$, such that the pair
$(\si(\ol F),\si(\ol G))$, again denoted as $(F,G)$, having the form \eqref{F-now-is}--\eqref{f0===} $($with $x,y$ replaced by $-u,v)$ such that $2\le m<n$ and $m\mbox{$\not|$}\,n$.
\end{theo}


{}From now on, we assume that
$(F,G)$ is a Dixmier pair in $\C[u,v]$ such that the Newton polygon of $F$ has
the a vertex $(m_0,m)$ with $m_0>m>0$.

Regarding the edge at the right bottom side of $\Supp\,F$ with top vertex $(m_0,m)$ as the prime line,
we can expression $G$ as in
 (\ref{A-equa2.5}). We rewrite it as
\begin{equation}
\label{A-equa2.5+} G=\mbox{$\sum\limits_{i=0}^{m+n-1}$}
b'_iF^{\frac{n-i}{m}}+R,
\end{equation} where $b'_i\in\C$ is defined similarly as in
(\ref{denote-b'-i}). Then all elements $G,F^{\frac{n-i}{m}},R$ are
in $\C[u^{\pm1}]((v^{-1}))$. We denote \equa{ww===}{w=uv.} It is easy to verify
that \equa{A-uw}{w^i u^j=u^j(w+j)^i\mbox{ \ \  for \ \ }i,j\in\Z.}
Note that for $i,j\in\Z$,
\equa{A-uv+}{u^iv^j=\left\{\begin{array}{ll}u^{i-j}w(w-1)\cdots(w-j+1)&\mbox{if \ }j\ge0,\\[4pt]
u^{i-j}\big((w+j')(w+j'-1)\cdots(w+1)\big)^{-1}&\mbox{if \
}j=-j'<0.\end{array}\right.} Thus we can regard the above elements
as in $\C[u^{\pm1}]((w^{-1}))$. Then
\begin{equation}
\label{A-equa2.5+1} G\ptl_w F=\mbox{$\dis\sum\limits_{i=0}^{m+n-1}$}
b'_iF^{\frac{n-i}{m}}\ptl_w F+R\ptl_w F.
\end{equation}
For $F\in\C[u^{\pm1}]((w^{-1}))$, as in \eqref{Res-x-y}, we define the {\it trace} of $F$
to be \equa{tr-F-def}{{\rm
tr}(F)=\Coeff(F,u^0w^{-1})=\Coeff(F,u^{-1}v^{-1}),} where the second
equality follows from (\ref{A-uv+}) by regarding $F$ as an element
in $\C[u^{\pm1}]((v^{-1}))$. Using (\ref{A-uw}), we see
\equa{A-tr}{{\rm tr}(HK)={\rm tr}(KH)\mbox{ \ \ for \ \
}H,K\in\C[u^{\pm1}]((w^{-1})).} Because of this property, we call it ``trace''.
Now we shall compute the trace of
(\ref{A-equa2.5+1}). Obviously ${\rm tr}(G\ptl_w F)=0$
since $F,G\in\C[u^{\pm1}][w]$. Let $H=F^{\frac{1}{m}}$, then $${\rm
tr}(F^{\frac{n-i}{m}}\ptl_w F)={\rm
tr}(H^{n-i}\ptl_w(H^m))=\frac{m}{m+n-i}{\rm
tr}(\ptl_w( H^{m+n-i}))=0,$$ by noting that
$\ptl_w$ is a derivation of $\C[u^{\pm1}]((w^{-1}))$
(using (\ref{A-uw}) to verify) and for $i\le0$, we have
$\ptl_w H^i
=\sum_{s=0}^{i-1}H^s\ptl_w (H) H^{i-1-s}$, and
$\ptl_w(H^{-1})=-H^{-1}\ptl_w( H)H^{-1}$ by
using $HH^{-1}=1$, and using (\ref{A-tr}). This proves
\equa{Tr-A}{\dis{\rm tr}(R\ptl_w F)=\tr(G\ptl_w F)=0.}
Since we can assume $p\le-1$ by \eqref{p==Dix}, if we take $p'=-1\ge p$, then we can consider the
$p'$-type components. So up to a nonzero scalar, we can assume
\equa{eeF-R=====}{\begin{array}{ll}F\comp{0}=u^{m_0}v^m+\a u^{m_0-1}v^{m-1}+\mbox{(lower
terms)},\\[6pt]R\comp0=r_0(u^{1-m_0}v^{1-m}+\b u^{-m_0}v^{-m})+\mbox{(lower
terms)},\end{array}}
for
some $\a,\b\in\C$ (if $p<p'$ then $F\comp{0}=u^{m_0}v^m$ and so
$\a=0$) and $r_0=\frac1{m-m_0}$ (by condition $[F\comp0,R\comp0]=1$).
Thus
(using \eqref{A-uw}){
\begin{eqnarray}\label{F0-1}\!\!\!\!\!\!\!\!\!\!\!\!\!\!\!\!\!\!&\!\!\!\!\!\!\!\!\!\!\!\!\!\!\!\!&
F\comp{0}=u^{m_0-m}w(w-1)\cdots(w-m+1+\a)+...
=u^{m_0-m}\Big(w^m-\Big(\binom{m}{2}-\a\Big)w^{m-1}\Big)+...
,\\[7pt]\label{F1-1}
\!\!\!\!\!\!\!\!\!\!\!\!\!\!\!\!\!\!&\!\!\!\!\!\!\!\!\!\!\!\!\!\!\!\!&R\comp{0}=
r_0u^{m-m_0}\Big(((w+m-1)\cdots(w+1))^{-1}+\b((w+m)\cdots(w+1))^{-1}\Big)+...
\nonumber\\
\!\!\!\!\!\!\!\!\!\!\!\!&\!\!\!\!\!\!\!\!\!\!\!\!\!\!\!\!& \phantom{R\comp{0}}=
r_0u^{m-m_0}\Big(w^{1-m}-\Big(\binom{m}{2}-\b\Big)w^{-m}\Big)+...,
\end{eqnarray}}%
where the omitted terms do not contribute to the trace. Hence
(using \eqref{A-uw} and \eqref{A-tr}){
$$\begin{array}{lll}\dis
(m-m_0)R\ptl_w F&\!\!\!=\!\!\!&\dis
m\Big((w+m_0-m)^{1-m}\dis-\Big(\binom{m}{2}-\b\Big)(w+m_0-m)^{-m}\Big)w^{m-1}\\[6pt]
&\!\!\!\!\!\!&\dis-(m-1)\Big(
\binom{m}{2}-\a\Big)w^{-1}+...\\[9pt]
&\!\!\!=\!\!\!&\dis-\Big(\mbox{$\dis (2m_0-1)\binom{m}{2}$}-m\b-(m-1)\a\Big)w^{-1}+...\end{array}
$$}%
Note that for $r<0$, $(R\ptl_w F)\comp{r}$ does not
contain the term $u^0w^{-1}$, thus \equa{AMAMAM}{\dis 0=(m_0\!-\!m){\rm
tr}(R\ptl_w F)=(m_0\!-\!m){\rm tr}(R\comp{0}\ptl_w
F\comp{0})=(2m_0-1)\binom{m}{2}-m\b-(m-1)\a.} Since we
can also regard elements in (\ref{A-equa2.5+}) as in
 $\C[w^{\pm1}]((v^{-1}))$ (to distinguish the difference, we use $\ptl^v_w$ to denote the derivative $\ptl_w$ in $\C[w^{\pm1}]((v^{-1}))$). Computing as above, we obtain
\equa{AMAMAM1}{\dis 0=(m_0-m)\tr(F\ptl_w^vR)=(2m-1)\binom{m_0}{2}-m_0\b-(m_0-1)\a.}
This together with (\ref{AMAMAM}) proves Theorem \ref{a-b====}(1) below.

\begin{theo}\label{a-b====}\begin{itemize}\parskip-3pt\item[\rm(1)]
Assume $(F,G)$ is a Dixmier pair such that $\Supp\,F$ has a vertex $(m_0,m)$ with $m_0>m>0$. Then the edge $L$
at the right bottom side of $\Supp\,F$ with top vertex $(m_0,m)$ always has slope $1$. Furthermore, if we denote $F\comp0$ to be the part of $F$ corresponding to $L$, then $F\comp0=u^{m_0}v^m+\frac{m_0m}{2}u^{m_0-1}v^{m-1}+\cdots$. More precisely, in \eqref{eeF-R=====}, $\a,\b$ are equal to
\equa{Al-Be}{\dis\a=\frac{m_0 m}{2}, \ \ \ \ \b=\frac{(1 - m_0)(1 -
m)}{2}.} 
\item[\rm(2)]If we write $F\comp0,\,R\comp0$ as $F\comp0=u^{m_0-m}f(w),\,R\comp0=r(w)u^{m-m_0}$. Then \equa{f-r-w==}{\dis f(w)r(w)=\frac{1}{m-m_0}(w+\frac{m_0-m+1}{2}).}
\item[\rm(3)]Assume $[F,G]=1$. An analogous result to Theorem $\ref{Jacobian-el}(1)$ holds, namely, for any automorphism $\si$ of $\BBu$, we have \equa{MMMMSMS}{\tr(\si(G)\ptl_w\si(F))=\tr(\si(v)\ptl_w\si(u)).}
\end{itemize}\end{theo}

\begin{rema}\rm\label{deri-w}\begin{itemize}\parskip-3pt\item[(1)]We remark that Lemma \ref{a-b====} is the place where the great difference between Newton polygons of Jacobi pairs and Dixmier pairs occurs; for the Jacobi pairs, an edge of the Newton polygon can never have slope $1$ (cf.~Theorem \ref{C[x,y]'}).\item[(2)]
The operators $\ptl_w,\ptl^v_w$ are in fact the unique
derivatives in $\BBu$ such that
\equa{dermmm}{
\ptl_w(u)=0,\ \ptl_w(v)=u^{-1},\ \ \ \ptl^v_w(u)=v^{-1},\ \ptl^v_w(v)=0.}
To compute $\ptl^v_w(u^{\frac{p}{q}})$ with  $\frac{p}{q}\notin\Z$, $p,q\!>\!0$, we  set $h=u^{\frac{p}{q}}$ and
assume
$\ptl^v_w(h)=\sum_{i=1}^\infty c_iu^{\frac{p}{q}-i}v^{-i}$ for some $c_i\in\C,$
 and
 use $\ptl_w (u^p)=\ptl_w (h^q)=\sum_{i=0}^{q-1}h^i(\ptl_w h)h^{q-i-1}$ to determine $c_i$. Similarly, we can use $hh^{-1}=1$ to determine $\ptl^v_w(h^{-1})$. One can obtain
 \equa{eew---}{\dis\ptl^v_w(u^a)=au^{a-1}v^{-1}-\binom{a}{2}u^{a-2}v^{-2}+\cdots\mbox{ \ for any \ }a\in\Q.}
\end{itemize}\end{rema}
\noindent{\it Proof of Theorem \ref{a-b====}.~}~(1) has been proved. To prove (2), using $$1\!=\![F\comp0,R\comp0]\!=\!u^{m_0-m}f(w)r(w)u^{m-m_0}\!-\!r(w)f(w)\!=\!f(w\!+\!m\!-\!m_0)r(w\!+\!m\!-\!m_0)\!-\!r(w)f(w),$$
from this and \eqref{Al-Be}, we obtain \eqref{f-r-w==}.

(3) Denote $\bar u=\si(u),\bar w=\si(w)$ (then $\bar w\bar u=\bar u(\bar w+1)$), $a=m_0-m$,
and $'=\ptl_w$.
Using \eqref{Al-Be} and \eqref{f-r-w==}, one can see $f'(\bar w)r(\bar w)=0$.
Then (``\,$\equiv$\,'' means equality under taking $\tr$)
$$\begin{array}{lll}\si(G)\ptl_w\si(F)\!\!\!&\equiv
\ptl_w(\si(F\comp0))\si(R\comp0)\equiv\Big((\bar u^a)'f(\bar w)+\bar u^af'(\bar w)\Big)r(\bar w)\bar u^{-a}\\[7pt]
&\equiv\sum\limits_{i=0}^{a-1}\bar u^{a-i-1}\bar u'\bar u^i f(\bar w)r(\bar w)\bar u^{-a}+f'(\bar w)r(\bar w)\\[11pt]
&
\equiv\sum\limits_{i=0}^{a-1}\bar u'\bar u^i f(\bar w)r(\bar w)\bar u^{-i-1}
\equiv\bar u'\bar u^{-1}\sum\limits_{i=0}^{a-1}f(\bar w\!-\!i\!-\!1)r(\bar w\!-\!i\!-\!1)\\[11pt]
&\equiv\bar u'\bar u^{-1}\bar w\equiv\bar u'\bar v\equiv\si(v)\ptl_w\si(u).
\end{array}$$

{{\small\footnotesize\lineskip=1pt
\section*{\normalsize Acknowledgement\VS{-7pt}s}\label{ACKN}
The author would like
to thank Professor Zhexian Wan for encouragement;
Professors Arno van den Essen, Leonid Makar-Limanov, T.T.~Moh, Jingen Yang, Jietai Yu, and Dr.~Victor Zurkowski
for comments and suggestions
\vs{-7pt}.
%
%
}

 \small\footnotesize \lineskip=2pt

}
\end{document}